\documentclass[notitlepage,leqno,10pt]{article}

\usepackage{latexsym}
\usepackage{amsmath}
\usepackage{amsfonts}
\usepackage{amssymb}
\usepackage{psfig}
\usepackage[dvips]{graphicx}
\renewcommand{\a }{\alpha }
\renewcommand{\b }{\beta }
\renewcommand{\d}{\delta }
\newcommand{\D }{\Delta }

\newcommand{\e }{\varepsilon }
\newcommand{\g }{\gamma}

\newcommand{\G }{\Gamma }
\newcommand{\IO }{\int_\O }
\renewcommand{\l }{\lambda }

\newcommand{\n }{\nabla }
\newcommand{\var }{\varphi }
\newcommand{\rh }{\rho }

\newcommand{\s }{\sigma }
\newcommand{\Sig }{\Sigma}
\renewcommand{\t }{\tau }
\renewcommand{\th }{\theta }
\renewcommand{\o }{\omega }
\renewcommand{\O }{\Omega }

\newcommand{\ov}{\overline}

\newcommand{\be}{\begin{equation}}
\newcommand{\ee}{\end{equation}}
\newenvironment{pf}{\noindent{\sc Proof}.\enspace}{\rule{2mm}{2mm}\medskip}
\newenvironment{pfn}{\noindent{\sc Proof}}{\rule{2mm}{2mm}\medskip}

\newcommand{\R}{\mathbb{R}}
\newcommand{\C}{\mathbb{C}}

\newcommand{\N}{\mathbb{N}}

\author{Sagun CHANILLO$^{\rm a}$ \and Andrea MALCHIODI$^{\rm b}$}

\date{}

\title{Asymptotic Morse theory for the equation $\D v = 2 v_x \wedge v_y$}

\begin{document}

\newtheorem{lem}{Lemma}[section]
\newtheorem{pro}[lem]{Proposition}
\newtheorem{thm}[lem]{Theorem}
\newtheorem{rem}[lem]{Remark}
\newtheorem{cor}[lem]{Corollary}
\newtheorem{df}[lem]{Definition}

\maketitle

\begin{center}

$^{\rm a}${Rutgers University, Department of Mathematics \\ 110
Frelinghuysen Road, Piscataway, NJ 08854 USA}

$^{\rm b}${\small School of Mathematics Institute for Advanced
Study \\ 1 Einstein Drive, Princeton, NJ 08540 USA}

\end{center}

\footnotetext[1]{E-mail addresses: chanillo@math.rutgers.edu (S.
Chanillo), malchiod@ias.edu (A. Malchiodi)}

\noindent {\sc abstract}. - Given a smooth bounded domain $\O
\subseteq \R^2$, we consider the equation $\D v = 2 v_x \wedge
v_y$ in $\O$, where $v: \O \to \R^3$. We prescribe Dirichlet
boundary datum, and consider the case in which this datum
converges to zero. An asymptotic study of the corresponding Euler
functional is performed, analyzing multiple-bubbling phenomena.
This allows us to settle a particular case of a question raised by
H. Brezis and J.M. Coron in \cite{brc2}.

\

\begin{center}

\end{center}

\section{Introduction}\label{s:i}

Let $\O \subseteq \R^2$ be a smooth bounded domain. We shall
denote by $v, \tilde{g}$ two maps such that $v : \O \to \R^3$ and
$\tilde{g} : \partial \O \to \R^3$, with $\tilde{g}$ smooth.
Consider the problem
\begin{equation}\label{eq:in1}
  \begin{cases}
    \D v = H(\xi,v,\n v) v_x \wedge v_y & \text{ in } \O, \\
    v = \tilde{g}  & \text{ on } \partial \O,
  \end{cases}
\end{equation}
where $H$ is a smooth scalar function, $v_x$, $v_y$ are the $x$
and $y$-derivatives of $v$, $\xi = (x,y)$ and $\wedge$ denotes the
cross-product in $\R^3$.

Equation \eqref{eq:in1} has been the subject of several works, see
for example the survey paper \cite{st} by K. Steffen and the
recent paper \cite{cm}. Existence of solutions of \eqref{eq:in1}
when $\tilde{g} \equiv 0$ strongly depends on the topology of the
domain. In fact we show using a Pohozahev-type identity, see
Proposition \ref{p:poh}, that equation \eqref{eq:in1} has no
solution in any simply connected domain when $\tilde{g} = 0$. When
$H(\xi,v,\n v) \equiv H$, a non-zero constant, such a result was
proved by H. Wente, \cite{hw}, using reflection techniques and the
Kelvin transformation. In the same paper, Wente also showed that
if $\O$ is an annulus then the study of \eqref{eq:in1} can be
reduced to an ordinary differential equation and \eqref{eq:in1}
does have a non-trivial solution when $v = 0$ on $\partial \O$.
Thus equation \eqref{eq:in1} presents features similar to the
Yamabe equation on domains with Dirichlet boundary condition,
studied in particular by A. Bahri and J.M. Coron, \cite{bc}. In
fact, part of the difficulty in studying \eqref{eq:in1} is that it
is invariant under conformal transformations. This invariance
forces the associated variational problem to exhibit
non-compactness phenomena, like in the Yamabe problem on domains.
We point out that in our case, contrary to the Yamabe problem,
simply connected domains always admit only trivial solutions. For
the Yamabe problem in dimension greater or equal than three, there
are indeed examples of contractible domains which admit
non-trivial positive solutions, see \cite{pas}.

From now on we consider the case of constant $H$, precisely
$H(\xi,v,\n v) \equiv 2$. So problem \eqref{eq:in1} reduces to
\begin{equation}\label{eq:in2}
  \begin{cases}
    \D v = 2 v_x \wedge v_y,  & \text{ in } \O, \\
    v = \tilde{g} & \text{ on } \partial \O.
  \end{cases}
\end{equation}
Under the assumption $\|\tilde{g}\|_\infty < 1$, S. Hildebrandt,
\cite{hil}, constructed a solution of \eqref{eq:in2} with minimal
energy called the {\em small solution}, while Brezis and Coron,
\cite{brc}, K.Steffen, \cite{st0} and M. Struwe, \cite{str},
constructed a second solution, referred to as the {\em large
solution}. We remark that the assumption $\|\tilde{g}\|_\infty <
1$ is sharp, see \cite{hei}.

Results similar to those regarding the Dirichlet problem hold for
the {\em Plateau problem}, in which one looks for solutions of $\D
u = H u_x \wedge u_y$ which are conformal and which map the
boundary to a given curve (with free parametrization). As a result
one obtains surfaces with constant mean curvature.

We mainly focus on the following problem
\begin{equation}\label{eq:inge}
  \begin{cases}
    \D v = 2 v_x \wedge v_y,  & \text{ in } \O, \\
    u = \e \tilde{g} & \text{ on } \partial \O.
  \end{cases}
\end{equation}
We will study \eqref{eq:inge} turning it into a variational
problem. In view of the non-existence result in \cite{hei}, it is
natural to assume that the boundary datum is small. T. Isobe in
particular, \cite{isc}-\cite{isa2}, analyzed the behavior of the
large solutions of Brezis and Coron in the limit $\e \to 0$ (the
small solutions converges to the trivial one $v \equiv 0$ as $\e
\to 0$).

Let $g$ denote the harmonic extension of $\tilde{g}$ in $\O$, i.e.
\begin{equation}\label{eq:ogg}
  \begin{cases}
    \D g = 0 & \text{ in } \O;  \\
    g = \tilde{g} & \text{ on } \partial \O.
  \end{cases}
\end{equation}
If $v$ is a solution of \eqref{eq:inge} and if we set $v = u + \e
\, g$, the function $u$ solves
\begin{equation}\label{eq:pe}\tag{$P_\e$}
  \begin{cases}
    \D u = \D v = 2 \, (u_x + \e \, g_x) \wedge
    (u_y + \e \, g_y) & \text{ in } \O;  \\
    u = 0 & \text{ on } \partial \O.
  \end{cases}
\end{equation}
Problem \eqref{eq:pe} admits the Euler functional $I_\e : H^1_0(\O
; \R^3) \to \R$, which has the following expression
\begin{equation}\label{eq:ie}
  I_\e(u) = \frac{1}{2} \int_\O |\n u|^2 + \frac{2}{3} \int_\O u
  \cdot (u_x \wedge u_y) + \e \int_{\O} u \cdot (u_x \wedge g_y +
  g_x \wedge u_y) + 2 \e^2 \int_\O u \cdot (g_x \wedge g_y).
\end{equation}
The aim of this paper is to develop a Morse theory for the
functional $I_\e$ when $\e$ is small. In order to do this we take
advantage of the perturbative approach in \cite{ab}. We first
recall from \cite{brc2} that the {\em fundamental solution} ({\em
bubble}) of the equation
\begin{equation}\label{eq:ue}
  \D u = 2 u_x \wedge u_y,  \qquad \qquad \text{ in } \R^2
\end{equation}
is the stereographic projection $\pi : \R^2 \to S^2 \subseteq
\R^3$
\begin{equation}\label{eq:stpr}
  \pi(x,y) = \left( \frac{2x}{1 + x^2 + y^2},
  \frac{2y}{1 + x^2 + y^2}, \frac{x^2 + y^2 - 1}{1 + x^2 + y^2}
  \right), \qquad (x,y) \in \R^2.
\end{equation}
Our analysis will use translations, dilations and rotations of the
function in \eqref{eq:stpr} and we set
\begin{equation}\label{eq:in4}
  R \d_{a,\l} (x,y) \equiv  R \circ \pi (\l (x - a_1, y - a_2) ),
\end{equation}
for $R \in SO(3)$, $a = (a_1,a_2) \in \R^2$ and $\l > 0$. The
functions $R \d_{a,\l}$ are mountain-pass critical points of the
functional
\begin{equation}\label{eq:ovin}
  \ov{I}(u) = \frac{1}{2} \int_{\R^2} |\n u|^2 + \frac{2}{3}
  \int_{\R^2} u \cdot (u_x \wedge u_y), \qquad u \in \mathcal{D},
\end{equation}
where $\mathcal{D}$ denotes the functional space
$$
\mathcal{D} = \left\{ u \in L^2_{loc}(\R^2;\R^3)  \, : \,
\|u\|_{\mathcal{D}}^2 = \int_{\R^2} |\n u|^2 + \int_{\R^2}
\frac{|u|^2}{(1+|\xi|^2)^2} < + \infty \right\}.
$$
The space $\mathcal{D}$ coincides with $H^1(S^2;\R^3)$ after
inverse stereographic projection. We point out that the
functionals $\ov{I}$ and $I_\e$ are well defined and smooth on
$\mathcal{D}$ and $H^1_0(\O,\R^3)$ respectively, see Section
\ref{s:pre}. It turns out that the manifold constituted by the
$\d$'s is non-degenerate for the functional $\ov{I}$ (modulo
constants), as proved in \cite{isa1} Lemma 5.5, using an
isoperimetric inequality. Proving the non-degeneracy condition is
equivalent to classify the solutions of
\begin{equation}\label{eq:in5}
  \D w = 2 \left( w_x \wedge \d_y + \d_x \wedge w_y \right) \qquad
  \hbox{ in } \R^2,
\end{equation}
which is the linearization of \eqref{eq:ue} around $\d$, and to
show that the only solutions are the tangent vectors to $\ov{Z}$
at $\d$, see equation \eqref{eq:ovz}. We remark that equation
\eqref{eq:ue} admits solutions of the form $\pi(z^k)$ (in complex
notation) for any integer $k$, see \cite{brc2}. We will refer to
these solutions as {\em higher degree bubbles}. For this reason we
give in the Appendix an alternative proof of the non-degeneracy,
which we believe could adapt naturally to the higher-degree case.

To analyze the problem in $\O$, it is convenient to consider the
functions $P \d = \d - \var$, where $\var$ is defined in
\eqref{eq:fi}. $P \d$ is the element of $H^1_0(\O;\R^3)$ closest
to $\d$ in the Dirichlet norm. We may write
\begin{equation}\label{eq:in6}
  u = \sum_{i=1}^k P R_i \d_{p_i,\l_i} + w,
\end{equation}
where $R_i \in O(3)$, $\l_i > 0$, $p_i \in \R^2$ for all $i$, and
$w$ is orthogonal to the manifold $\sum_{i=1}^k P R_i
\d_{p_i,\l_i}$. Once we have the non-degeneracy property for
$\ov{I}$, then it is standard to prove that for suitable values of
$a$ and $\l$ also the manifold of projected bubbles is
non-degenerate for $I_\e$, and the same holds true for a finite
sum of bubbles. This property allows us to solve the equation
$I'_\e(u) = 0$ in $w$ (see Proposition \ref{p:ambbad}), and thus
our problem is reduced to a finite-dimensional one which involves
an auxiliary functional $\tilde{I}_\e(z)$ (see Section \ref{s:4})
depending only on $\{p_i\}_i, \{l_i\}_i$ and $\{R_i\}_i$.
Substituting \eqref{eq:in6} into $I_{\e}$ and letting $\e \to 0$,
we expand $\tilde{I}_\e(z)$ for large values of $\l_i$ (roughly of
order $\e^{-1}$).

The large solution of Brezis and Coron corresponds to a {\em one
bubble solution} when $\e \to 0$, and has been studied in detail
by T. Isobe, \cite{isc}-\cite{isa2}. However, from Theorem 0.3 in
\cite{brc2} it is clear that a more complicated configuration may
occur. Thus to manufacture this type of solutions we are naturally
led to a variational analysis of the functional \eqref{eq:ie} for
multiple bubbles. We point out that from the work of Brezis-Coron
the bubbles will not necessarily be all of degree 1. However the
variational analysis is more difficult if we allow bubbles of
arbitrary degree, and we will return to this point in a subsequent
article.

To state our results we need some notation. Given $(a, \xi) =
((a_1,a_2), (x,y)) \in \O \times \O$, let $h_1, h_2 : \O \times \O
\to \R$ be the solutions of the problems
\begin{equation}\label{eq:h1h2}
  \begin{cases}
    \D_\xi h_1 (a,\xi) = 0 & \text{ in } \O, \\
    h_1(a,\xi) = \frac{\xi_1-a_1}{|\xi-a|^2} & \text{ on } \partial
    \O;
  \end{cases}
  \qquad
\begin{cases}
    \D_\xi h_2 (a,\xi) = 0& \text{ in } \O, \\
    h_2(a,\xi) = \frac{\xi_2-a_2}{|\xi-a|^2} & \text{ on } \partial
    \O,
  \end{cases}
\end{equation}
see Remark \ref{r:pos} $(b)$. If $G(a,\xi)$ denotes the Green's
function of $\O$, normalized so that $G (a,\xi) \sim - \log
|a-\xi|$ for $a \sim \xi$, and if $H(a,\xi)$ denotes the regular
part of $G$ ($G(a, \xi) = - \log |a - \xi| - H(a,\xi)$), then we
have
\begin{equation}\label{eq:hhi}
 h_1(a,\xi) = \frac{\partial H(a,\xi)}{\partial a_1}; \qquad
 \qquad
 h_2(a,\xi) = \frac{\partial H(a,\xi)}{\partial a_2}.
\end{equation}
Let also
\begin{equation}\label{eq:hhi2}
  \tilde{H}(a) = \left( \frac{\partial h_1}{\partial x} +
\frac{\partial h_2}{\partial y} \right)|_{\xi = a}.
\end{equation}
It has been proved in \cite{isc}-\cite{isa2} (see also \cite{s})
that the function $\tilde{H}$ plays a crucial role in studying the
location of blowing-up solutions of \eqref{eq:in2}, when the
boundary datum converges to zero. In fact, the function
$\tilde{H}$ appears in the expansion of $I_\e(u)$, when $u$ is of
the form \eqref{eq:in6} with $k = 1$, as a self-interaction term,
see Proposition \ref{p:}. The expansion for $k = 1$ is essentially
performed in the works of Isobe, but we derive it in a framework
which is convenient to treat the case of $k > 1$, see Section
\ref{s:2b}. We have the following result, regarding the function
$\tilde{H}$.

\begin{thm}\label{t:in1}
(a) For a simply connected domain $\O$ there holds
$$
\tilde{H} (a) = 2 e^{2 H(a,a)}, \qquad a \in \O,
$$
where $H(a,\xi)$ is the regular part of the Green's function of
$\O$.

(b) For a multiply connected domain $\O$ there not exist in
general a function $F$ such that $\tilde{H}(a) = F(H(a,a))$. In
particular, for some annulus of the form $\left\{ \rh^{-1} < |z| <
\rh \right\}$, $\rh
> 1$, the critical points of $H(a,a)$ and of $\tilde{H}(a)$ do not
coincide.
\end{thm}

\noindent Theorem \ref{t:in1} is proved in Section \ref{s:robin}.
The function $H(a,a)$ is called the Robin function of the domain
$\O$, see \cite{baf}. In dimension $2$ it also appears in extremal
problems related to the Moser-Trudinger inequality, where the
critical points of $H(a,a)$ are shown to be related to the
conformal incenter of $\O$. Isobe showed that $\tilde{H} > 0$ on
any domain, see Remark \ref{r:pos}, but did not analyze it
further. Since $\tilde{H}$ is defined by means of second
derivatives of $H$, we need to use a global argument (the Riemann
mapping theorem) to compare the two functions $H$ and $\tilde{H}$.
The new feature of Theorem \ref{t:in1} is that the Robin function
plays a role in concentration phenomena only for the case of
simply connected domains.

The regular part of the Green's function plays an important role
in many problems with critical exponent in dimension larger than
two, see \cite{bab}, \cite{blr}, \cite{bk}, \cite{han}, \cite{r0},
\cite{r1}. The difference here is that the regular part does not
appear directly in the expansion, we find the above function
$\tilde{H}$ instead, and we recover the regular part from the
Riemann mapping Theorem. The Robin function is also related to the
notion of conformal incenter, see \cite{fl}.

%

The expansion of $I_\e(u)$ for multiple bubbles is performed in
Section \ref{s:2b}, see Proposition \ref{p:km}. It turns out that
when $\l_i \sim \e^{-1}$ for all $i$, the mutual interaction among
the bubbles is of the same order as the interaction with the
boundary (through both the geometry of $\O$ and the datum $g$). We
observe that the interaction among the bubbles depend on their
mutual orientation.

There is a by-product of the expansion in Proposition \ref{p:km}.
It allows us to settle a particular case of a question raised by
Brezis and Coron, see Section \ref{s:bou}. In \cite{brc2} the
authors consider a sequence of solutions $u_n$ of \eqref{eq:in1}
and a sequence $g_n$ of boundary data which converge to zero in
$H^{\frac 12}(\partial \O) \cap L^\infty(\partial \O)$. Under
these conditions they prove that the sequence $u_n$ {\em splits}
into a finite number of bubbles, and their image converge to a
finite and connected union of spheres of radius 1. They ask
whether every configuration of spheres can be obtained as a limit
of solutions $u_n$ for a suitable sequence of boundary data $g_n$.
We have an affirmative answer if all the spheres pass through the
origin.

\begin{thm}\label{t:in2}
Let $D$ denote the unit disk in $\R^2$, and let $A = \{S_1 \cup
\dots \cup S_k\}$ be any configuration of unit spheres, each
passing through the origin of $\R^3$. Then there exist a sequence
$\tilde{g}_n : S^1 \to \R^3$ and a sequence of functions $u_n$
solving
\begin{equation*}\label{eq:in7}
  \begin{cases}
    \D v = 2 v_x \wedge v_y & \text{ in } D, \\
    v = \tilde{g}_n  & \text{ on } \partial D,
  \end{cases}
\end{equation*}
such that the image of the function $u_n$ converge to $A$ in the
Hausdorff sense.
\end{thm}

\noindent The functions $u_n$ in Theorem \ref{t:in2} are
constructed studying the interactions of the bubbles (of degree
$1$) with the boundary datum and among themselves. Choosing
boundary data with an appropriate {\em strong concentration} at
$k$ points on $\partial D$, we show that the self interaction
among the bubbles becomes negligible. Hence we can find solutions
$u_n$ which are highly concentrated at $k$ points close to the
boundary of $D$ and with prescribed orientations in $\R^3$. We
remark that the order of concentration, roughly the parameter $\l$
in \eqref{eq:in4}, turns out to be the same for all the bubbles.

The case of spheres not passing through the origin is not treated
here. We believe that it could be possible to achieve such
configurations by considering bubbles with higher degree. In fact,
in the recent paper by A. Bahri and S. Chanillo, \cite{bach}, the
authors showed that when considering changing-sign solutions of
the Yamabe problem, the bubbles can exhibit different orders of
concentration. If there is an analogy between higher-degree
bubbles and changing-sign solutions of the Yamabe equation, then
one could obtain bubbles with higher and higher concentration and
with image not passing through the origin. This will be the object
of a future work.

\begin{center}
{\bf Acknowledgments}
\end{center}

\noindent S.C. is supported by a NSF grant. A. M. has been
supported by a Fulbright fellowship for the academic year
2000-2001 and by MURST, under the project {\em Variational Methods
and Nonlinear Differential Equations}. He is also grateful to the
Mathematics Department at Rutgers University and IAS for the kind
hospitality, where part this work has been accomplished. The
authors wish to thank A. Bahri for useful discussions, and are
very grateful to P. Caldiroli and R. Musina for their helpful
comments.

\section{Notation and preliminary facts}\label{s:pre}

In this Section we introduce some notation and preliminary facts
in order to tackle problem \eqref{eq:inge}.

In the following $D$ will denote the unit disk in $\R^2$
$$
D = \left\{ \xi = (x,y) \in \R^2 \, : \, x^2 + y^2 < 1 \right\}.
$$
Let $\ov{I} : \mathcal{D} \to \R$ be defined by \eqref{eq:ovin}.
From \cite{cl} the last term in $\ov{I}(u)$ is well defined on
$\mathcal{D}$, together with its Frechet derivatives. This makes
$\ov{I}$ a smooth functional on $\mathcal{D}$. The same argument
provides regularity of the functional $I_\e$ on $H^1_0(\O;\R^3)$.

Using a finite-dimensional reduction, we are going to treat the
functional $I_\e$ as a perturbation of $\ov{I}$. In order to do
this, it is essential to consider the critical points of the
functional $\ov{I}$, namely the solutions of
\begin{equation}\label{eq:pinf}
  \D u = 2 u_x \wedge u_y \qquad \hbox{ in } \R^2, \qquad
  u \in \mathcal{D}.
\end{equation}
The stereographic projection \eqref{eq:stpr} is indeed a solution
of \eqref{eq:pinf}, which we call {\em fundamental solution} or
{\em bubble}. By invariance its translations, dilations and
rotations are also solutions of \eqref{eq:pinf}. We set
\begin{equation}\label{eq:ovz}
  \ov{Z} = \left\{ R \, \d_{a,\l} (\cdot) = R \d_\l (\cdot - a)
  \, : \, \l > 0, a \in \R^2, R \in SO(3) \right\}.
\end{equation}
We remark that, since $SO(3)$ is a
three-dimensional manifold, $\ov{Z}$ is a six-dimensional
manifold.


We list now some useful expressions. Note that the function $\d
\left( \l \, (\xi - a) \right)$ has the explicit form
$$
\d \left( \l \, (\xi - a) \right) = \left( \frac{2 \l (\xi - a)}{1
+ \l^2 |\xi - a|^2}, \frac{\l^2 |\xi-a|^2 - 1}{\l^2 |\xi-a|^2 + 1}
\right),
$$
from which, if $a$ is bounded away from $\partial \O$, one can
deduce
\begin{equation}\label{eq:apd}
\pi \left( \l \, (\xi - a) \right) \sim \left( \frac{2}{\l}
\frac{(\xi - a)}{|\xi - a|^2}, 1 - \frac{2}{\l^2}
\frac{1}{|\xi-a|^2} \right) + O(\l^{-3}), \qquad \hbox{ for } \l
\hbox{ large}.
\end{equation}
Writing for brevity $\d$ instead of $\d_{a,\l}$, we compute some
derivatives of $\d$. We emphasize that throughout the paper,
unless explicitly stated, the point $a$ will always be bounded
away from $\partial \O$, namely we will assume dist$(a,\partial
\O) \geq \t_0$ for some fixed $\t_0 > 0$. We have
\begin{eqnarray}
(\d_x)_1 = 2 \l \frac{1 + \l^2 (y^2-x^2)}{(1 + \l^2 (x^2+y^2))^2};
\qquad (\d_x)_2 = - \frac{4 \l^3 x y}{(1 + \l^2 (x^2+y^2))^2};
\qquad (\d_x)_3 = \frac{4 \l^2 x}{(1 + \l^2 (x^2+y^2))^2};
\nonumber \\ \label{eq:dxy} \\
(\d_y)_1 = - \frac{4 \l^3 x y}{(1 + \l^2 (x^2+y^2))^2}; \qquad
(\d_y)_2 = 2 \l \frac{1 + \l^2 (x^2-y^2)}{(1 + \l^2 (x^2+y^2))^2};
\qquad (\d_y)_3 = \frac{4 \l^2 y}{(1 + \l^2 (x^2+y^2))^2}.
\nonumber
\end{eqnarray}
From the last two formulas we deduce
\begin{equation}\label{eq:dxdyxy}
  (\d_x \wedge \d_y)_1 = - \frac{8 \l^3 x}{(1 + \l^2 (x^2+y^2))^3};
  \qquad (\d_x \wedge \d_y)_2 = - \frac{8 \l^3 y}{(1 + \l^2 (x^2+y^2))^3};
\end{equation}
\begin{equation}\label{eq:dxdyz}
  (\d_x \wedge \d_y)_3 = 4 \l^2 \frac{1 - \l^2 (x^2+y^2)}{(1 + \l^2
  (x^2+y^2))^3}.
\end{equation}
The functions $R \d_{a,\l}|_\O$ do not belong to $H^1_0(\O;\R^3)$
since they are non-zero at the boundary. Following \cite{bab},
\cite{r1}, it is convenient to project these functions on the
space $H^1_0(\O;\R^3)$, by subtracting the harmonic function on
$\O$ with the same boundary data. Let $\varphi : \O \to \R^3$ be
the unique solution of the problem
\begin{equation}\label{eq:fi}
  \begin{cases}
    \D \varphi  = 0 & \text{ in } \O, \\
    \varphi = \d & \text{ on } \partial \O,
  \end{cases}
\end{equation}
and set $P \d = \d - \varphi$. We will often omit the dependence
of $\varphi$ on the parameters $a, \l, R$, as for $\d$.

From \eqref{eq:dxy} and some standard computations it is easy to
find that, in the case $R = Id$
\begin{equation}\label{eq:estfi}
  \var = \left( \frac{2}{\l} h_1(\xi,a) + o(\l^{-1}),
  \frac{2}{\l} h_2(\xi,a) + o(\l^{-1}), 1 - \frac{2}{\l^2}
  h_3(\xi,a) + o(\l^{-2}) \right),
\end{equation}
and
\begin{equation}\label{eq:estdf}
  (\d - \var) = \begin{pmatrix}
    \frac{2 \l (x - a_1)}{1 +
  \l^2 |\xi - a|^2} - \frac{2}{\l} h_1(\xi,a) + o(\l^{-1}) \\
   \frac{2 \l (y - a_2)}{1 + \l^2 |\xi - a|^2} - \frac{2}{\l}
  h_2(\xi,a) + o(\l^{-1}) \\ - \frac{2}{1 + \l^2 |\xi - a|^2} +
  \frac{2}{\l^2} h_3(\xi,a) + o(\l^{-2}) \
  \end{pmatrix},
\end{equation}
where $h_1$ and $h_2$ are defined in \eqref{eq:h1h2}, and where
$h_3$ is the solution of
\begin{equation}\label{eq:h3}
\begin{cases}
    \D_\xi h_3 (a,\xi) = 0 & \text{ in } \O, \\
    h_3(a,\xi) = \frac{1}{|\xi-a|^2} & \text{ on } \partial
    \O.
  \end{cases}
\end{equation}
The quantities $o(\l^{-1})$ and $o(\l^{-2})$ in formulas
\eqref{eq:estfi} and \eqref{eq:estdf} denote functions which
$C^k(\ov{\O})$-norm, for any $k \in \N$, is of order $o(\l^{-1})$
and $o(\l^{-2})$ respectively.

We collect some further estimates, whose proof are trivial, and
which we will use later. Given a fixed positive constant $\t \leq
\frac{\t_0}{2}$, for $\l$ sufficiently large there holds
\begin{equation}\label{eq:est}
  \int_{B_\t} |P \d| \leq \frac{C}{\l}; \qquad
  \int_{B_\t} |\n \d| \leq C \frac{\log \l}{\l}; \qquad
  \int_{B_\t} |\n P \d| \leq C \frac{\log \l}{\l}; \qquad
  \int_{B_\t} |\d|^2 \leq C \frac{\log \l}{\l^2};
\end{equation}
\begin{equation}\label{eq:est2}
  |P \d|(x) + |\n \d|(x) \leq \frac{C}{\l}, \quad \forall
  x \in \O \setminus B_{\t}; \qquad |\var - (0,0,1)|(x) +
  |\n \var|(x) \leq \frac{C}{\l}, \quad \forall x \in \O.
\end{equation}
In the following, for brevity of notation, the constant $C$ will
be allowed to vary from formula to formula and from line to line.

For $k \geq 1$ and for $i, j \in \{1, \dots, k\}$, $i \neq j$, we
will use the following notation
$$
\tilde{e}(\e,\l_1, \dots\l_k) = O(\e^2) + \sum_{i=1}^k O(\e
\l_i^{-1} |\log \l_i|) + \sum_{i,j} O \left( \frac{1}{\l_i \l_j}
\right);
$$
$$
e(\e, \l_i) = o(\e^2) + o(\e \l_i^{-1}); \qquad e(\l_i, \l_j) =
O\left( (\log \l_i + \log \l_j) \left( \frac{1}{\l_i^3} +
\frac{1}{\l_j^3} + \frac{1}{\l_i^2 \l_j} + \frac{1}{\l_i \l_2^j}
\right) \right);
$$
$$
e(\e,\l_1,\dots,\l_k) = o(\e^2) + \sum_{i=1}^k o(\e \l_i^{-1}) +
\sum_{i<j} e(\l_i,\l_j) + \sum_{i<j<k} O \left( \frac{1}{\l_i \l_j
\l_k} \right).
$$
We will often make use of the identity
\begin{equation}\label{eq:id}
  \int_{\R^2} \frac{1 - |\xi|^2}{(1+|\xi|^2)^3} = 0,
\end{equation}
which is immediate to verify (the integrand is indeed the third
component of $\D \d$ up to a constant).

\section{A non-existence result via a Pohozahev-type identity}

In this section we prove a Pohozaev-type identity for the
$H$-surface equation. The proof is elementary and extends a
previous result of Wente, see \cite{w1}.

\begin{pro}\label{p:poh}
Let $\O \subseteq \R^2$ be a smooth bounded and simply-connected
domain, and let $v \in C^2(\ov{\O};\R^3)$ be a solution of
\begin{equation}\label{eq:eqh}
  \begin{cases}
    \D v = H(\xi,v,\n v) v_x \wedge v_y, & \text{ in } \O, \\
    v = 0 & \text{ on } \partial \O,
  \end{cases}
\end{equation}
for some continuous function $H(\xi,v,\n v)$. Then $v \equiv 0$ in
$\O$.
\end{pro}

\begin{pf}
We assume first that the domain $\O$ is the unit disk $D$. In the
spirit of the Pohozahev identity, we consider the quantity
$\sum_{i=1}^3 (\xi \cdot \n v_i) \D v_i$ and integrate on $D$. We
claim that
\begin{equation}\label{eq:vl}
  \sum_{i=1}^3
  (\xi \cdot \n v_i) (v_x \wedge v_y)_i = 0, \qquad \qquad
  \xi = (x,y).
\end{equation}
Once \eqref{eq:vl} is proved, we have
\begin{equation}\label{eq:vl2}
  \sum_{i=1}^3 (\xi \cdot \n v_i) \D v_i \equiv 0.
\end{equation}
Integrating \eqref{eq:vl2} over $D$ and taking into account that
the dimension is 2, we find
$$
\frac{1}{2} \sum_{i=1}^3 \int_{\partial D} \left( \frac{\partial
v_i}{\partial \nu} \right)^2 = \frac{1}{2} \sum_{i=1}^3
\int_{\partial D} (\xi \cdot \nu) \left( \frac{\partial
v_i}{\partial \nu} \right)^2 = 0,
$$
where $\nu$ denoted the exterior unit normal to $\partial D$. As a
consequence we have $\frac{\partial v}{\partial \nu} = 0$ on
$\partial D$. Thus, extending $v$ to zero on the complement of $D$
and also extending $H$ continuously outside $D$ we obtain a $C^1$
solution of $\D v = H(\xi,v,\n v) v_x \wedge v_y$ in $\R^2$.
Hence, applying Theorem 1 in \cite{hw} we obtain $v \equiv 0$ in
$D$. Let us now verify \eqref{eq:vl}: using simple computation we
find
\begin{eqnarray*}
 \sum_{i=1}^3 (\xi \cdot \n v_i) (v_x \wedge v_y)_i & = &
 x \left[ (v_1)_x (v_2)_x (v_3)_y - (v_1)_x (v_3)_x (v_2)_y +
 (v_2)_x (v_1)_y (v_3)_x \right. \\ & - & \left.
 (v_2)_x (v_1)_x (v_3)_y +
 (v_3)_x (v_1)_x (v_2)_y - (v_3)_x (v_2)_x (v_1)_y \right] \\
  & + &
 y \left[ (v_1)_x (v_2)_x (v_3)_y - (v_1)_x (v_3)_x (v_2)_y +
 (v_2)_x (v_1)_y (v_3)_x \right. \\ & - & \left.
 (v_2)_x (v_1)_x (v_3)_y +
 (v_3)_x (v_1)_x (v_2)_y - (v_3)_x (v_2)_x (v_1)_y \right] = 0.
\end{eqnarray*}
This concludes the proof in the case $\O = D$. For the general
case of a simply-connected domain, it is sufficient to use the
Riemann Mapping Theorem and the transformation rule of
\eqref{eq:eqh} under conformal mappings. We recall that for $\O$
smooth, the Riemann map is also smooth up to the boundary, see
\cite{war}.
\end{pf}

\section{The finite-dimensional reduction}\label{s:4}

In this section we show how problem \eqref{eq:inge} can be reduced
to a finite-dimensional one for small values of $\e$. The starting
point is the following Proposition, proven in \cite{isa1} (Lemma
5.5) using an isoperimetric inequality. We give an alternative
proof in the Appendix, using the stereographic projection and
shifting the problem from $\R^2$ to $S^2$. We believe that our
proof could be naturally extended to the case of higher degree
bubbles.

\begin{pro}\label{p:up}
There exists a constant $C_0 > 0$ such that
$$
\ov{I}''(R \d_{a,\l})[R
\d_{a,\l},R \d_{a,\l}] \leq - C_0 \|R \d_{a,\l}\|^2, \qquad \hbox{
for all } R \d_{a,\l} \in \ov{Z}. $$
and $$
  \ov{I}''(R \d_{a,\l})[v,v] \geq C_0 \|\n v\|^2_{L^2(\R^2)}, \qquad
  \hbox{ for all } R \d_{a,\l} \in \ov{Z} \hbox{ and } v \perp
  \left( T_{R \d_{a,\l}} \ov{Z} \oplus \{t R \d_{a\,\l}\}_t \right).
$$
In particular, the equation $\ov{I}''(R \d_{a,\l}) [v] = 0$
implies $v - c \in T_{R \d_{a,\l}} \ov{Z}$ for some $c \in \R^3$.
\end{pro}

\noindent We are going to consider now problem \eqref{eq:pinf} on
the domain $\O$. Given $\ov{C} > 0$ we set
\begin{equation}\label{eq:manz}
  Z = \left\{ \sum_{i=1}^k P \d_{p_i,\l_i} \, : \,
  dist(p_i, \partial \O) \geq \ov{C}^{-1}, \quad dist(p_i,p_j)
  \geq \ov{C}^{-1} \; \forall i \neq j, \quad \l_i \, \e \in
  [\ov{C}^{-1},\ov{C}], R_i \in SO(3) \right\}.
\end{equation}
Proposition \ref{p:up} asserts that the manifold $\ov{Z}$, see
\eqref{eq:ovz} is non-degenerate for the functional $\ov{I}$
module translations. As a consequence, it is easy to extremize
$I_\e$ in the direction perpendicular to $Z$. This is stated in
the following Proposition \ref{p:ambbad}, in the same spirit as
\cite{ab}. We need first a preliminary Lemma (see also \cite{bab},
Proposition 3.1).

\begin{lem}\label{l:inv}
Let $k \in \N$, $\ov{C} > 0$, and let $Z$ be as in
\eqref{eq:manz}. Then there exists a positive constant $C$ such
that
$$
\hbox{ if } \qquad v \in H^1_0(\O), v \perp T_z Z, v \perp P \d_i
\, \forall i, \qquad \hbox{ then } \qquad I''_\e(z) [v,v] \geq
C^{-1} \|v\|_{H^1_0(\O)}^2.
$$
\end{lem}

\begin{pf} For $i = 1, \dots, k$, let $B_i$ be the ball
of radius $\frac{\ov{C}^{-1}}{2}$ around $p_i$, and let also
$\tilde{B}_i$ be the ball of radius $\frac{\ov{C}^{-1}}{4}$ around
$p_i$. Let us denote by $P_i$ the orthogonal projection of
$H^1_0(\O)$ onto $H^1_0(B_i)$, and for any $v \in H^1_0(\O)$ set
$v_1 = v - \sum_{i=1}^k P_i v$. It follows immediately that
\begin{equation}\label{eq:decv}
  \|v\|_{H^1_0(\O)}^2 = \sum_{i = 1}^k \|P_i v\|_{H^1_0(\O)}^2 +
  \|v_1\|_{H^1_0(\O)}^2.
\end{equation}
From standard regularity results, since the function $v_1$ is
harmonic in each $B_i$, and since it coincides with $v$ on each
$\partial B_i$, there holds
\begin{equation}\label{eq:stv1}
  \|v_1\|_{C^2(\tilde{B}_i)} \leq C \|v\|_{H^1_0(\O)}, \qquad
  \hbox{ for all } i = 1, \dots, k,
\end{equation}
where $C$ is a constant independent of $v$. Since $v$ is
orthogonal to $P \d_i$, from \eqref{eq:est2} we deduce
$$
(P_i v, P \d_i) = \left( v - v_1 - \sum_{j \neq i} P_j v, P \d_i
\right) = - (v_1, P \d_i) + \sum_{j \neq i} O \left(
\frac{1}{\l_j} \right) \|v\|_{H^1_0(\O)}.
$$
To evaluate the scalar product $(v_1, P \d_i) = \int_\O \n v_1
\cdot \n P \d_i$, we divide the integral in the regions
$\tilde{B}_i$ and $\O \setminus \tilde{B}_i$. We have clearly
$\int_{\O \setminus \tilde{B}_i} \n v_1 \cdot \n P \d_i =
O(\l_i^{-1}) \|v\|_{H^1_0(\O)}$. On the other hand, using
\eqref{eq:est} and \eqref{eq:stv1} we find
$$
\left| \int_{\tilde{B}_i} \n v_1 \cdot \n P \d_i \right| = O
\left( \int_{\tilde{B}_i} |\n P \d_i| \right) \|v\|_{H^1_0(\O)}
\leq O \left( \frac{\log \l_i}{\l_i} \right) \|v\|_{H^1_0(\O)}.
$$
Using these formulas and \eqref{eq:est2} we obtain
\begin{equation}\label{eq:ortzpzi}
  \left|(P_i v, \d_i) \right| = O \left( \frac{\log \l_i}{\l_i}
  \right) \|v\|_{H^1_0(\O)} + \sum_{j \neq i} O \left( \frac{1}{\l_j}
  \right) \|v\|_{H^1_0(\O)} \leq C \e |\log \e| \, \|v\|_{H^1_0(\O)}.
\end{equation}
In the same way as \eqref{eq:ortzpzi}, using the explicit
expression of the function $\d_i$ and taking the scalar product of
$v$ with $\frac{\partial \d_i}{\partial p_i}$, $\frac{\partial
\d_i}{\partial R_i}$ and $\frac{\partial \d_i}{\partial \l_i}$,
one finds $\|\Pi_i P_i v\| \leq C \e |\log \e| \|v\|_{H^1_0(\O)}$,
where $\Pi_i$ denotes the orthogonal projection onto the space
spanned by $\d_i$, $\frac{\partial \d_i}{\partial p_i}$,
$\frac{\partial \d_i}{\partial R_i}$ and $\frac{\partial
\d_i}{\partial \l_i}$. The functional $I''_\e(z)$ is given by
$$
I''_\e(z) [v,\tilde{v}] = \int_\O \n v \cdot \n \tilde{v} - 2
\int_\O z \cdot \left( v_x \wedge \tilde{v}_y + \tilde{v}_x \wedge
v_y \right) + 2 \e \int_\O \tilde{v} \cdot \left( g_x \wedge v_y +
v_x \wedge g_y \right), \quad v, \tilde{v} \in H^1_0(\O).
$$
It follows easily from the expression of $I''_\e$ and from
Proposition \ref{p:up} that
\begin{equation}\label{eq:i''piv}
  I''_\e(z) [P_i v, P_i v] \geq C_0 \|P_i v\|_{H^1_0(\O)}^2 - C \e |\log \e|
  \, \|v\|_{H^1_0(\O)}.
\end{equation}
For an arbitrary function $v$ there holds
\begin{equation}\label{eq:i''de}
  I''_\e (z) [v,v] = \sum_{i = 1}^k I''_\e (z) [P_i v, P_i v] +
  I''_\e (z) [v_1,v_1] + 2 \sum_{i = 1}^k I''_\e (z) [P_i v, v_1].
\end{equation}
From the orthogonality of $P_i v$ and $v_1$ it follows that
\begin{eqnarray*}
  I''_\e (z) [P_i v, v_1] & = & - 2 \int_\O z \cdot \left( (P_i v)_x
  \wedge (v_1)_y + (v_1)_x \wedge (P_i v)_y \right) + 2 \e \int_\O
  P_i v \cdot \left( g_x \wedge (v_1)_y + (v_1)_x \wedge g_y
  \right) \\ & = & O \left( \int |z| |\n P_i v| |\n v_1| \right) +
  O(\e) \|v\|_{H^1_0(\O)}^2.
\end{eqnarray*}
Dividing again the integral into the regions $\tilde{B}_i$ and $\O
\setminus \tilde{B}_i$ we deduce
\begin{equation}\label{eq:i''piv1}
  I''_\e (z) [P_i v, v_1] = O \left( \frac{\log \l_i}{\l_i} \right)
  \|v\|^2_{H^1_0(\O)} + \sum_{j \neq i} O \left( \frac{1}{\l_j} \right)
  \|v\|^2_{H^1_0(\O)} + O(\e) \|v\|^2_{H^1_0(\O)}.
\end{equation}
Similarly, we obtain
\begin{equation}\label{eq:i''piv2}
  I''_\e (z) [v_1, v_1] = \|v_1\|^2_{H^1_0(\O)} +
  O \left( \frac{\log \l_i}{\l_i} \right)
  \|v\|^2_{H^1_0(\O)} + \sum_{j \neq i} O \left( \frac{1}{\l_j} \right)
  \|v\|^2_{H^1_0(\O)} + O(\e) \|v\|^2_{H^1_0(\O)}.
\end{equation}
From \eqref{eq:decv}, \eqref{eq:i''piv}, \eqref{eq:i''de},
\eqref{eq:i''piv1} and \eqref{eq:i''piv2} the Lemma follows.
\end{pf}

\begin{pro}\label{p:ambbad}
Let $\ov{C}$ be a fixed positive constant, let $k \in \N$, let $\e
> 0$, and $Z$ be defined as above. Then, if $\e$ is sufficiently
small, for every $z \in Z$ there exist a function $w_\e(z) \in
H^1(\O;\R^3)$ and $C > 0$ with the following properties
\begin{description}
\item{i)} $w_\e(z)$ is orthogonal to $T_z Z$, \qquad for all
  $z \in Z$;
\item{ii)} $I'_\e (z + w_\e(z)) \in T_z Z$, \qquad for all
  $z \in Z$;
\item{iii)} $\| w_\e(z) \| \leq C \|I'_\e(z)\|$, \qquad for all
  $z \in Z$;
\end{description}
By $i)$ and $ii)$, the manifold
$$
Z_\e = \{z + w_\e(z) \, : \, z \in Z \}
$$
is a natural constraint for $I'_\e$. Namely if $u \in Z_\e$ and
$I'_\e|_{Z_\e}(u) = 0$, then $I'_\e(u) = 0$.
\end{pro}

\

\begin{pf}
Given Proposition \ref{p:up}, the arguments are quite standard.
For convenience, we give a brief sketch in the case $k = 1$. In
the proof, we simply write $\d$ for $R \d_{a,\l}$.

Let us define $F_\e:  Z \times H^1(\O; \R^3) \times T_z Z \to
H^1(\O; \R^3) \times \mathbb{R}$ by setting
$$
F_\e (z,w,q) = \left( \begin{array}{cc} I'_\e(z+w) - q
\\ (w, q) \end{array} \right).
$$
With this notation, the unknown $(w, q) = (w_\e,
I'_\e(z+w_\e(z)))$ can be implicitly defined via the equation
$F_\e (z,w,q)=(0,0)$. Setting $G_\e(z,w,q) = F_\e(z,w,q) -
\partial_{(w,q)} F_\e(z,0,0)[(w,q)]$ we have that
$$
F_\e(z,w,q) = 0 \quad \Leftrightarrow \quad
\partial_{(w,q)} F_\e(z,0,0)[(w,q)] + G_\e(z,w,q) = 0.
$$
Reasoning as in \cite{ab}, using Lemma \ref{l:inv} one can prove
that $\partial_{(w,q)} F_\e(z,0,0)$ is uniformly invertible for $z
\in Z$ and $\e$ sufficiently small. Hence we can write
$$
F_\e(z,w,q) = 0 \quad \Leftrightarrow \quad (w,q) = W_\e(z,w,q) :=
- \left(\partial_{(w,q)} F_\e(z,0,0) \right)^{-1} \left[
F_\e(z,0,0) + Q_\e(z,w,q) \right],
$$
where
$$
Q_\e(z,w,q) = F_\e(z,w,q) - F_\e(z,0,0) - \partial_{(w,q)}
F_\e(z,0,0)[(w,q)].
$$
It is also standard to prove that $Q_\e(z,w,q)$ satisfies
\begin{equation}\label{eq:q}
  \begin{cases}
  \| Q_\e(z,w,q) \| \leq C \| (w,q) \|^2 & \\
  \| Q_\e(z,w,q) -
  Q_\e(z,\tilde{w},\tilde{q}) \| \leq C \left( \| (w,q) \| +
  \| (\tilde{w},\tilde{q}) \| \right) \| (w,q) -
  (\tilde{w},\tilde{q}) \|, &
  \end{cases}
\end{equation}
where $\| (w,q) \|$ and where $C = C(\O, g, \ov{C})$ is a constant
depending on $\O$, $g$, $\ov{C}$, and independent of $z \in Z$ and
$\e$. Using \eqref{eq:q} it is possible to prove that the function
$W_\e$ is a contraction in a ball of radius $\tilde{C}
\|F_\e(z,0,0)\|$ for some positive constant $\tilde{C}(\O, g,
\ov{C})$. Since $\|F_\e(z,0,0)\| \leq C \|I'_\e(z)\|$ for some
constant $C$, the conclusion follows.
\end{pf}

\

\noindent We estimate now the quantity $\|I'_\e(\sum P \d_i)\|$ in
order to control the norm of $w_\e(z)$, see $iii)$ in Proposition
\ref{p:ambbad}.

\begin{lem}\label{l:gr}
Let $C$ be a fixed positive constant, let $k \in \N$, let $\e >
0$, and let $Z$ be as in \eqref{eq:manz}. Then there holds
$$
\|I'_\e(z)\| \leq \tilde{e}(\e,\l_1, \dots, \l_k), \qquad \hbox{
for } \e \hbox{ sufficiently small and for all } z \in Z,
$$
where $\tilde{e}(\e,\l_1, \dots, \l_k)$ is defined in Section
\ref{s:pre}.
\end{lem}

\begin{pf}
Let $v \in H^1_0(\O; \R^3)$ and $z \in Z$. Using integration by
parts we deduce easily
\begin{equation}\label{eq:grieu}
I'_\e(z)[v] = \int_{\O} \n z \cdot \n v + 2 \int_{\O} v \cdot (z_x
\wedge z_y) + 2 \e \int_{\O} z \cdot (v_x \wedge g_y + g_x \wedge
v_y) + 2 \e^2 \int_{\O} v \cdot (g_x \wedge g_y).
\end{equation}
From H\"older's inequality we get
$$
\left| \int_{\O} z (v_x \wedge g_y + g_x \wedge v_y) \right| \leq
\|g\|_{C^1(\ov{\O})} \left( \int_{\O} |z|^2 \right)^{\frac 12}
\|v\|, \qquad \hbox{ for all } v \in H^1(\O; \R^3).
$$
From \eqref{eq:est}, \eqref{eq:est2} it is easy to check that
$\left( \int_{\O} |z|^2 \right)^{\frac 12} \leq C \sum_{i=1}^k
\frac{1}{\l_i} |\log \l_i|^{\frac{1}{2}}$, hence we have
$$
\left| \e \int_{\O} z \cdot (v_x
\wedge g_y + g_x \wedge v_y) \right| \leq C \sum_{i=1}^k
\frac{\e}{\l_i} |\log \l_i|^{\frac{1}{2}} \|v\|, \qquad \hbox{ for
all } v \in H^1_0(\O; \R^3).
$$
On the other hand, it is also immediate to verify the inequality
$$
\e^2 \left| \int_{\O} v \cdot
(g_x \wedge g_y) \right| \leq C \e^2 \|v\|, \qquad \hbox{ for all
} v \in H^1_0(\O) \hbox{ and all } z \in Z.
$$
It
remains to estimate the first two terms in \eqref{eq:grieu}.
Writing for brevity $\d_i = R_i \d_{p_i,\l_i}$, we have also
$$
\int_{\O} \n z \cdot \n v + 2 \int_{\O} v \cdot (z_x \wedge z_y) =
\int_{\O} v \cdot \left[ - \D \left(\sum P \d_i\right) + 2
\left(\left(\sum P \d_i\right)_x \wedge \left(\sum P
\d_i\right)_y\right) \right].
$$
Using the equation $\D \d_i = 2 ((\d_i)_x \wedge (\d_i)_y)$ and
the fact that $\D \d_i = \D (P \d_i)$, the above quantity becomes
$$
2 \int_{\O} v \cdot \left[ \left(\sum P \d_i\right)_x \wedge
\left(\sum P \d_i\right)_y - \sum (\d_i)_x \wedge (\d_i)_y
\right],
$$
which can be written as
\begin{equation}\label{eq:grd3}
2 \sum_{i \neq j} \int_{\O} v \cdot \left( (P \d_i)_x \wedge (P
\d_j)_y \right) + 2 \sum_i \int_{\O} v \cdot \left[ (P \d_i -
\d_i)_x \wedge (\d_i)_y + (\d_i)_x \wedge (P \d_i - \d_i)_y
\right].
\end{equation}
Let us estimate first the term $\int_{\O} v \cdot (P \d_i)_x
\wedge (P \d_j)_y$. Let $\g < \frac{1}{2} \ov{C}^{-1}$ (recall the
definition of $Z$) be a fixed positive number and divide the
integral in the three regions $B_{\g}(p_i)$, $B_{\g}(p_j)$ and $\O
\setminus (B_{\g}(p_i) \cup B_{\g}(p_j))$. Integrating by parts on
the balls $B_{\g}(p_i)$ and $B_{\g}(p_j)$, the quantity $\int_{\O}
v \cdot (P \d_i)_x \wedge (P \d_j)_y$ becomes
\begin{eqnarray*}
& - & \int_{B_\g(p_i)} v_x \cdot (P \d_i \wedge (P \d_j)_y) -
\int_{B_\g(p_i)} v \cdot (P \d_i \wedge (P \d_j)_{xy}) +
\int_{\partial B_\g(p_i)} v \cdot (P
\d_i \wedge (P \d_j)_y) \nu_x  \\
& - &  \int_{B_\g(p_j)} v_x \cdot(P \d_i \wedge (P \d_j)_y) -
\int_{B_\g(p_j)} v \cdot (P \d_i \wedge (P \d_j)_{xy}) +
\int_{\partial B_\g(p_j)} v \cdot (P \d_i \wedge (P \d_j)_y) \nu_x  \\
& + & \int_{\O \setminus (B_{\g}(p_i) \cup B_{\g}(p_j))} v \cdot
\left( (P \d_i)_x \wedge (P \d_j)_y \right).
\end{eqnarray*}
Hence, since $\d_i$ and its derivatives are of order $\l_i^{-1}$
in $\O \setminus B_\g(p_i)$ one finds
$$
\left| \int_{\O} v \cdot \left( (P \d_i)_x \wedge (P \d_j)_y
\right) \right| \leq C \frac{1}{\l_i \l_j} \|v\|.
$$
Using similar estimates we find that the whole expression in
\eqref{eq:grd3} is of order $\e^2$. So we obtain the conclusion.
\end{pf}

\begin{lem}\label{l:i2zp}
Let $\ov{C} > 0$ and let $Z$ be as in \eqref{eq:manz}. Then there
holds
$$
\left| I''_\e(z) \left[ \frac{\partial z}{\partial p_i} \right]
\right| + \left| I''_\e(z) \left[ \frac{\partial z}{\partial R_i}
\right] \right| \leq C \left( \e + \sum_{j=1}^{k} \l_j^{-1}
\right); \qquad \left| I''_\e(z) \left[ \frac{\partial z}{\partial
\l_i} \right] \right| \leq \frac{1}{\l_i} C \left( \e +
\sum_{j=1}^{k} \l_j^{-1} \right),
$$
for every $i = 1, \dots, k$.
\end{lem}

\begin{pf}
From \eqref{eq:grieu} and some integration by parts it follows
that
$$
I''_\e(z) [v,\tilde{v}] = \int_\O \n v \cdot \n \tilde{v} + 2
\int_\O \tilde{v} \cdot \left( z_x \wedge v_y + v_x \wedge z_y
\right) + 2 \e \int_\O \tilde{v} \cdot \left( g_x \wedge v_y + v_x
\wedge g_y \right),
$$
where $v, \tilde{v}$ are arbitrary functions in $H^1_0(\O, \R^3)$.
We choose now $\tilde{v} = \frac{\partial z}{\partial p_i}$, and
we let $v$ be an arbitrary test function. We have clearly
$\frac{\partial z}{\partial p_i} = \frac{\partial \d_i}{\partial
p_i} - \frac{\partial \var_i}{\partial p_i}$, where $\var_i$ is
the function in \eqref{eq:fi} corresponding to $\d_i$. From the
estimates in \eqref{eq:est2} and from the explicit expression of
$\frac{\partial \d_i}{\partial p_i}$ we find
\begin{equation}\label{eq:add1}
  \e \left| \int \frac{\partial z}{\partial p_i} \cdot \left( g_x \wedge
  v_y + v_x \wedge g_y \right) \right| \leq C \frac{\e}{\l_i} \|v\| +
  C \e \left\| \frac{\partial \d_i}{\partial p_i} \right\|_{L^2(\O)}
  \|v\| \leq C \e \|v\|.
\end{equation}
Turning to the remaining two terms, we have
\begin{eqnarray}\label{eq:add2}
  & & \int_\O \n v \cdot \n \frac{\partial z}{\partial p_i} + 2 \int_\O
\frac{\partial z}{\partial p_i} \cdot \left( z_x \wedge v_y + v_x
\wedge z_y \right) = \int_\O \n v \cdot \n \frac{\partial
\d_i}{\partial p_i} + 2 \int_\O \frac{\partial \d_i}{\partial p_i}
\cdot \left( (\d_i)_x \wedge v_y + v_x \wedge (\d_i)_y \right)
\nonumber \\ & & + 2 \sum_{j \neq i} \int_\O \frac{\partial
\d_i}{\partial p_i} \cdot \left( (\d_j)_x \wedge v_y + v_x \wedge
(\d_j)_y \right) + O(\l_i^{-1}) \|v\|.
\end{eqnarray}
Integrating by parts and using the fact that $\n \frac{\partial
\d_i}{\partial p_i}$ is of order $\l_i^{-1}$ on $\partial \O$ we
find
\begin{equation}\label{eq:add3}
  \int_\O \n v \cdot \n \frac{\partial \d_i}{\partial p_i} + 2
\int_\O \frac{\partial \d_i}{\partial p_i} \cdot \left( (\d_i)_x
\wedge v_y + v_x \wedge (\d_i)_y \right) = \int_{\partial \O} v
\cdot \frac{\partial}{\partial \nu} \left( \frac{\partial
\d_i}{\partial p_i} \right) = O(\l_i^{-1}) \|v\|,
\end{equation}
since $\frac{\partial \d_i}{\partial p_i}$ satisfies
\eqref{eq:in5}. To estimate $\int_\O \frac{\partial \d_i}{\partial
p_i} \cdot \left( (\d_j)_x \wedge v_y + v_x \wedge (\d_j)_y
\right)$ for $j \neq i$, we proceed as follows. Let $B_i$ and
$B_j$ denote the balls of radius $\frac{1}{2 \ov{C}}$ centered at
$p_i$ and $p_j$ respectively. By the definition of $Z$ these two
balls are disjoint and moreover, by \eqref{eq:est2} $\n \d_i$ and
$\n \d_j$ are of order $\l_i^{-1}$ and $\l_j^{-1}$ respectively
outside $B_i$ and $B_j$. Hence we have
\begin{eqnarray}\label{eq:add4}
  \left| \int_\O \frac{\partial \d_i}{\partial
p_i} \cdot \left( (\d_j)_x \wedge v_y + v_x \wedge (\d_j)_y
\right) \right| \leq C \left( \frac{1}{\l_i} + \frac{1}{\l_j} +
\frac{1}{\l_i \l_j} \right) \|v\|.
\end{eqnarray}
Hence \eqref{eq:add1}-\eqref{eq:add4} imply $\left| I''_\e(z)
\left[ \frac{\partial z}{\partial p_i} \right] \right| \leq C
\left( \e + \sum_{j=1}^{k} \l_j^{-1} \right)$. The remaining part
of the statement follows from similar arguments.
\end{pf}

\noindent From Proposition \ref{p:ambbad}, critical points of
$I_\e$ restricted to $Z_\e$ are true critical points of $I_\e$. We
define $\tilde{I}_\e : Z \to \R$ as $\tilde{I}_\e(z) =
I_\e(z+w_\e(z))$. We now analyze the reduced functional
$\tilde{I}_\e$.

\begin{pro}\label{p:exp}
Let $\ov{C} > 0$, let $Z$ be as in \eqref{eq:manz} and let
$w_\e(z)$ be as in Proposition \ref{p:ambbad}. Then we have
\begin{equation}\label{eq:tieie}
  \left| \tilde{I}_\e(z) - I_\e(z) \right| \leq C \tilde{e}^2(\e, \l_1,
  \dots, \l_k), \qquad \forall z \in Z;
\end{equation}
Moreover, for all $z = \sum_{i=1}^k R_i P_{p_i,\l_i} \in Z$ there
holds
\begin{equation}\label{eq:gtieie}
  \begin{cases}
  \left| \frac{\partial \tilde{I}_\e(z)}{\partial p_i} -
  \frac{\partial I_\e(z)}{\partial p_i} \right| +
  \left| \frac{\partial \tilde{I}_\e(z)}{\partial R_i} -
  \frac{\partial I_\e(z)}{\partial R_i} \right| \leq
  \left( \e + \sum_j \l_j^{-1} \right) \tilde{e}(\e, \l_1,
  \dots, \l_k); & \\
  \left| \frac{\partial \tilde{I}_\e(z)}{\partial \l_i} -
  \frac{\partial I_\e(z)}{\partial \l_i} \right| \leq
  \frac{1}{\l_i} \left( \e + \sum_j \l_j^{-1} \right)
  \tilde{e}(\e, \l_1, \dots, \l_k)
  \end{cases}
\end{equation}
\end{pro}

\begin{pf}
We have
\begin{eqnarray*}
  \tilde{I}_\e(z) - I_\e(z) = I_\e(z+w) - I_\e(z) = \int_0^1
  I'_\e(z + s w)[w] \, d s = I'_\e(z) [w] + \int_0^1
  (I'_\e(z + s w) - I'_e(z))[w] \, d s.
\end{eqnarray*}
Since the functional $I''_\e$ is locally bounded, we have the
following estimate
$$
\left| I'_\e(z + s w) - I'_\e(z) \right| \leq \left| \int_0^1
  I''_\e(z + t s w)[w] dt \right| \leq \sup_{t,s \in [0,1]}
  \| I''_\e(z + t s w) \| \|w\| \leq C \|w\|
$$
for some fixed constant $C$ depending on $\O$,
$\|g\|_{C^2(\partial \O)}$ and the above constant $\ov{C}$. Using
the last three equations, Lemma \ref{l:gr} and the property $iii)$
in Proposition \ref{p:ambbad} we find
$$
\left| \tilde{I}_\e(z) - I_\e(z) \right| \leq \| I'_\e(z) \| \, \|
w_\e(z) \| + C \|w_\e(z)\|^2 \leq \tilde{e}^2(\e, \l_1,
  \dots, \l_k).
$$
This concludes the proof of \eqref{eq:tieie}. We just sketch the
proof of \eqref{eq:gtieie}. Differentiating the equation
$F_\e(z,w,q) = 0$ with respect to $p_i$ we obtain
$$
  0 = \frac{\partial F_\e}{\partial z} \frac{\partial z}{\partial
  p_i} + \frac{\partial F_\e}{\partial (w,q)} \frac{\partial
  (w,q)}{\partial p_i} = I''_\e(z+w_\e(z)) \frac{\partial z}{\partial
  p_i} + \frac{\partial F_\e}{\partial (w,q)} \frac{\partial
  (w,q)}{\partial p_i}.
$$
Similarly as before, one finds that
$\partial_{(w,q)} F_\e$ is uniformly invertible, and hence
\begin{equation}\label{eq:wp}
  \left\| \frac{\partial w}{\partial p_i} \right\| \leq C \left\|
I''_\e(z+w_\e(z)) \frac{\partial z}{\partial p_i} \right\| \leq C
\left\| I''_\e(z) \frac{\partial z}{\partial p_i} \right\| + C
\|w_\e(z)\| \left\| \frac{\partial z}{\partial p_i} \right\|,
\end{equation}
where we have used the fact that $I''_\e$ is locally Lipschitz. We
have
\begin{eqnarray}\label{eq:wp2}
  \frac{\partial \tilde{I}_\e(z)}{\partial p_i} - \frac{\partial
  I_\e(z)}{\partial p_i} & = & I''_\e(z) \left[
  \frac{\partial z}{\partial p_i}, w_\e(z) \right] + \int_0^1
  \left( I''_\e(z + s w_\e(z)) - I''_\e(z) \right) \left[
  \frac{\partial z}{\partial p_i}, w_\e(z) \right] \nonumber \\
  & + & I'_\e(z) \left[ \frac{\partial w}{\partial p_i} \right] +
  I''_\e(z) \left[ \frac{\partial w}{\partial p_i}, w_\e(z) \right]
  + \int_0^1 \left( I''_\e(z + s w_\e(z)) - I''_\e(z) \right) \left[
  \frac{\partial w}{\partial p_i}, w_\e(z) \right].
\end{eqnarray}
Equation \eqref{eq:wp2} implies
\begin{eqnarray*}
  \left| \frac{\partial \tilde{I}_\e(z)}{\partial p_i} - \frac{\partial
  I_\e(z)}{\partial p_i} \right| & \leq & C \left\| I''_\e(z) \left[
  \frac{\partial z}{\partial p_i} \right] \right\| \|w_\e(z)\| +
  C \|w_\e(z)\|^2 \left( \left\| \frac{\partial z}{\partial p_i} \right\|
  + \left\| \frac{\partial w}{\partial p_i} \right\| \right) \\ &
  + & C \left\| I'_\e(z) \right\| \left\| \frac{\partial w}{\partial p_i}
  \right\| + C \left\| w_\e(z) \right\| \left\|
  \frac{\partial w}{\partial p_i}  \right\|.
\end{eqnarray*}
Then the estimate of $\frac{\partial \tilde{I}_\e(z)}{\partial
p_i} - \frac{\partial I_\e(z)}{\partial p_i}$ in \eqref{eq:gtieie}
follows from Lemma \ref{l:gr}, \eqref{eq:wp} and Lemma
\ref{l:i2zp}. The remaining part of \eqref{eq:gtieie} follows from
similar estimates.
\end{pf}

\section{The expansion for one bubble}\label{s:1b}

In this section we compute the expansion of $I_\e(z)$, with $z \in
Z$, for $\e$ small and in the case $k = 1$. This is essentially
performed in \cite{isc}-\cite{isa1}, in order to construct
blowing-up solutions of \eqref{eq:in2}, and in order to
characterize the mountain-pass solutions in the limit $\e \to 0$.
We derive the expansion here, in a form which is useful for us in
the expansion for multiple bubbles in Section \ref{s:2b}. Let us
first introduce some notation. We recall that $g : \O \to \R^3$
denotes the solution of \eqref{eq:ogg}, and letting $R \in SO(3)$,
define $d_R g : \O \to \R$ by
$$
  d_R g (a) = \frac{\partial}{\partial x} (R \circ g)_1 (a) +
  \frac{\partial}{\partial y} (R \circ g)_2 (a), \qquad x \in \O.
$$
%
For a fixed boundary datum $\tilde{g}$, we are interested in
expanding the functional $I_\e(P \d)$ as a function of the
parameters $a$, $\l$, $R$ and $\e$. We have the following
Proposition.
\begin{pro}\label{p:}
Let $\ov{C} > 0$ be fixed, and let $a, \l, R$ be such that $P R
\d_{a,\l} \in Z$. Then, setting
$$
  A_0 = \int_{\R^2} \frac{|\xi|^2}{(1+|\xi|^2)^3};  \qquad
  F_{\O,g} (\e,a,\l,R) = 8 A_0 \left( \frac{1}{\l^2}
\tilde{H}(a) - \frac{\e}{\l} d_{R^{-1}}g (a) \right),
$$
there holds
$$
I_\e(P R \d_{a,\l}) = \frac{8}{9} A_0 + F_{\O,g} (a,\l,R) +
o(\e^2) + e(\e,\l);
$$
$$
\frac{\partial I_\e(P  R \d_{a,\l})}{\partial a} = \frac{\partial
F_{\O,g}}{\partial a} + e(\e,\l); \qquad \frac{\partial I_\e(P  R
\d_{a,\l})}{\partial \l} = \frac{\partial F_{\O,g}}{\partial \l} +
\frac{e(\e,\l)}{\l}; \qquad \frac{\partial I_\e(P R
\d_{a,\l})}{\partial R} = \frac{\partial F_{\O,g}}{\partial R} +
e(\e,\l),
$$
where $e(\e,\l)$ is defined in Section \ref{s:pre}.
\end{pro}

\begin{pf}
We assume that $R = Id$, and we write $\d$ for $\d_{a,\l}$. Let
also $\var$ be the solution of \eqref{eq:fi}. We have
\begin{eqnarray}\label{eq:exp1}
I_\e(P \d) & = & \frac{1}{2} \int_{\O} |\n \d|^2 + \frac{1}{2}
\int_{\O} |\n \varphi|^2 - \IO \n \var \cdot \n \d + \frac{2}{3}
\IO (\d - \var) \cdot \left( (\d - \var)_x \wedge (\d - \var)_y
\right) \nonumber
\\ & + & \e \IO (\d - \var) \cdot \left( (\d - \var)_x \wedge g_y + g_x
\wedge (\d - \var)_y \right) + 2 \e^2 \IO (\d - \var) \cdot (g_x
\wedge g_y).
\end{eqnarray}
Integrating by parts we can write
\begin{equation}\label{eq:fd}
\frac{1}{2} \int_{\O} |\n \d|^2 + \frac{1}{2} \int_{\O} |\n
\varphi|^2 - \IO \n \var \cdot \n \d = \IO (\var -\d) \cdot (\d_x
\wedge \d_y).
\end{equation}
We expand first the expression in \eqref{eq:fd}. Let us evaluate
the $z$-component in the scalar product of the integral on the
right-hand side in \eqref{eq:fd}. From formulas \eqref{eq:dxy} and
\eqref{eq:estdf} we deduce
\begin{eqnarray*}
(\var - \d)_3 (\d_x \wedge \d_y)_3 & = & - \left( 1 -
\frac{2}{\l^2} h_3(a,a) - \frac{\l^2 |\xi-a|^2 - 1}{\l^2 |\xi-a|^2
+ 1} + o(\l^{-2}) \right) 4 \l^2 \frac{1 - \l^2 |\xi-a|^2}{(\l^2
|\xi-a|^2 + 1)^3} \\
& = & - 8 \l^2 \left( \frac{1}{\l^2 |\xi-a|^2 + 1} -
\frac{1}{\l^2} h_3(a,a) + o(\l^{-2}) \right) \frac{1 - \l^2
|\xi-a|^2}{(\l^2 |\xi-a|^2 + 1)^3}.
\end{eqnarray*}
Integrating on $\O$ we get
\begin{eqnarray}\label{eq:43'}
\IO (\var - \d)_3 (\d_x \wedge \d_y)_3 = - 8 \l^2 \IO \left(
\frac{1}{\l^2 |\xi-a|^2 + 1} - \frac{1}{\l^2} h_3(a,a) +
o(\l^{-3}) \right) \frac{1 - \l^2 |\xi-a|^2}{(\l^2 |\xi-a|^2 +
1)^3}.
\end{eqnarray}
Using a change of variable we obtain
\begin{eqnarray}\label{eq:int1}
\l^2 \IO \frac{1 - \l^2 |\xi-a|^2}{(\l^2 |\xi-a|^2 + 1)^4} & = &
\l^2 \int_{\R^2} \frac{1 - \l^2 |\xi-a|^2}{(\l^2 |\xi-a|^2 + 1)^4}
- \l^2 \int_{\R^2 \setminus \O} \frac{1 - \l^2 |\xi-a|^2}{(\l^2
|\xi-a|^2 + 1)^4} \nonumber \\ & = & \int_{\R^2} \frac{1 -
|\xi|^2}{(|\xi|^2 + 1)^4} + O(\l^{-4}),
\end{eqnarray}
and also
\begin{eqnarray}\label{eq:int2}
\IO \frac{1 - \l^2 |\xi-a|^2}{(\l^2 |\xi-a|^2 + 1)^3} & = &
\int_{\R^2} \frac{1 - \l^2 |\xi-a|^2}{(\l^2 |\xi-a|^2 + 1)^3} -
\int_{\R^2 \setminus \O} \frac{1 - \l^2 |\xi-a|^2}{(\l^2 |\xi-a|^2
+ 1)^3} \nonumber \\ & = & \frac{1}{\l^2} \int_{\R^2} \frac{1 -
|\xi|^2}{(|\xi|^2 + 1)^3} + O(\l^{-4}) = O(\l^{-4}).
\end{eqnarray}
The last identity follows from \eqref{eq:id}. In conclusion, from
\eqref{eq:43'}, \eqref{eq:int1} and \eqref{eq:int2} we get
\begin{equation}\label{eq:intz}
  \IO (\var - \d)_3 (\d_x \wedge \d_y)_3 = 8 \int_{\R^2}
  \frac{|\xi|^2 - 1}{(|\xi|^2 + 1)^4} + o(\l^{-2}).
\end{equation}
We consider now the $x$ and $y$ components of the integral on the
right-hand side in \eqref{eq:fd}. We have, using \eqref{eq:apd}
and \eqref{eq:dxdyxy}
\begin{eqnarray*}
\IO (\var - \d)_1 (\d_x \wedge \d_y)_1 = - \IO \left( \frac{2}{\l}
h_1(\xi,a) - \frac{2 \l (x - a_1)}{1 + \l^2 |\xi-a|^2} +
o(\l^{-1}) \right) \frac{8 \l^3 (x-a_1)}{(1 + \l^2 |\xi-a|^2)^3}.
\end{eqnarray*}
We can write
\begin{eqnarray*}
\l^4 \IO \frac{(x - a_1)^2}{(1 + \l^2 |\xi-a|^2)^4}  =  \l^4
\int_{\R^2} \frac{(x - a_1)^2}{(1 + \l^2 |\xi-a|^2)^4} - \l^4
\int_{\R^2 \setminus \O} \frac{(x - a_1)^2}{(1 + \l^2
|\xi-a|^2)^4}
\\ = \int_{\R^2} \frac{x^2}{(1 + |\xi|^2)^4} + O(\l^{-4}).
\end{eqnarray*}
From the smoothness of $h_i$ we have also
\begin{eqnarray*}
\left| \IO \left(h_1(\xi,a) - h_1(a,a) - (\xi-a) \cdot \n h_1
(a,a) \right) \frac{x-a_1}{(1 + \l^2 |\xi-a|^2)^3} \right| \leq
O\left( \IO \frac{|\xi-a|^3}{(1 + \l^2 |\xi-a|^2)^3} \right).
\end{eqnarray*}
As a consequence we deduce
\begin{eqnarray*}
\IO h_1(x,a) \frac{(x-a_1)}{(1 + \l^2 |\xi-a|^2)^3} & = & \IO
h_1(a,a) \frac{(x-a_1)}{(1 + \l^2 |\xi-a|^2)^3} \\ & + & \IO
((\xi-a) \cdot \n h_1(a,a)) \frac{(x-a_1)}{(1 + \l^2 |\xi-a|^2)^3} \\
& + & O\left( \IO \frac{|\xi-a|^3}{(1 + \l^2 |\xi-a|^2)^3} \right)
= \frac{1}{\l^4} \, \partial_x h_1(a,a) \int_{\R^2} \frac{x^2}{(1
+ |\xi|^2)^3} + o(\l^{-4}).
\end{eqnarray*}
From the last equation we deduce
\begin{equation}\label{eq:intx}
\IO (\var - \d)_1 (\d_x \wedge \d_y)_1 = 16 \int_{\R^2}
\frac{x^2}{(1 + |\xi|^2)^4} - 16 \frac{1}{\l^2} \, \partial_x
h_1(a,a) \int_{\R^2} \frac{x^2}{(1 + |\xi|^2)^3} + o(\l^{-2}).
\end{equation}
In the same way we obtain
\begin{equation}\label{eq:inty}
\IO (\var - \d)_2 (\d_x \wedge \d_y)_2 = 16 \int_{\R^2}
\frac{x^2}{(1 + |\xi|^2)^4} - 16 \frac{1}{\l^2} \, \partial_y
h_3(a,a) \int_{\R^2} \frac{x^2}{(1 + |\xi|^2)^3} + o(\l^{-2}).
\end{equation}
From \eqref{eq:intz}, \eqref{eq:intx} and \eqref{eq:inty} it
follows that
\begin{eqnarray*}
  \IO (\var - \d) \cdot (\d_x \wedge \d_y) & = & 8 \int_{\R^2}
  \frac{|\xi|^2 - 1}{(1 + |\xi|^2)^4} + 32 \int_{\R^2} \frac{x^2}{(1 + |\xi|^2)^4}
  \nonumber \\ & - & 16 \, \frac{1}{\l^2}
  \left( \int_{\R^2} \frac{x^2}{(1 + |\xi|^2)^3} \right)
  \left( \frac{\partial h_1}{\partial x} +
  \frac{\partial h_2}{\partial y} \right)(a,a) + o(\l^{-2}).
\end{eqnarray*}
It is standard to check that
$$
8 \int_{\R^2} \frac{|\xi|^2 - 1}{(1 + |\xi|^2)^4} + 32 \int_{\R^2}
\frac{x^2}{(1 + |\xi|^2)^4} = \frac{8}{3} A_0,
$$
Hence from the last two equations we find
\begin{eqnarray}\label{eq:intf1}
  \IO (\var - \d) \cdot (\d_x \wedge \d_y) = \frac{8}{3} A_0 -
  \frac{8}{\l^2} A_0 \left( \frac{\partial h_1}{\partial x} +
  \frac{\partial h_2}{\partial y} \right)(a,a) + o(\l^{-2}).
\end{eqnarray}
We turn now to the fourth term in \eqref{eq:exp1}. We have clearly
$$
\IO (\d - \var) \cdot \left( (\d - \var)_x \wedge (\d - \var)_y
\right) = \IO (\d - \var) \cdot \left( \d_x \wedge \d_y \right) -
\IO (\d - \var) \cdot ( \d_x \wedge \var_y + \var_x \wedge \d_y ).
$$
Let us consider the term $\d_x \wedge \var_y$. Using the above
formulae we deduce
\begin{eqnarray}\label{eq:star}
  (\d_x \wedge \var_y)_1 = \frac{8 \l (x-a_1)(y-a_2)}{(1 + \l^2
|\xi-a|^2)^2} \left( \frac{\partial h_3}{\partial y} + O(\l^{-1})
\right) - \frac{8 \l (x-a_1)}{(1 + \l^2 |\xi-a|^2)^2} \left(
\frac{\partial h_2}{\partial y} + O(\l^{-1}) \right); \nonumber \\
(\d_x \wedge \var_y)_2 = \frac{8 \l (x-a_1)}{(1 + \l^2
|\xi-a|^2)^2} \left( \frac{\partial h_1}{\partial y} + O(\l^{-1})
\right) + \frac{4}{\l} \frac{1 + \l^2 (y^2-x^2)}{(1 + \l^2
|\xi-a|^2)^2} \left( \frac{\partial h_3}{\partial y} + O(\l^{-1})
\right); \\
(\d_x \wedge \var_y)_3 = 4 \frac{1 + \l^2
((y-a_2)^2-(x-a_1)^2)}{(1 + \l^2 |\xi-a|^2)^2} \left(
\frac{\partial h_2}{\partial y} + O(\l^{-1}) \right) + \frac{8
\l^2 (x-a_1)(y-a_2)}{(1 + \l^2 |\xi-a|^2)^2} \left( \frac{\partial
h_1}{\partial y} + O(\l^{-1}) \right). \nonumber
\end{eqnarray}
and also
\begin{eqnarray}\label{eq:dstar}
  (\var_x \wedge \d_y)_1 = \frac{8 \l (y-a_2)}{(1 + \l^2
|\xi-a|^2)^2} \left( \frac{\partial h_2}{\partial x} + O(\l^{-1})
\right) + \frac{4}{\l} \frac{1 + \l^2((x-a_1)^2- (y-a_2)^2)}{(1 +
\l^2 |\xi-a|^2)^2} \left( \frac{\partial h_3}{\partial x} +
O(\l^{-1}) \right); \nonumber \\
(\var_x \wedge \d_y)_2 = \frac{8 \l (x-a_1)(y-a_2)}{(1 + \l^2
|\xi-a|^2)^2} \left( \frac{\partial h_3}{\partial x} + O(\l^{-1})
\right) - \frac{8 \l (y-a_2)}{(1 + \l^2 |\xi-a|^2)^2} \left(
\frac{\partial h_1}{\partial x} + O(\l^{-1}) \right); \\
(\var_x \wedge \d_y)_3 = 4 \frac{1 + \l^2
((x-a_1)^2-(y-a_2)^2)}{(1 + \l^2 |\xi-a|^2)^2} \left(
\frac{\partial h_1}{\partial x} + O(\l^{-1}) \right) + \frac{8
\l^2 (x-a_1)(y-a_2)}{(1 + \l^2 |\xi-a|^2)^2} \left( \frac{\partial
h_2}{\partial x} + O(\l^{-1}) \right). \nonumber
\end{eqnarray}
The functions $\frac{\partial h_i}{\partial x_j}$ in the last
formulas, as before, are evaluated at the point $(a,a)$. Using
\eqref{eq:estdf}, \eqref{eq:star} and \eqref{eq:dstar} we find
\begin{equation}\label{eq:m2}
\IO (\d - \var) \cdot ( \d_x \wedge \var_y) \sim - 32
\frac{1}{\l^2} \left( \int_{\R^2} \frac{x^2}{(1 + |\xi|^2)^3}
\right) \frac{\partial h_2}{\partial y}(a,a) + o(\l^{-2}),
\end{equation}
and similarly
$$
\IO (\d - \var) \cdot ( \var_x \wedge \d_y) \sim - 32
\frac{1}{\l^2} \left( \int_{\R^2} \frac{x^2}{(1 + |\xi|^2)^3}
\right) \frac{\partial h_1}{\partial x}(a,a) + o(\l^{-2}),
$$
Since $\left|(\var_x \wedge \var_y)_1\right|, \left|(\var_x \wedge
\var_y)_1\right| \leq O(\l^{-3})$ and $\left|(\var_x \wedge
\var_y)_3\right| \leq O(\l^{-2})$, one can check that
$$
  \IO (\d - \var) \cdot (\var_x \wedge \var_y) = o(\l^{-2}).
$$
Let us now turn to the fifth term in
\eqref{eq:exp1}. The quantity $(\d - \var) \cdot \left( (\d -
\var)_x \wedge g_y \right)$ can be estimated as
\begin{eqnarray}\label{eq:e1}
 &  &
\left( \frac{2 \l (x-a_1)}{1 + \l^2 |\xi-a|^2} - \frac{2}{\l} h_1
\right) \left[ \left( - \frac{4 \l^3(x-a_1)(y-a_2)}{(1 + \l^2
|\xi-a|^2)^2} - \frac{2}{\l} \partial_x h_2 \right) (g_3)_y
\right. \nonumber \\ & - & \left. \left( \frac{4\l^2(x-a_1)}{(1 +
\l^2 |\xi-a|^2)^2} + \frac{2}{\l^2} \partial_x h_3 \right) (g_2)_y
\right] \nonumber \\
& + & \left( \frac{2 \l (y-a_2)}{1 + \l^2 |\xi-a|^2} -
\frac{2}{\l} h_2 \right) \left[ \left( \frac{4 \l^2 (x-a_1)}{(1 +
\l^2 |\xi-a|^2)^2} + \frac{2}{\l^2} \partial_x h_3 \right) (g_1)_y
\right. \nonumber
\end{eqnarray}
\begin{eqnarray}\label{eq:ee1}
& - & \left. \left( 2 \l \frac{1 + \l^2((y-a_2)^2 - (x-a_1)^2)}{(1
+ \l^2 |\xi-a|^2)^2} - \frac{2}{\l}
\partial_x h_1 \right) (g_3)_y \right] \nonumber
\\ & + & \left( \frac{2}{\l^2} h_3 - \frac{2}{1 + \l^2 |\xi-a|^2} \right)
\left[ \left( 2 \l \frac{1 + \l^2((y-a_2)^2-(x-a_1)^2)}{(1 + \l^2
|\xi-a|^2)^2} - \frac{2}{\l}
\partial_x h_1 \right) (g_2)_y \right. \nonumber \\ & + & \left. \left(
\frac{4\l^3(x-a_1)(y-a_2)}{(1 + \l^2 |\xi-a|^2)^2} + \frac{2}{\l}
\partial_x h_2 \right) (g_1)_y \right] \nonumber.
\end{eqnarray}
Integrating on $\O$ and reasoning as before we get
\begin{equation}\label{eq:ef1}
  \IO (\d - \var) \cdot \left( (\d - \var)_x \wedge g_y \right)
  = - \frac{16}{\l} \int_{\R^2} \frac{x^2}{(1 + |\xi|^2)^3} \,
  (g_2)_y + O(\l^{-2}),
\end{equation}
and similarly
\begin{equation}\label{eq:ef2}
  \IO (\d - \var) \cdot \left( g_x \wedge (\d - \var)_y \right)
  = - \frac{16}{\l} \int_{\R^2} \frac{x^2}{(1 + |\xi|^2)^3} \,
  (g_1)_x + O(\l^{-2}).
\end{equation}
Finally, the last term in \eqref{eq:exp1} is easily seen to be of
order $o(\e^2)$. This concludes the proof in the case of $R = Id$.
For a generic rotation $R$ it is sufficient, by invariance, to
consider the boundary datum $R^{-1} g$ and to substitute $(g_1)_x
+ (g_2)_y$ with $d_{R^{-1}} g$.
\end{pf}

\begin{rem}\label{r:is}
Propositions \ref{p:exp} and \ref{p:} allow us to find critical
points of $I_\e$ extremizing the reduced functional $\tilde{I}_\e$
on $Z$. Differentiating with respect to $R$, $\l$, $a$ we get
\begin{equation}\label{eq:d2}
  \frac{\partial}{\partial R} d_{R^{-1}} g = 0;
  \quad 2 \tilde{H}(a) - \e \, \l d_{R^{-1}} g = 0; \quad
  \frac{1}{\l^2} \n \tilde{H}(a) - \frac{\e}{\l} \n d_{R^{-1}} g = 0.
\end{equation}
Using the second and third equations in \eqref{eq:d2} we deduce
$$
  \n \log \tilde{H}(a) = 2 \n \log d_{R^{-1}} g (a).
$$
The extremization with respect to $R$ (the first
equation in \eqref{eq:d2}) is performed in \cite{isa1} Lemma 5.4
and \cite{isc} Lemma 3.1.2 and, requiring $d_{R^{-1}} g$ to be
positive (from the second equation in \eqref{eq:d2}) yields
\begin{equation}\label{eq:d4}
  \frac{\partial}{\partial R} d_{R^{-1}} g = 0 \qquad \Rightarrow
  \qquad d_{R^{-1}} g = \left( |\n g|^2 \pm 2 | g_x \wedge g_y |
  \right)^{\frac{1}{2}}.
\end{equation}
Hence, under the conditions $\n g \neq 0$, the extremization in
\eqref{eq:d2} becomes
\begin{equation}\label{eq:b1fin}
  \l = \frac{2}{\e} \frac{\tilde{H}(a)}{d_{R^{-1}} g(a)};
  \qquad \frac{\partial}{\partial R} {\rm div } \pi_R g = 0;
  \qquad \n \frac{|\n g|^2 \pm 2 | g_x \wedge
  g_y |}{\tilde{H}}(a) = 0.
\end{equation}
In particular the mountain-pass solution of \eqref{eq:inge} (see
\cite{isa1}) has minimal energy on $Z$, and one has to chose the
$+$ sign in \eqref{eq:d4}-\eqref{eq:b1fin}, and the last function
in \eqref{eq:b1fin} is maximized on $\O$.
\end{rem}

\begin{rem}\label{r:pos}
(a) We point out that the expansions in \eqref{eq:fd} and
\eqref{eq:intf1} yield
$$
\|P \d\|^2_{H^1(\O)} = \|\d\|^2_{\mathcal{D}} - \frac{16
A_0}{\l^2} \tilde{H}(a) + o(\l^{-2}).
$$
Since the norm of $P \d$, being the projection of $\d$, is smaller
than the norm of $\d$, the above formula implies $\tilde{H} \geq
0$ on $\O$. A little more calculation shows that indeed $\tilde{H}
> 0$, as proved in \cite{isa1}.

(b) The functions $h_1$, $h_2$, and $h_3$ defined in
\eqref{eq:h1h2}, \eqref{eq:h3} are related to the boundary values
of $\d$. Since $h_1$ and $h_2$ are of order $\l^{-1}$, while $h_2$
is of order $\l^{-2}$, $h_3$ appears in the expansion only as a
lower order term.
\end{rem}

\section{The role of the Robin function}\label{s:robin}

In this Section we investigate the relation between the Robin
function $H$ and the function $\tilde{H}$ defined in
\eqref{eq:hhi2}. Since $\tilde{H}$ consists of second derivatives
of the regular part of the Green's function, while the Robin
function involves the regular part itself, we need to use global
arguments, based on the Riemann mapping theorem. See \cite{baf}
for some properties of the Robin function.

\subsection{Simply connected domains}

In this subsection we prove the first assertion of Theorem
\ref{t:in1}.

\begin{pro}\label{p:hth}
Let $\O \subseteq \R^2$ be a smooth simply connected domain. Then,
if $f : \O \to D$ is a Riemann map, there holds
\begin{equation}\label{eq:hth}
  e^{2 H(a,a)} = \frac{|f'(a)|^2}{\left( 1 - |f'(a)|^2 \right)^2};
  \qquad \tilde{H}(a) = 2 \frac{|f'(a)|^2}{\left( 1 - |f'(a)|^2
  \right)^2}.
\end{equation}
\end{pro}

\begin{pf}
The first part of the statement is well-known, see e.g.
\cite{baf}, Table 2. Letting $G_D(a,\xi)$ denote the Green's
function for $D$, and letting $\var : D \to \R$ being any smooth
function with compact support, we have
$$
\var(a') = \int_{D} G_D(a',\xi) \D \var(\xi) d \xi, \qquad a' \in
D.
$$
Let $\psi : \O \to \R$ be defined by $\psi = \var \circ f$. We
have $\D \psi (\xi) = \frac{1}{|f'(\xi)|^2} \D \var (\xi)$ and
hence, letting $G(a,\xi)$ be the Green's function for $\O$ and
using a change of variables we get
$$
\psi (a) = \var(f(a)) = \var(a') = \int_D G_D(a',\xi) \D \var(\xi)
d \xi = \int_{\O} G_D(f(a),f(\xi)) \D \psi(\xi) d \xi,
$$
where $a' = f(a)$. Hence, from the explicit expression of $G_D$ it
turns out that
$$
G(a,z) = G_D(f(a),f(z)) = - \frac{1}{2} \log \frac{|f(z) -
f(a)|^2}{|1 - f(z) \ov{f(a)}|^2}; \qquad a, z \in \O,
$$
where we have identified $\xi$ with the point $z$ in the complex
plane. It follows that
$$
H(a,z) = \frac 1 2 \left[ \log |z - a|^2 - \log \frac{|f(z) -
f(a)|^2}{|1 - f(z) \ov{f(a)}|^2} \right]; \qquad a, z \in \O.
$$
The last expression can be rewritten as
$$
  H(a,z) = \frac{1}{2} \left[ \log |1 - f(z) \ov{f(a)}|^2 -
  \log \frac{|f(z) - f(a)|^2}{|z - a|^2} \right];
  \qquad a, z \in \O.
$$
In particular, taking the limit $z \to a$, we
deduce immediately the first equality in \eqref{eq:hth}.

Using complex notation, we have
$$
  \frac{\partial}{\partial a_1} = \left( \frac{\partial}{\partial a} +
  \frac{\partial}{\partial \ov{a}} \right); \quad
  \frac{\partial}{\partial a_2} = i \left( \frac{\partial}{\partial a}
  - \frac{\partial}{\partial \ov{a}} \right); \qquad
  \frac{\partial}{\partial x} = \left( \frac{\partial}{\partial z} +
  \frac{\partial}{\partial \ov{z}} \right); \quad
  \frac{\partial}{\partial y} = i \left( \frac{\partial}{\partial z}
  - \frac{\partial}{\partial \ov{z}} \right).
$$
It follows that
\begin{equation}\label{eq:defl}
  \frac{\partial}{\partial x} \frac{\partial}{\partial a_1} +
\frac{\partial}{\partial y} \frac{\partial}{\partial a_2} = 2
\frac{\partial}{\partial z} \frac{\partial}{\partial \ov{a}} + 2
\frac{\partial}{\partial \ov{z}} \frac{\partial}{\partial a} = 4
{\bf Re} \frac{\partial}{\partial \ov{z}} \frac{\partial}{\partial
a} := L.
\end{equation}
To derive the expression of $\tilde{H}$,
recall \eqref{eq:hhi}, we apply $L$ to $H(a,z)$ and evaluate at $z
= a$. We have, still in complex notation
$$
\frac{\partial}{\partial a} H(a,z) = - \frac{1}{2} \frac{\ov{f(z)}
f'(a)}{1 - \ov{f(z)} f(a)} - \frac{1}{2} \frac{\partial}{\partial
a} \log \frac{f(z)-f(a)}{z-a}.
$$
When we apply the operator $\frac{\partial}{\partial \ov{z}}$ the
second term vanishes and we get
$$
\frac{\partial}{\partial \ov{z}} \frac{\partial}{\partial a}
H(a,z) = - \frac{1}{2} \frac{f'(a) \ov{f'(z)}}{(1 - \ov{f(z)}
f(a))^2}.
$$
Taking the real part we find
$$
4 {\bf Re} \frac{\partial}{\partial \ov{z}}
\frac{\partial}{\partial a} H(a,z) = - \left[
\frac{f'(a)\ov{f'(z)}}{(1 - \ov{f(z)}f(a))^2} + \frac{\ov{f'(a)}
f'(z)}{(1 - \ov{f(a)}f(z))^2} \right]
$$
Choosing $z = a$ in the last formula, we obtain the second
identity in \eqref{eq:hth}. This concludes the proof.
\end{pf}

\begin{rem}\label{r:ic}
(a) The expression of $\tilde{H}$ in \eqref{eq:hth} does not
depend on the choice of the conformal map $f$ of $\O$ into $D$.

(b) From the explicit description of $\tilde{H}$ in Proposition
\ref{p:hth} we obtain $\tilde{H}(a) \to + \infty$ as $a \to
\partial \O$. This is true for any domain, as proved in
\cite{isa1}, Lemma 5.7.

(c) In the case of simply connected domains, the function
$\tilde{H}$ coincides with the square of the reciprocal of the
conformal radius and the hyperbolic radius, see \cite{baf},
Definitions 1, 7 and Theorem 8. See also Remark \ref{r:hhr}.

(e) Since every convex domain has a single conformal incenter, see
\cite{fl} Proposition 11, it follows that $\tilde{H}$ possesses a
unique critical point in this case. For a general simply connected
domain $\tilde{H}$ will have multiple critical points, see
\cite{fl} page 483. We also point out that, even if a conformal
transformation of the domain affects the number of critical points
of $\tilde{H}$, the topology of critical points at infinity (see
\cite{bab}) at the first level of non-compactness should be an
invariant.
\end{rem}

\subsection{Multiply connected domains}

In this subsection we derive a general formula for $\tilde{H}$ on
multiply connected domains. This formula makes use of the covering
map and the deck transformation.

Let us recall that the Green's function in the unit disk with pole
$z_0 \in D$ is given by
$$
G_D(z,z_0) = - \log \left| \frac{z - z_0}{1 - \ov{z}_0 z} \right|,
\qquad z, z_0 \in D.
$$
Let us pick a point $w \in \O$ and consider the Green's function
for $\O$ with pole at $a$. From \cite{baf}, Theorem 4, one has
$$
G_\O(a,w) = - \sum_k \log \left| \frac{z_k - z}{1 - \ov{z}_k z}
\right|, \qquad \hbox{ where } f(z) = w, f(z_k) = a \hbox{ for all
} k.
$$
Thus the regular part $G_\O(a,w)$ is
\begin{equation}\label{eq:mc1}
  H_\O(a,w) = \log |f(z_0) - f(z)| - \sum_{k} \log |z_k - z| + \sum_k
  \log |1 - \ov{z}_k z|.
\end{equation}
Now, as in the previous subsection, it is sufficient to apply the
operator $L$ defined in \eqref{eq:defl}. The first two terms
vanish when $L$ is applied. To handle the third term, note that
$f^{-1}$ is a {\em local} diffeomorphism. So $f^{-1}$ is defined
from $f^{-1} : \mathcal{U} \to \mathcal{V}_k$, where $\mathcal{U}$
is a neighborhood of $a$ and $\mathcal{V}_k$ is a neighborhood of
$z_k$. Thus we have
$$
\frac{\partial H}{\partial a} = \frac{\partial H}{\partial z_k}
\frac{\partial z_k}{\partial a} = \frac{1}{f'(z_k)} \frac{\partial
H}{\partial z_k}.
$$
Similarly, there holds
\begin{equation}\label{eq:mc2}
  \frac{\partial^2 H}{\partial \ov{w} \partial a} =
  \frac{1}{f'(z_k)} \frac{1}{\ov{f'(z)}}
  \frac{\partial^2 H}{\partial \ov{z} \partial z_k}.
\end{equation}
Hence, using \eqref{eq:mc1} and \eqref{eq:mc2} we find
$$
  4 {\bf Re} \frac{\partial^2 H}{\partial \ov{w} \partial a}|_{w=a}
  = - 2 {\bf Re} \sum_k \frac{1}{f'(z_k) \ov{f'(z_0)}} \frac{1}{(1 -
  z_k \ov{z}_0)^2}.
$$
Let $T_k$ be the Mobius (deck) transformation that maps $z_0$ into
$z_k$. We have
$$
f = f \circ T_k \qquad \Rightarrow \qquad f'(z_0) = f(z_k)
T'_k(z_0).
$$
Using the last equation and factoring the term $(1-|z_0|^2)^2$, we
deduce
\begin{equation}\label{eq:ffin}
  2 {\bf Re} \frac{\partial^2 H}{\partial \ov{w} \partial a}|_{w=w_1}
  = - 2 \frac{1}{|f'(z_0)|^2(1-|z_0|^2)^2} {\bf Re} \sum_k
  \frac{T'_k(z_0)(1-|z_0|^2)^2}{(1 - z_k \ov{z}_0)^2}.
\end{equation}
From \eqref{eq:ffin} we obtain immediately the following result.

\begin{pro}\label{p:zk}
Let $\O \subseteq \R^2$ be a multiply connected domain, and let $f
: D \to \O$ be a conformal covering map. Given $a \in \O$, let
$\{z_k\}_k$ be the pre-image of the point $a$ under the map $f$,
and let $T_k : D \to D$ denote the deck transformation mapping
$z_0$ into $z_k$. Then there holds
\begin{equation}\label{eq:ffin2}
  \tilde{H}(a) =
  \frac{1}{|f'(z_0)|^2(1-|z_0|^2)^2} \sum_k \left(
  \frac{T'_k(z_0)(1-|z_0|^2)^2}{(1 - z_k \ov{z}_0)^2} +
  \frac{\ov{T'_k(z_0)}(1-|z_0|^2)^2}{(1 - \ov{z_k} z_0)^2} \right).
\end{equation}
\end{pro}

\

\noindent We note that when $f^{-1}(a) = \{z_0\}$ we recover the
formula for the simply connected domain.

\subsection{Some numerical computation}

In this subsection we prove that in general, for a multiply
connected domain, the two functions $\tilde{H}$ and $2e^{2H}$ do
not coincide (this is the case for simply connected domains, see
Proposition \ref{p:hth}). We consider in particular the case of an
annulus of inner radius $\frac{1}{\rh}$ and outer radius $\rh$,
where $\rh > 1$. Our numerical computations show that the critical
points of these two functions do not coincide, hence we obtain the
statement $(b)$ in Theorem \ref{t:in1}.

\begin{figure}
\centerline{\rotatebox{270}{\psfig{figure=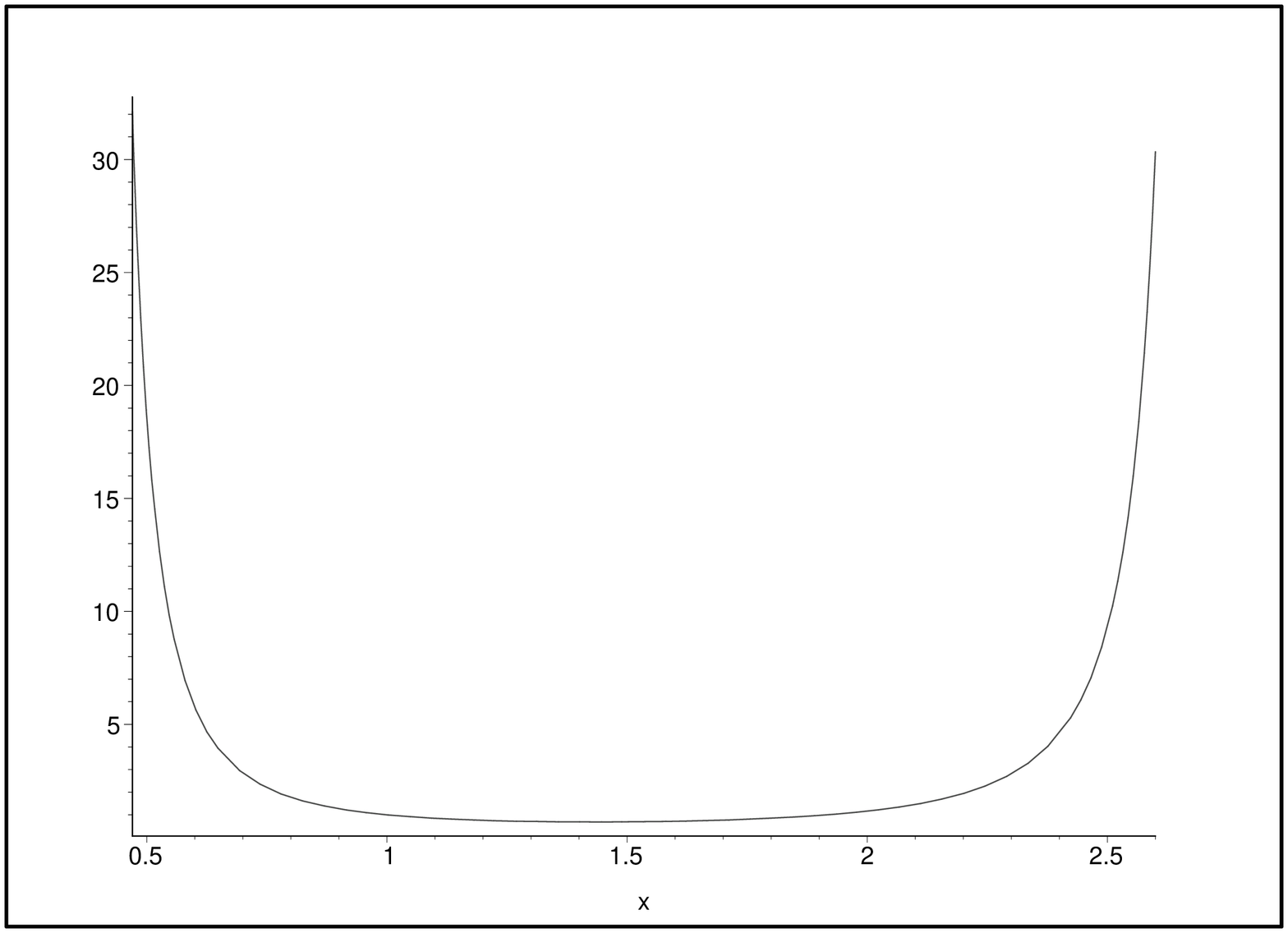,width=5.5cm}}
\quad \rotatebox{270}{\psfig{figure=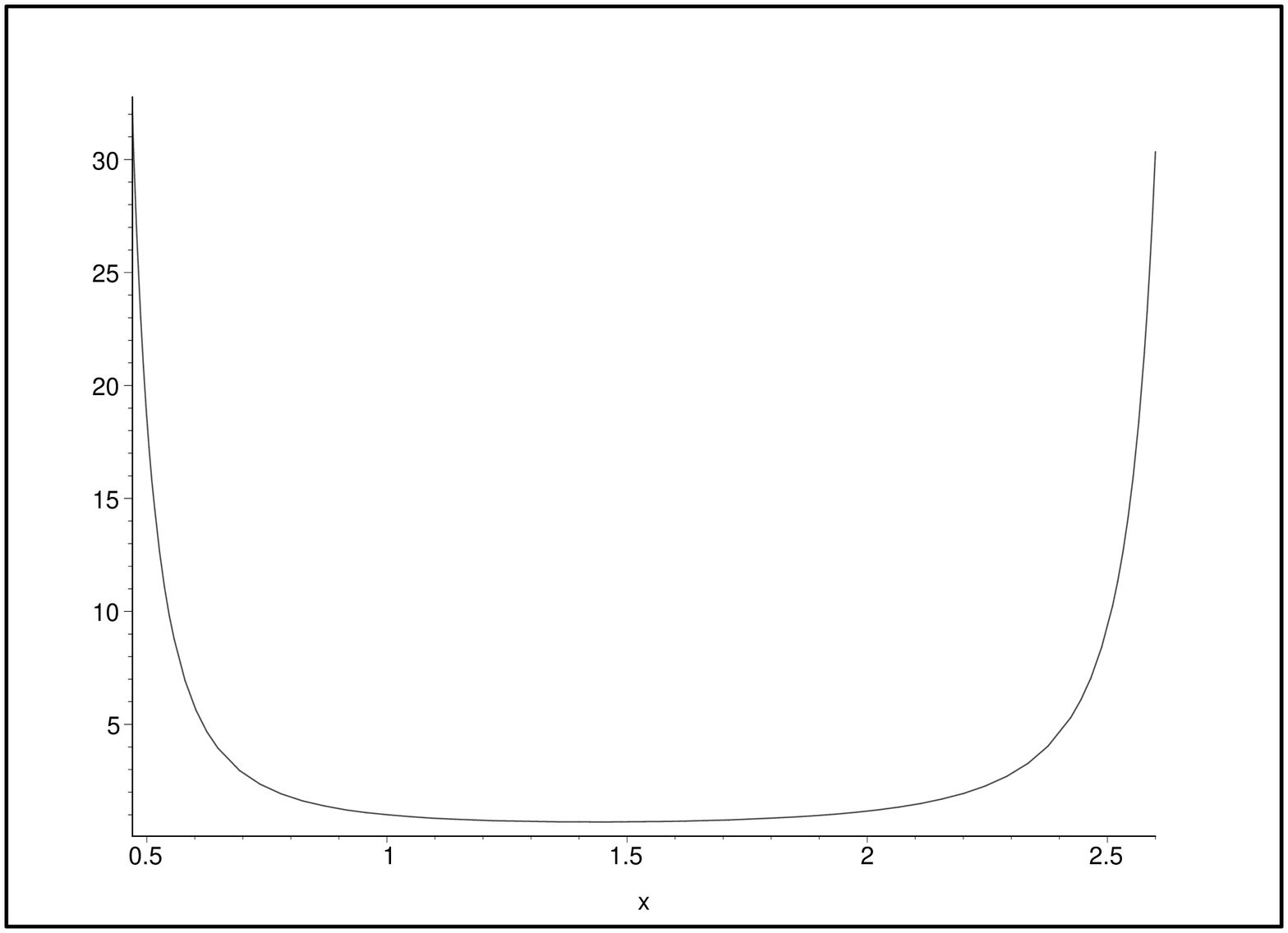,width=5.5cm}}}

\caption{the functions $2 e^{2H}$ and $\tilde{H}$ for $\rh = e$,
full picture} \label{fig2}
\end{figure}

\begin{figure}
\centerline{\rotatebox{270}{\psfig{figure=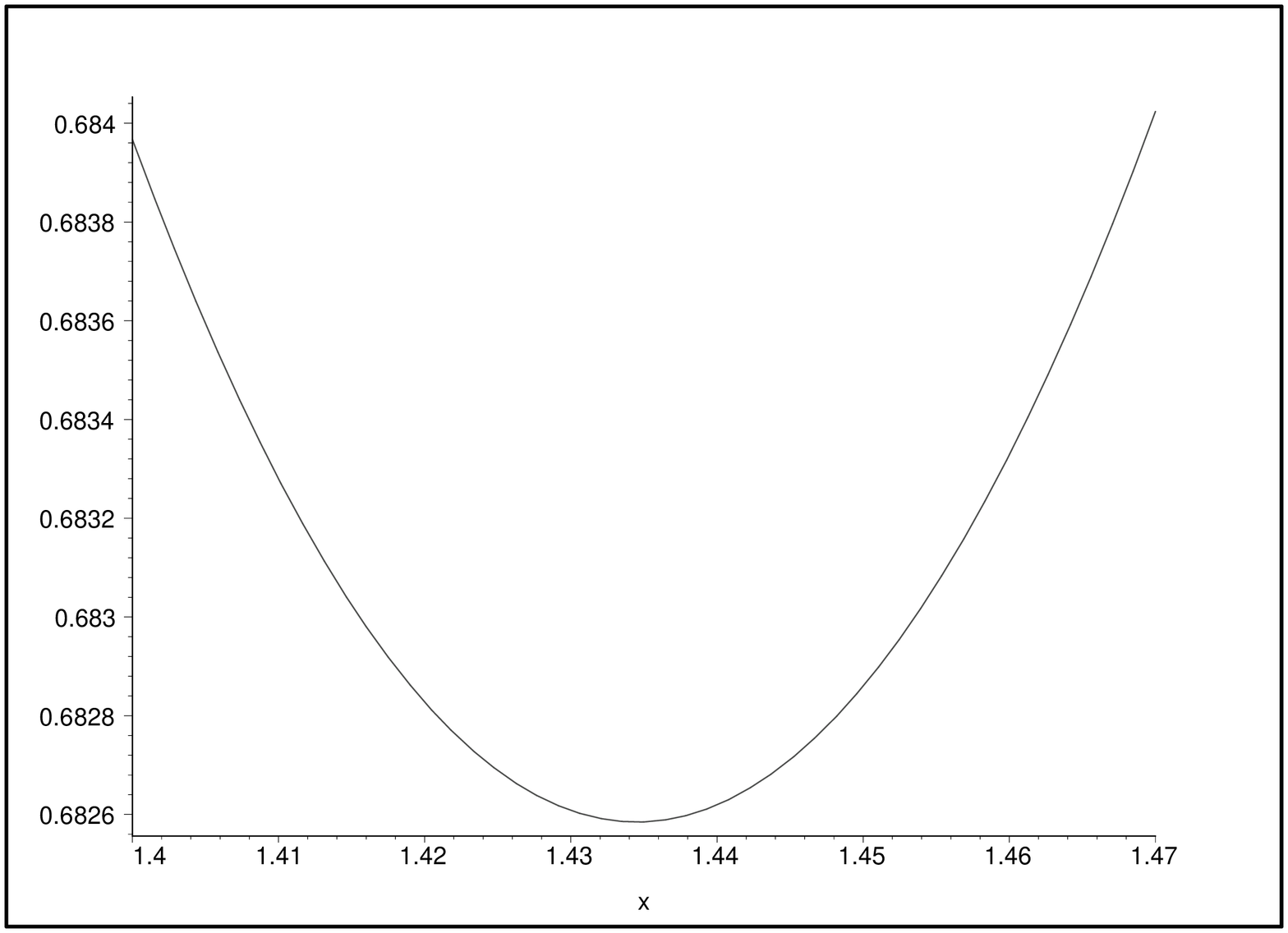,width=5.5cm}}
\quad \rotatebox{270}{\psfig{figure=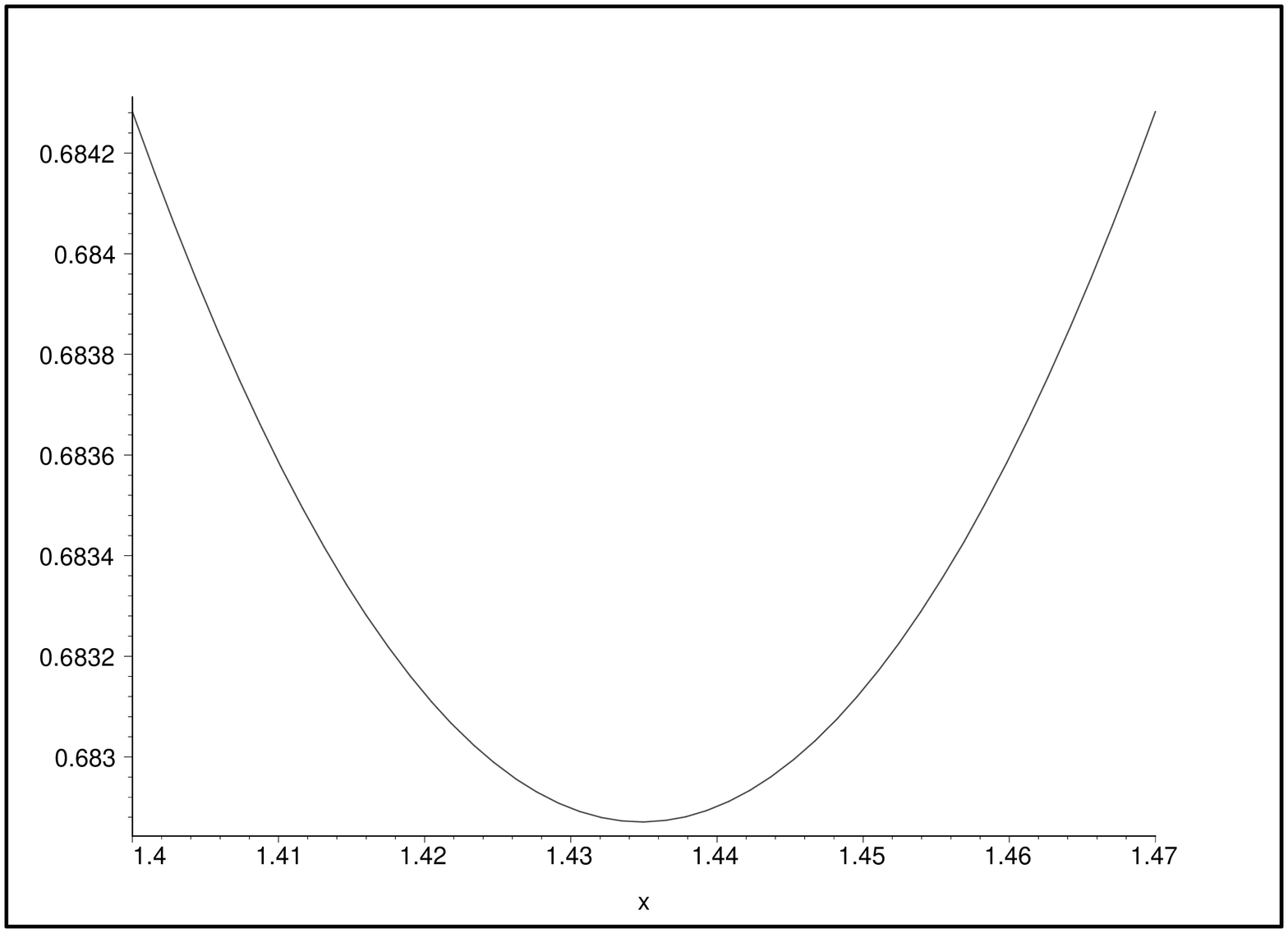,width=5.5cm}}}

\caption{the functions $2 e^{2H}$ and $\tilde{H}$ for $\rh = e$,
detail} \label{fig3}
\end{figure}

\begin{figure}
\centerline{\rotatebox{270}{\psfig{figure=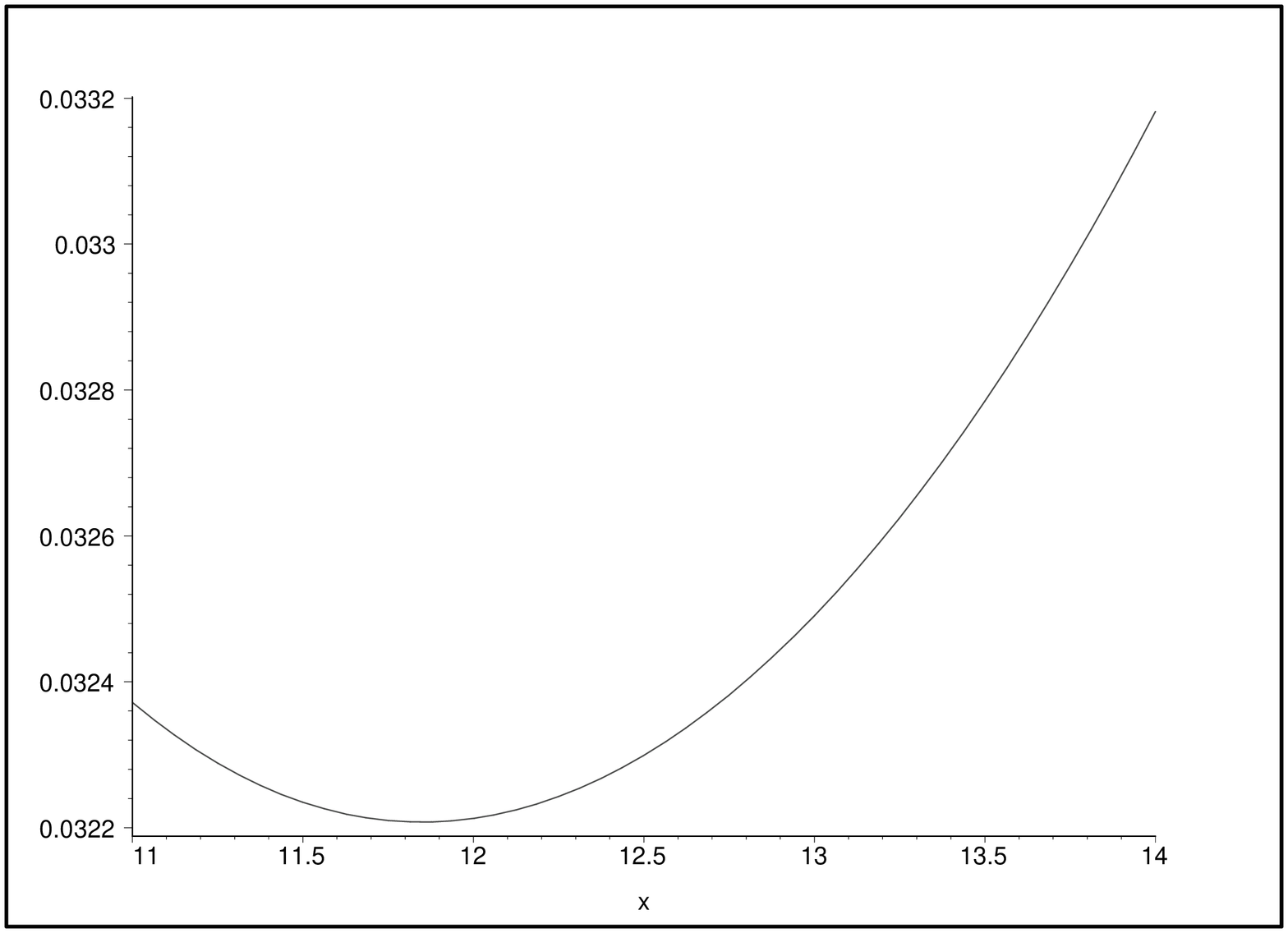,width=5.5cm}}
\quad \rotatebox{270}{\psfig{figure=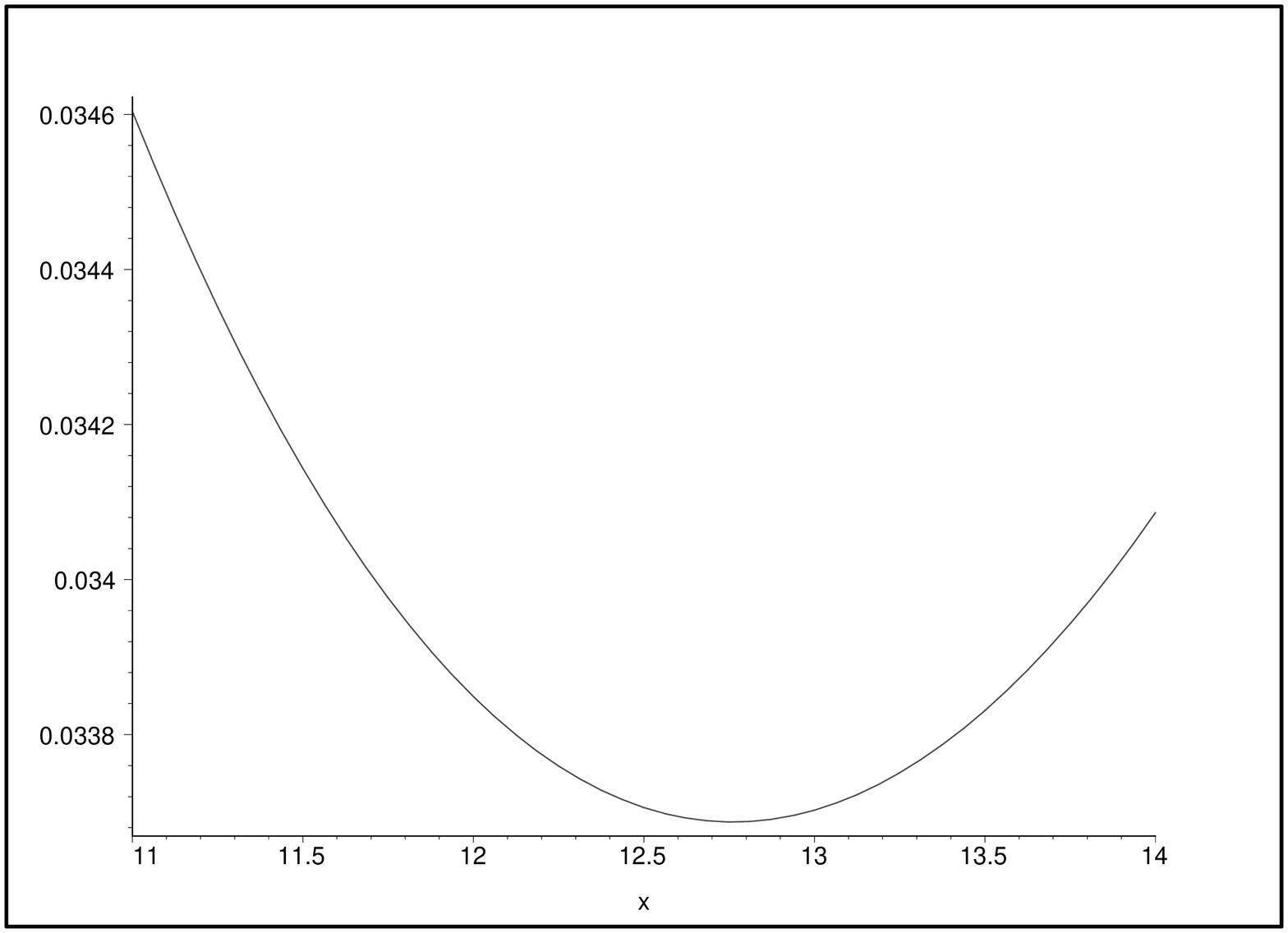,width=5.5cm}}}

\caption{the functions $2 e^{2H}$ and $\tilde{H}$ for $\log \rh =
3.5$} \label{fig4}
\end{figure}

For $\rh > 1$ we set
$$
A_\rh = \left\{ (x,y) \, : \, \frac{1}{\rh^2} < x^2 + y^2 < \rh^2
\right\}; \qquad S_\rh = \left\{ (x,y) \, : \, - \log \rh < x <
\log \rh \right\}.
$$
It is clear that $S_\rh$ is a covering of $A_\rh$ through the
exponential map. We also define $\a \in \C$,  $h_\rh : S_\rh \to
D$ and $f_\rh : D \to A_\rh$ by
$$
  \a = - \frac{i \pi}{2 \log \rh}; \qquad h_\rh(w) = \frac{e^{\a w} -
  1}{e^{\a w} + 1}, \qquad h_\rh^{-1}(z) = \frac{1}{\a} \log
  \left( \frac{1+z}{1-z} \right),  \qquad f = \exp \circ h_\rh^{-1}.
$$
where $w \in S_\rh$ and $z \in D$. Our aim is to compute formula
\eqref{eq:ffin2} for this particular case. Fixing $z_0 \in D$, the
points $z_k$ and the corresponding points $w_k = h_\rh^{-1} z_k$
are given by
$$
  z_k = \frac{e^{2 k \pi i \a} \left( \frac{1+z_0}{1-z_0}
  \right) - 1}{e^{2 k \pi i \a} \left( \frac{1+z_0}{1-z_0}
  \right) + 1}; \qquad \qquad
  w_k = \frac{1}{\a} \log \left( \frac{1+z_0}{1-z_0} \right) + 2 k \pi
  i.
$$
Using some elementary computations we obtain
$$
  z_k = \frac{z_0 + M_k}{z_0 M_k + 1}; \qquad \qquad \hbox{ where }
  M_k = \tanh \left( \frac{k \pi^2}{2 \log \rh} \right).
$$
If $T_k$ denotes as before the deck transformation, then there
holds
$$
  T_k (z) = \frac{z + M_k}{z M_k + 1}; \qquad \qquad T'_k (z) =
  \frac{1 - M_k^2}{(1+M_k z)^2}.
$$
By symmetry, it is sufficient to compute \eqref{eq:ffin} for
$f(z_0)$ real and positive. It is convenient to use the following
parametrization for the points $a \in A_\rh$ and $z_0 \in D$
$$
  a = \log x; \qquad z_0 = - i \tan \left( \frac{\pi}{4 \log \rh}
  \log x \right), \qquad x \in (- \log \rh, \log \rh).
$$
Using this notation, from equation \eqref{eq:ffin2} we are left
with
$$
\tilde{H}(a) = \frac{1}{|f'(z_0)|^2(1-|z_0|^2)^2} \sum_k \left(
\frac{(1-M_k^2)(1-|z_0|^2)^2}{(1 + M_k \ov{z}_0)^2(1 - \ov{z}_k
z_0)^2} + \frac{(1-M_k^2)(1-|z_0|^2)^2}{(1 + M_k z_0)^2(1 - z_k
\ov{z}_0)^2} \right).
$$
From the above formuls it follows
$$
\tilde{H}(a) = \frac{2}{|f'(z_0)|^2(1-|z_0|^2)^2} \sum_k
\frac{(1-M_k^2)(1-|z_0|^2)^2 \left( (1-|z_0|^2)^2 - 4 |z_0|^2
M_k^2\right)}{\left( (1-|z_0|^2)^2 + 4 |z_0|^2 M_k^2 \right)^2}.
$$
From \cite{baf} we have
$$
e^{2H(a,a)} = \frac{1}{|f'(z_0)|^2(1-|z_0|^2)^2} \prod_{z_k \neq
z_0} \left| \frac{z_k-z_0}{1-z_k\ov{z}_0} \right|^2.
$$
Using elementary computations it turns out that
\begin{equation}\label{eq:hex}
  \tilde{H}(a) = \frac{\pi^2}{8 (\log \rh)^2} \frac{1}{\cos^2
  \left( \frac{\pi}{2 \log \rh} \log x \right) |x|^2} \left( 1 + 2
  \sum_{k=1}^\infty W(k,x) \right);
\end{equation}
\begin{equation}\label{eq:gex}
  2 e^{2H(a,a)} = \frac{\pi^2}{8 (\log \rh)^2} \frac{1}{\cos^2
  \left( \frac{\pi}{2 \log \rh} \log x \right) |x|^2} \prod_{k=1}^\infty
  Z(k,x)^2,
\end{equation}
where
\begin{eqnarray*}
  W(k,x) & = & \left( 1 - \tanh^2 \left( \frac{k \pi^2}{2 \log \rh} \right)
  \right)^2 \left( 1 - \tan^2 \left( \frac{\pi}{4 \log \rh} \log x \right)
  \right)^2 \\ & \times & \frac{ \left( 1 - \tan^2 \left(
  \frac{\pi}{4 \log \rh} \log x \right) \right)^2  -
  4 \tanh^2 \left( \frac{k \pi^2}{2 \log \rh} \right)
  \tan^2 \left( \frac{\pi}{4 \log \rh} \log x \right)}{\left( \left( 1 -
  \tan^4 \left( \frac{\pi}{4 \log \rh} \log x \right) \right)^2 +
  4 \tanh^2 \left( \frac{k \pi^2}{2 \log \rh} \right)
  \tan^2 \left( \frac{\pi}{4 \log \rh} \log x \right) \right)^2}, \qquad k
  \geq 1;
\end{eqnarray*}
\begin{eqnarray*}
  Z(k,x) & = & \frac{\tanh^2 \left( \frac{k \pi^2}{2 \log \rh} \right)
  \sec^4 \left( \frac{\pi}{4 \log \rh} \log x \right)}{\left( 1 -
  \tan^4 \left( \frac{\pi}{4 \log \rh} \log x \right) \right)^2 +
  4 \tanh^2 \left( \frac{k \pi^2}{2 \log \rh} \right)
  \tan^2 \left( \frac{\pi}{4 \log \rh} \log x \right)}, \qquad k
  \geq 1.
\end{eqnarray*}
In Figures \ref{fig2}-\ref{fig4} we plot the functions $2 e^{2H}$
and $\tilde{H}$ (modulo the irrelevant factor $\frac{\pi^2}{8
(\log \rh)^2}$) for $\rh = e$ and for $\rh = e^{3.5}$. We note
that, roughly, $W(k,x) \sim e^{-\frac{2k\pi^2}{\log \rh}}$ and
$Z(k,x) \sim 1 - e^{-\frac{k\pi^2}{\log \rh}}$ so for small values
of $\rh$ the terms with $k \neq 0$ are almost negligible. This
accounts for the fact that for $\rh = e$ the graphs are very
similar, see Figure \ref{fig2}, even on a fine scale, see Figure
\ref{fig2}. For large values of $\rh$ the difference between the
two functions is mpre pronounced, see Figure \ref{fig4}.

\begin{rem}\label{r:hhr}
We recall that the {\em harmonic and hyperbolic radii} are defined
by
$$
r_{har} (\xi) = e^{-H(\xi,\xi)}; \qquad r_{hyp} (f(z)) = |f'(z)|
(1-|z|^2),
$$
where $f : D \to \O$ denotes a conformal covering map. See
\cite{baf} Definition 1 and page 15. Note that the function
$\tilde{H}$, in the case of general non-simply connected domains,
do not even coincide with $r_{hyp}^{-2}$. However, in the case of
small annuli, numerical computation show that $r_{har}$ and
$r_{hyp}$ are very close, see \cite{baf} Figure 8. We also point
out that the harmonic radius is related to the Bergman kernel, see
\cite{baf} Section 8.4.
\end{rem}

\section{The expansion for multiple bubbles}\label{s:2b}

In this section we consider the case of multiple bubbles. We begin
considering only two bubbles $R_1 \d_{a,\l_1}$ and $R_2
\d_{b,\l_2}$, which we denote for simplicity by $\d_1$ and $\d_2$
respectively. We assume that $R_1 = Id$, namely the first bubble
is not rotated, and we simply write $R$ for $R_2$.

For $\ov{C} > 0$, $k = 2$ and $P \d_1 + P \d_2 \in Z$, our aim is
to expand the functional $I_\e$ (see \eqref{eq:pe}) on $Z$ in
terms of the parameters $a = p_1$, $b = p_2$, $\l_1, \l_2$ and
$R$. In the following, for brevity, we set (see Section
\ref{s:pre})
$$
\s = b - a, \qquad e(\l_1, \l_2) = O\left( (\log \l_1 + \log \l_2)
\left( \frac{1}{\l_1^3} + \frac{1}{\l_2^3} + \frac{1}{\l_1^2 \l_2}
+ \frac{1}{\l_1 \l_2^2} \right) \right).
$$
We recall the explicit form of the functional $I_\e(u)$, for $u
\in H^1_0(\O;\R^3)$
$$
I_\e(u) = \frac{1}{2} \int_\O |\n u|^2 + \frac{2}{3} \int_\O u
  \cdot (u_x \wedge u_y) + \e \int_{\O} u \cdot (u_x \wedge g_y +
  g_x \wedge u_y) + 2 \e^2 \int_\O u \cdot (g_x \wedge g_y).
$$

\subsection{Interaction with $g$}

We consider first the interaction term $\e \IO u (u_x \wedge g_y +
g_x \wedge u_y)$ in \eqref{eq:pe}, with $u = P \d_1 + P \d_2$. We
recall that throughout this section we assume that $\ov{C}^{-1}
\l_i^{-1} \leq \e \leq \ov{C} \l_i^{-1}$ for some fixed constant
$\ov{C}$ and for $i = 1, 2$. We have
\begin{eqnarray}\label{eq:14}
\IO u \cdot (u_x \wedge g_y + g_x \wedge u_y) & = &  \IO (P \d_1 +
P \d_2) \cdot \left[ ((P \d_1)_x + (P \d_2)_x) \wedge g_y + g_x
\wedge ((P \d_1)_y + (P \d_2)_y) \right] \nonumber \\ & = &  \IO P
\d_1 \cdot ((P \d_1)_x \wedge g_y + g_x \wedge (P \d_1)_y) + \IO P
\d_1 \cdot ((P \d_2)_x \wedge g_y + g_x \wedge (P \d_2)_y) \\ & +
& \IO P \d_2 \cdot ((P \d_2)_x \wedge g_y + g_x \wedge (P \d_2)_y)
+ \IO P \d_2 \cdot ((P \d_1)_x \wedge g_y + g_x \wedge (P
\d_1)_y). \nonumber
\end{eqnarray}
The first term in \eqref{eq:14} has been estimated in Section
\ref{s:1b}, formulas \eqref{eq:ef1}-\eqref{eq:ef2}, and gives
\begin{equation}\label{eq:2b1}
\e \IO P \d_1 ((P \d_1)_x \wedge g_y + g_x \wedge (P \d_1)_y) = -
16 \frac{\e}{\l_1} \IO \frac{x^2}{(1+|\xi|^2)^3} d_{Id}g(a) +
O(\l_1^{-2}).
\end{equation}
The third term in \eqref{eq:14} can be estimated similarly, using
the invariance of the problem under rotation, and gives
\begin{equation}\label{eq:2b3}
\e \IO P \d_2 ((P \d_2)_x \wedge g_y + g_x \wedge (P \d_2)_y) = -
16 \frac{\e}{\l_2} \IO \frac{x^2}{(1+|\xi|^2)^3} d_{R^{-1}}g(b) +
O(\l_2^{-2}).
\end{equation}
Let us compute now the remaining two terms in \eqref{eq:14},
starting from the fourth. We write $a = (a_1,a_2)$ and $b =
(b_1,b_2)$. Up to an error of order $\e^2$, we have
\begin{equation}\label{eq:2b41}
  (P \d_1)_x \wedge g_y =
  \begin{pmatrix}
    \left( - \frac{4 \l_1^3 (x - a_1) (y - a_2)}{(1 + \l_1^2 |\xi-a|^2)^2}
    - \frac{2}{\l_1} \partial_x h_2 \right) (g_3)_y - \left(
    \frac{4 \l_1^2 (x - a_1)}{(1 + \l_1^2 |\xi-a|^2)^2} + \frac{2}{\l_1^2}
    \partial_x h_3 \right) (g_2)_y \\
    \left( \frac{4 \l_1^2 (x - a_1)}{(1 + \l_1^2 |\xi-a|^2)^2} + \frac{2}{\l_1^2}
    \partial_x h_3 \right) (g_1)_y - \left( 2 \l_1
    \frac{1 + \l_1^2 ((y - a_2)^2 - (x - a_1)^2)}{(1 + \l_1^2 |\xi-a|^2)^2}
    - \frac{2}{\l_1} \partial_x h_1 \right) (g_3)_y \\
    \left( 2 \l_1
    \frac{1 + \l_1^2 ((y - a_2)^2 - (x - a_1)^2)}{(1 + \l_1^2 |\xi-a|^2)^2}
    - \frac{2}{\l_1} \partial_x h_1 \right) (g_2)_y + \left(
    \frac{4 \l_1^3 (x - a_1) (y - a_2)}{(1 + \l_1^2 |\xi-a|^2)^2}
    + \frac{2}{\l_1} \partial_x h_2 \right) (g_1)_y
  \end{pmatrix}
\end{equation}
and
\begin{eqnarray}\label{eq:2b42}
  P \d_2 =
  \begin{pmatrix}
    r_{11} & r_{12} & r_{13} \\
    r_{21} & r_{22} & r_{23} \\
    r_{31} & r_{32} & r_{33} \
  \end{pmatrix}
  \nonumber
  \begin{pmatrix}
    \frac{2 \l_2 (x - b_1)}{1 + \l_2^2 |\xi - b|^2} - \frac{2}{\l_2} h_1(b,\xi)
    \\ \frac{2 \l_2 (y - b_2)}{1 + \l_2^2 |\xi - b|^2} - \frac{2}{\l_2} h_2(b,\xi)
    \\ \frac{2}{\l_2^2} h_3 - \frac{2}{1 + \l_2^2 |\xi - b|^2}
  \end{pmatrix},
\end{eqnarray}
where $r_{ij}$ are the entries of the matrix $R$. We are going to
prove that
\begin{equation}\label{eq:cl}
  \IO P \d_2 ((P \d_1)_x \wedge g_y + g_x \wedge (P \d_1)_y) =
  O\left( \e^2 | \log \e | \right).
\end{equation}
If one uses \eqref{eq:2b41} and \eqref{eq:2b42}, the integrals
involved in the above expression are of the form
\begin{eqnarray}\label{eq:i1}
  \l_1^3 \l_2 \int_{\O} \frac{(x - a_1) (y - a_2)
  (x - b_1)}{\left( 1 + \l_1^2|\xi-a|^2 \right)^2
  \left( 1 + \l_2^2|\xi-b|^2 \right)}; \quad
  \l_1^3 \l_2 \int_{\O} \frac{(x - a_1) (y - a_2)
  (y - b_2)}{\left( 1 + \l_1^2|\xi-a|^2 \right)^2
  \left( 1 + \l_2^2|\xi-b|^2 \right)};   \\ \label{eq:i2}
  \l_1^3 \int_{\O} \frac{(x - a_1)
  (y - a_2)}{\left( 1 + \l_1^2|\xi-a|^2 \right)^2
  \left( 1 + \l_2^2|\xi-b|^2 \right)}; \quad
  \frac{\l_1^3}{\l_2} \int_{\O} \frac{(x - a_1)
  (y - a_2)}{\left( 1 + \l_1^2|\xi-a|^2 \right)^2}; \quad
  \frac{\l_1^3}{\l_2^2} \int_{\O} \frac{(x - a_1)
  (y - a_2)}{\left( 1 + \l_1^2|\xi-a|^2 \right)^2}
\end{eqnarray}
\begin{eqnarray} \label{eq:i3}
  \l_1^2 \l_2 \int_{\O} \frac{(x - a_1)
  (x - b_1)}{\left( 1 + \l_1^2|\xi-a|^2 \right)^2
  \left( 1 + \l_2^2|\xi-b|^2 \right)}; \quad
  \l_1^2 \l_2 \int_{\O} \frac{(x - a_1)
  (y - b_2)}{\left( 1 + \l_1^2|\xi-a|^2 \right)^2
  \left( 1 + \l_2^2|\xi-b|^2 \right)};   \\ \label{eq:i4}
  \l_1^2 \int_{\O} \frac{(x - a_1)
  }{\left( 1 + \l_1^2|\xi-a|^2 \right)^2
  \left( 1 + \l_2^2|\xi-b|^2 \right)}; \quad
  \frac{\l_1^2}{\l_2} \int_{\O} \frac{(x - a_1)
  }{\left( 1 + \l_1^2|\xi-a|^2 \right)^2}; \quad
  \frac{\l_1^2}{\l_2^2} \int_{\O} \frac{(x - a_1)
  }{\left( 1 + \l_1^2|\xi-a|^2 \right)^2}
  \end{eqnarray}
\begin{eqnarray}\label{eq:i5}
  \qquad \l_1 \l_2 \int_{\O} \frac{1 + \l_1^2 ((y-a_2)^2 -
  (x-a_1)^2) (x - b_1)}{\left( 1 + \l_1^2|\xi-a|^2 \right)^2
  \left( 1 + \l_2^2|\xi-b|^2 \right)}; \quad \!
  \l_1 \l_2 \int_{\O} \frac{1 + \l_1^2 ((y-a_2)^2 - (x-a_1)^2)
  (y - b_2)}{\left( 1 + \l_1^2|\xi-a|^2 \right)^2
  \left( 1 + \l_2^2|\xi-b|^2 \right)};   \\ \label{eq:i6}
  \l_1 \int_{\O} \frac{1 + \l_1^2 ((y-a_2)^2 -
  (x-a_1)^2)}{\left( 1 + \l_1^2|\xi-a|^2 \right)^2
  \left( 1 + \l_2^2|\xi-b|^2 \right)}; \quad
  \frac{\l_1}{\l_2} \int_{\O} \frac{1 + \l_1^2 ((y-a_2)^2
  - (x-a_1)^2)}{\left( 1 + \l_1^2|\xi-a|^2 \right)^2};
  \end{eqnarray}
\begin{eqnarray}\label{eq:i7}
  \frac{\l_1}{\l_2^2} \int_{\O} \frac{1 + \l_1^2 ((y-a_2)^2
  - (x-a_1)^2)}{\left( 1 + \l_1^2|\xi-a|^2 \right)^2} \quad
  \frac{\l_2}{\l_1} \int_{\O} \frac{x-b_1}{\left( 1 +
  \l_2^2|\xi-b|^2 \right)}; \quad \frac{\l_2}{\l_1} \int_{\O}
  \frac{y-b_2}{\left( 1 + \l_2^2|\xi-b|^2 \right)};  \\ \label{eq:i8}
  \frac{1}{\l_1} \int_{\O} \frac{1}{\left( 1 +
  \l_1^2|\xi-b|^2 \right)}; \quad \int_{\O}
  \frac{1}{\l_1 \l_2}; \quad \int_{\O} \frac{1}{\l_1 \l_2^2};
  \quad \frac{\l_2}{\l_1^2} \int_{\O} \frac{x-b_1}{\left( 1 +
  \l_2^2|\xi-b|^2 \right)};  \\ \label{eq:i9}
  \frac{\l_2}{\l_1^2} \int_{\O}
  \frac{y-b_2}{\left( 1 + \l_2^2|\xi-b|^2 \right)}; \quad
  \frac{1}{\l_1^2} \int_{\O} \frac{1}{\left( 1 +
  \l_2^2|\xi-b|^2 \right)}; \quad \int_{\O} \frac{1}{\l_1^2 \l_2};
  \quad \int_{\O} \frac{1}{\l_1^2 \l_2^2}.
\end{eqnarray}
The errors in the expressions of $(P \d_1)_x \wedge g_y$ and $P
\d_2$ are negligible with respect to the quantities listed in
\eqref{eq:i1}-\eqref{eq:i9}, hence it is sufficient to consider
the above expressions.

\

\noindent {\bf Estimate of \eqref{eq:i1}}. Using the rescaling
$\xi \mapsto \l_1 (\xi - a)$ and setting $\O_{\l_1,a} = \l_1 (\O -
a)$, we get
\begin{eqnarray}\label{eq:na1}
  \l_1^3 \l_2 \int_{\O} \frac{(x - a_1) (y - a_2)
  (x - b_1)}{\left( 1 + \l_1^2|\xi-a|^2 \right)^2
  \left( 1 + \l_2^2|\xi-b|^2 \right)} & = & \frac{\l_2}{\l_1^2}
  \int_{\O_{\l_1,a}} \frac{x^2 y}{(1+|\xi|^2)^2
  \left( 1 + \frac{\l_2^2}{\l_1^2}|\xi-\l_1 \s|^2 \right)} \nonumber \\
  & + & \frac{\l_2}{\l_1} \s_1
  \int_{\O_{\l_1,a}} \frac{x y}{(1+|\xi|^2)^2
  \left( 1 + \frac{\l_2^2}{\l_1^2}|\xi-\l_1 \s|^2 \right)}.
\end{eqnarray}
Let us consider the first integral. We divide $\O_{\l_1,a}$ into
the regions $|\xi| \leq \frac{\l_1 \s}{2}$ and $|\xi| \geq
\frac{\l_1 \s}{2}$. If $|\xi| \leq \frac{\l_1 \s}{2}$, then
$|\xi-\l_1 \s|^2 \geq \frac{\l_1^2 |\s|^2}{4}$, so we have
$$
\frac{\l_2}{\l_1^2}
  \left| \int_{|\xi| \leq \frac{\l_1 \s}{2}} \frac{x^2 y}{(1+|\xi|^2)^2
  \left( 1 + \frac{\l_2^2}{\l_1^2}|\xi-\l_1 \s|^2 \right)} \right|
  \leq \frac{C \l_2}{\l_1^2 \l_2^2 |\s|^2} \int_{|\xi| \leq
  \frac{\l_1 \s}{2}} \frac{|\xi|^3}{(1+|\xi|^2)^2} \leq
  \frac{C \l_1 \l_2}{\l_1^2 \l_2^2 |\s|^2},
$$
where $C$ is a positive constant depending only on $\O$. When
$|\xi| \geq \frac{\l_1 |\s|}{2}$ then for $\l_1$ large there holds
$\frac{|\xi|^3}{(1+|\xi|^2)^2} \leq \frac{2}{\l_1 |\s|}$. As a
consequence, using also the change of variables $\frac{\l_2}{\l_1}
(\xi - \l_1 \s) \mapsto \xi$, we deduce
$$
\frac{\l_2}{\l_1^2} \left| \int_{|\xi| \geq \frac{\l_1 \s}{2}}
\frac{x^2 y}{(1+|\xi|^2)^2 \left( 1 +
\frac{\l_2^2}{\l_1^2}|\xi-\l_1 \s|^2 \right)} \right| \leq \frac{C
\l_2}{\l_1^3 |\s|} \left( \frac{\l_1}{\l_2} \right)^2
\int_{\O_{\l_2,b}} \frac{1}{(1+|\xi|^2)} \leq C \frac{\log
\l_2}{\l_1 \l_2 |\s|}.
$$
Turning to the second integral in the r.h.s. of \eqref{eq:na1}, we
again divide the domain into two regions $|\xi| \leq \frac{\l_1
\s}{2}$ and $|\xi| \geq \frac{\l_1 \s}{2}$. Reasoning as before we
find
\begin{eqnarray*}
  \frac{\l_2}{\l_1} \s_1 \left| \int_{|\xi| \leq \frac{\l_1 \s}{2}} \frac{x
y}{(1+|\xi|^2)^2 \left( 1 + \frac{\l_2^2}{\l_1^2}|\xi-\l_1 \s|^2
\right)} \right| \leq C \frac{\l_2 \s_1}{\l_1 \l_2^2 |\s|^2}
\int_{|\xi| \leq \frac{\l_1 \s}{2}} \frac{|\xi|^2}{(1+|\xi|^2)^2}
\leq C
\frac{|\s_1|}{|\s|^2} \frac{\log \l_1}{\l_1 \l_2}; \\
\frac{\l_2}{\l_1} \s_1 \left| \int_{|\xi| \geq \frac{\l_1 \s}{2}}
\frac{x y}{(1+|\xi|^2)^2 \left( 1 + \frac{\l_2^2}{\l_1^2}|\xi-\l_1
\s|^2 \right)} \right| \leq C \frac{\l_2 \s_1}{\l_1^3 |\s|^2}
\left( \frac{\l_1}{\l_2} \right)^2 \int_{\O_{\l_2,b}} \frac{1}{(1
+ |\xi|^2)} \leq C \frac{|\s_1|}{|\s|^2} \frac{\l_2 \log
\l_2}{\l_1^3}.
\end{eqnarray*}
Since in the definition of $Z$ we assume dist$(p_i,p_j) \geq
\ov{C}^{-1}$, $|\s|$ is uniformly bounded from below, the last
formulas imply
$$
\l_1^3 \l_2 \left| \int_{\O} \frac{(x - a_1) (y - a_2)
  (x - b_1)}{\left( 1 + \l_1^2|\xi-a|^2 \right)^2
  \left( 1 + \l_2^2|\xi-b|^2 \right)} \right| \leq C \e^2 | \log \e |.
$$
The second expression in \eqref{eq:i1} can be estimated in the
same way.

\

\noindent {\bf Estimate of \eqref{eq:i2}-\eqref{eq:i9}.} We just
treat some particular cases, since many terms are similar to
each-other. First of all, all the terms for which the quantity
$\left( 1 + \l_1^2|\xi-a|^2 \right)^2 \left( 1 + \l_2^2|\xi-b|^2
\right)$ appears in the denominator can be treated as before.

Next, we consider for example the second term in \eqref{eq:i6} and
the last term in \eqref{eq:i7}. Using the changes of variable
$\l_1 (\xi-a) \mapsto \xi$ and $\l_2 (\xi-b) \mapsto \xi$ we find
\begin{eqnarray*}
  \frac{\l_1}{\l_2} \int_{\O} \frac{1 + \l_1^2 ((y-a_2)^2
  - (x-a_1)^2)}{\left( 1 + \l_1^2|\xi-a|^2 \right)^2} =
  \frac{1}{\l_1 \l_2} \int_{\O_{\l_1,a}}
  \frac{1+y^2-x^2}{\left( 1 + |\xi|^2 \right)^2} \leq C
  \frac{\log \l_1}{\l_1 \l_2} \leq C \e^2 | \log \e |;
  \\ \frac{\l_2}{\l_1} \int_{\O}
  \frac{y-b_2}{\left( 1 + \l_2^2|\xi-b|^2 \right)} =
  \frac{1}{\l_1 \l_2^2} \int_{\O_{\l_2,b}} \frac{|\xi|}{1+|\xi|^2}
  \leq C \frac{1}{\l_1 \l_2} \leq C \e^2.
\end{eqnarray*}

\

\noindent {\bf Conclusion.} Using the estimates of \eqref{eq:i1}
and those of \eqref{eq:i2}-\eqref{eq:i9}, we obtain \eqref{eq:cl}.
In the same way, one can prove that
\begin{equation}\label{eq:cl2}
  \IO P \d_1 ((P \d_2)_x \wedge g_y + g_x \wedge (P \d_2)_y) =
  O \left( \e^2 | \log \e | \right).
\end{equation}
Hence, from equations \eqref{eq:2b1}, \eqref{eq:2b3},
\eqref{eq:cl} and \eqref{eq:cl2} we deduce

\begin{lem}\label{l:intg}
For $u = P \d_1 + P \d_2 \in Z$ there holds
\begin{equation}\label{eq:intg}
  \e \IO u (u_x \wedge g_y + g_x \wedge u_y) = - 8 A_0 d_{Id}g(a)
  - 8 A_0 d_{R^{-1}}g(b) + O \left( \e^2 | \log \e | \right).
\end{equation}
\end{lem}

\subsection{Mixed terms in $P \d_1$ and $P \d_2$}

For $u = P \d_1 + P \d_2$, we consider the first and the second
integrals in \eqref{eq:ie}. We are interested in the terms
involving both $P \d_1$ and $P \d_2$, namely
\begin{eqnarray*}
  \IO \n P \d_1 \cdot \n P \d_2 & + & \frac{2}{3} \IO P \d_1 \cdot
  \left( (P \d_2)_x \wedge (P \d_2)_y \right) + \frac{2}{3} \IO P
  \d_1 \cdot \left( (P \d_2)_x \wedge (P \d_1)_y \right) \\ & + &
  \frac{2}{3} \IO P \d_1 \cdot
  \left( (P \d_1)_x \wedge (P \d_2)_y \right) +
  \frac{2}{3} \IO P \d_2 \cdot
  \left( (P \d_2)_x \wedge (P \d_1)_y \right) \\ & + &
  \frac{2}{3} \IO P \d_2 \cdot
  \left( (P \d_1)_x \wedge (P \d_2)_y \right) +
  \frac{2}{3} \IO P \d_2 \cdot
  \left( (P \d_1)_x \wedge (P \d_1)_y \right).
\end{eqnarray*}
Integrating by parts it is easy to we see that the last expression
becomes
$$
  \IO \n P \d_1 \cdot \n P \d_2 + 2 \IO P \d_2 \cdot
  \left( (P \d_1)_x \wedge (P \d_1)_y \right) +
  2 \IO P \d_1 \cdot
  \left( (P \d_2)_x \wedge (P \d_2)_y \right).
$$
Integrating by parts the first term we get
$$
\IO \n P \d_1 \cdot \n P \d_2 = - \IO \D P \d_2 \cdot P \d_1 = - 2
\IO P \d_1 \cdot ((\d_2)_x \wedge (\d_2)_y),
$$
so we are left with
$$
  2 \IO P \d_2 \left( (P \d_1)_x \wedge (P \d_1)_y \right) +
  2 \IO P \d_1 \left( (P \d_2)_x \wedge (P \d_2)_y - (\d_2)_x
  \wedge (\d_2)_y \right).
$$

\begin{lem}\label{l:apex}
For $P \d_1 + P \d_2 \in Z$, there holds
\begin{eqnarray}\label{eq:apex}
  2 \IO P \d_2 \left( (P \d_1)_x \wedge (P \d_1)_y \right) & + &
  2 \IO P \d_1 \left( (P \d_2)_x \wedge (P \d_2)_y - (\d_2)_x
  \wedge (\d_2)_y \right) \nonumber \\ & = &
  2 \IO P \d_2 \left( (\d_1)_x
  \wedge (\d_1)_y \right) + e(\l_1, \l_2).
\end{eqnarray}
\end{lem}

\begin{pf}
Since the difference between the l.h.s. and the r.h.s. of
\eqref{eq:apex} is
$$
2 \IO P \d_1 \left( (P \d_2)_x \wedge (P \d_2)_y - (\d_2)_x
  \wedge (\d_2)_y \right) + 2 \IO P \d_2 \left( (P \d_1)_x
  \wedge (P \d_1)_y - (\d_1)_x \wedge (\d_1)_y \right),
$$
it is sufficient by symmetry to estimate one of the two terms in
the last expression. We have
\begin{eqnarray*}
  \IO P \d_1 \left( (P \d_2)_x \wedge (P \d_2)_y - (\d_2)_x
  \wedge (\d_2)_y \right) = \IO P \d_1 \cdot \left[ (\var_2)_x \wedge
  (\var_2)_y - (\var_2)_x \wedge (\d_2)_y - (\d_2)_x \wedge
  (\var_2)_y \right].
\end{eqnarray*}
Using equations \eqref{eq:est}-\eqref{eq:est2}, choosing $\t \leq
\frac{1}{4} \ov{C}$ and setting $\O_{a,b} = \O \setminus
(B_{\t}(a) \cup B_{\t}(b))$, we find
\begin{eqnarray*}
  \left| \IO P \d_1 \cdot \left( (\var_2)_x \wedge (\var_2)_y
  \right) \right| & \leq & \frac{C}{\l_1 \l_2^2}; \\
  \left| \IO P \d_1 \cdot \left( (\var_2)_x \wedge (\d_2)_y
  \right) \right| & \leq & \frac{C}{\l_2} \left[ \int_{B_{\t}(a)}
  |P \d_1| \, |(\d_2)_y|  +  \int_{B_{\t}(b)}
  |P \d_1| \, |(\d_2)_y|  +  \int_{\O_{a,b}} |P \d_1| \, |(\d_2)_y| \right]
  \\ & \leq & \frac{C}{\l_2} \left[ \frac{C}{\l_1 \l_2} + C
  \frac{\log \l_2}{\l_1 \l_2} + \frac{C}{\l_1 \l_2} \right],
\end{eqnarray*}
and an analogous estimate for the term $P \d_1 \cdot \left(
(\d_2)_x \wedge (\var_2)_y \right)$. This concludes the proof.
\end{pf}

\begin{lem}\label{l:apex2}
For $P \d_1 + P \d_2 \in Z$, there holds
\begin{eqnarray}\label{eq:apex2}
  2 \IO P \d_2 \cdot \left( (\d_1)_x \wedge (\d_1)_y \right) & = &
  \frac{16 A_0}{\l_1 \l_2} r_{11} \left( \frac{\s_1^2 - \s_2^2}{|\s|^4}
  + \frac{\partial h_1}{\partial x}(a,b) \right) +
  \frac{16 A_0}{\l_1 \l_2} r_{22} \left( \frac{\s_2^2 - \s_1^2}{|\s|^4}
  + \frac{\partial h_2}{\partial y}(a,b) \right) \nonumber \\
  & + & \frac{16 A_0}{\l_1 \l_2} r_{12} \left( 2 \frac{\s_1 \s_2}{|\s|^4}
  + \frac{\partial h_1}{\partial y}(a,b) \right) +
  \frac{16 A_0}{\l_1 \l_2} r_{21} \left( 2 \frac{\s_1 \s_2}{|\s|^4}
  + \frac{\partial h_2}{\partial x}(a,b) \right) \\ & + & e(\l_1,\l_2).
  \nonumber
\end{eqnarray}
\end{lem}

\begin{pf}
The left-hand side of \eqref{eq:apex2} is given explicitly by
\begin{equation}\label{eq:mtm}
  2 \IO P \d_2 \cdot ((\d_1)_x \wedge (\d_1)_y) = 2 \sum_{i,j=1}^3 r_{ij}
  \IO (\d_2)_j \cdot ((\d_1)_x \wedge (\d_1)_y)_i - 2 \sum_{i,j=1}^3 r_{ij}
  \IO (\var_2)_j \cdot ((\d_1)_x \wedge (\d_1)_y)_i,
\end{equation}
where $\{r_{ij}\}$ are the entries of the matrix $R$. We are now
going to estimate these integrals. We recall that, by
\eqref{eq:estfi}
$$
\var_2 (\xi) = \left( \frac{2}{\l_2} h_1(\xi,b) + O(\l_2^{-2}),
\frac{2}{\l_2} h_2(\xi,b) + O(\l_2^{-2}), 1 - \frac{2}{\l_2^2}
h_3(\xi,b) + O(\l_2^{-3}) \right).
$$
Taking this into account, we find that the terms in \eqref{eq:mtm}
involving the coefficients $r_{11}$, $r_{12}$, $r_{13}, r_{31}$
and $r_{33}$ are given respectively by
\begin{eqnarray}\label{eq:c11}
  - 32 \l_1^3 \l_2 \IO \frac{(x-b_1)(x-a_1)}{(1+\l_1^2|\xi-a|^2)^3
  (1+\l_2^2|\xi-b|^2)} + 32 \l_1^3 \IO
  \frac{(\l_2^{-1}h_1(\xi,b) +
  O(\l_2^{-2}))(x_1-a_1)}{(1+\l_1^2|\xi-a|^2)^3};
\\ \label{eq:c12}
  - 32 \l_1^3 \l_2 \IO \frac{(x-b_1)(y-a_2)}{(1+\l_1^2|\xi-a|^2)^3
  (1+\l_2^2|\xi-b|^2)} + 32 \l_1^3 \IO
  \frac{(\l_2^{-1}h_1(\xi,b) +
  O(\l_2^{-2}))(y-a_2)}{(1+\l_1^2|\xi-a|^2)^3};
  \\  \label{eq:c13}
  \qquad \quad 16 \l_1^2 \l_2 \IO \frac{(x-b_1)\left( 1 - \l_1^2 |\xi-a|^2
  \right)}{(1+\l_1^2|\xi-a|^2)^3 (1+\l_2^2|\xi-b|^2)} - 16 \l_1^2 \IO
  \frac{(\l_2^{-1}h_1(\xi,b) + O(\l_2^{-2}))\left( 1 - \l_1^2|\xi-a|^2
  \right)}{(1+\l_1^2|\xi-a|^2)^3};
\end{eqnarray}
\begin{eqnarray}\label{eq:c31}
  16 \l_1^3 \IO \left( \frac{2}{1+\l_2^2|\xi-b|^2} - \frac{2}{\l_2^2}
  h_3(\xi,b) + O(\l_2^{-3}) \right) \frac{(x-a_1)}{(1+\l_1^2|\xi-a|^2)^3};
  \\ \label{eq:c33} - 8 \l_1^2 \IO \left( \frac{2}{1+\l_2^2|\xi-b|^2} -
  \frac{2}{\l_2^2} h_3(\xi,b) + O(\l_2^{-3}) \right)
  \frac{(1-\l_1^2|\xi-a|^2)}{(1+\l_1^2|\xi-a|^2)^3}.
\end{eqnarray}
The terms involving the other coefficients of the matrix $R$ can
be estimated using the above ones, and will be taken into account
later.

\

\noindent {\bf Estimate of \eqref{eq:c11}.} Using the change of
variables $\l_1 (\xi-a) \mapsto \xi$, equation \eqref{eq:c11}
becomes
\begin{equation}\label{eq:c112}
  - 32 \frac{\l_2}{\l_1} \int_{\O_{\l_1,a}} \frac{x (x - \l_1
  \s_1)}{(1+|\xi|^2)^3 \left(1+\frac{\l_2^2}{\l_1^2}|\xi-\l_1 \s|^2\right)}
  + 32 \int_{\O_{\l_1,a}} \frac{(\l_2^{-1}h_1(\l_1^{-1}\xi + a,b) +
  O(\l_2^{-2}))x}{(1+|\xi|^2)^3}.
\end{equation}
We estimate the first term in \eqref{eq:c112}. Consider the
following subsets of the domain of integration
$$
\mathcal{B}_1 = \left\{ \xi \in \O_{\l_1,a} \, : \, |\xi| \leq
\frac{\l_1 |\s|}{4} \right\}, \quad \mathcal{B}_2 = \left\{ \xi
\in \O_{\l_1,a} \, : \, |\xi - \l_1 \s| \geq \frac{\l_1 |\s|}{4}
\right\}, \quad \mathcal{B}_3 = \O_{\l_1,a} \setminus
(\mathcal{B}_1 \cup \mathcal{B}_2).
$$
We can write
\begin{equation}\label{eq:s2}
  \left( 1+\frac{\l_2^2}{\l_1^2}|\xi-\l_1 \s|^2 \right)^{-1} = (1 +
  \l_2^2 |\s|^2)^{-1} \left( 1 + \frac{\frac{\l_2^2}{\l_1^2}|\xi|^2 -
  2 \frac{\l_2^2}{\l_1} \xi \cdot\s}{(1 + \l_2^2 |\s|^2)}
  \right)^{-1}.
\end{equation}
On the set $\mathcal{B}_1$ we have the following inequality
$$
\left| \frac{\l_2^2}{\l_1^2} |\xi|^2 - 2 \frac{\l_2^2}{\l_1} \xi
\cdot\s \right| \leq \frac{\l_2^2 |\s|^2}{16} + \frac{\l_2^2
|\s|^2}{2} \leq \frac{3}{4} \l_2^2 |\s|^2
$$
hence, from a Taylor expansion, we obtain the following uniform
estimate
\begin{equation}\label{eq:ue2}
  \left| \left( 1 + \frac{\frac{\l_2^2}{\l_1^2}|\xi|^2 -
  2 \frac{\l_2^2}{\l_1} \xi \cdot\s}{(1 + \l_2^2 |\s|^2)}
  \right)^{-1} - 1 + \frac{\frac{\l_2^2}{\l_1^2}|\xi|^2 -
  2 \frac{\l_2^2}{\l_1} \xi \cdot\s}{(1 + \l_2^2 |\s|^2)}
  \right| \leq C \left( \frac{|\xi|^4}{\l_1^4 |\s|^4} +
  \frac{|\xi|^2}{\l_1^2 |\s|^2} \right), \qquad \xi \in \mathcal{B}_1.
\end{equation}
Using equation \eqref{eq:ue2} and some elementary computations we
find
\begin{eqnarray}\label{eq:ue3}
  \int_{\mathcal{B}_1} \frac{x (x - \l_1 \s_1)}{(1+|\xi|^2)^3
  \left(1+\frac{\l_2^2}{\l_1^2}|\xi-\l_1 \s|^2\right)} & = &
  \int_{\mathcal{B}_1} \frac{x (x - \l_1 \s_1)}{(1+|\xi|^2)^3
  (1 + \l_2^2 |\s|^2)} \left( 1 + \frac{2 \frac{\l_2^2}{\l_1}
  \xi \cdot\s - \frac{\l_2^2}{\l_1^2} |\xi|^2}{1 + \l_2^2 |\s|^2}
  \right) \nonumber \\ & + & e(\l_1,\l_2) = \frac{A_0}{2}
  \frac{1}{(1 + \l_2^2 |\s|^2)} - \frac{A_0}{2}
  \frac{2 \s_1^2}{(1+\l_2^2|\s|^2)|\s|^2} + e(\l_1,\l_2)
  \\ \nonumber & = & \frac{A_0}{2} \frac{1}{\l_2^2 |\s|^2} \left( 1 -
  2 \frac{\s_1^2}{|\s|^2} \right) + e(\l_1,\l_2) = \frac{A_0}{2}
  \frac{\s_2^2-\s_1^2}{\l_2^2 |\s|^4} + e(\l_1,\l_2).
\end{eqnarray}
%
%
%
%
%
On the set $\mathcal{B}_2$ we have
$$
|\xi| \geq \l_1 |\s| - |\xi - \l_1 \s| \geq \frac{3}{4} \l_1 |\s|,
$$
and hence we deduce easily
\begin{eqnarray}\label{eq:dede}
  \int_{\mathcal{B}_2} \frac{|\xi| |x - \l_1 \s_1|}{(1+|\xi|^2)^3
  \left(1+\frac{\l_2^2}{\l_1^2}|\xi-\l_1 \s|^2\right)} \leq
  \frac{C}{\l_1^4 |\s|^4} \int_{B_{\frac{\l_1 |\s|}{4}}}
  \frac{1}{1 + \frac{\l_2^2}{\l_1^2} |\xi|^2} \leq C
  \frac{\log \l_1}{\l_1^4 |\s|^4}.
\end{eqnarray}
In $\mathcal{B}_3$ we have $|\xi| \geq \frac{\l_1 |\s|}{4}$ and
$|\xi - \l_1 \s| \geq \frac{\l_1 |\s|}{4}$, and hence
\begin{eqnarray}\label{eq:dede2}
  \int_{\mathcal{B}_3} \frac{|\xi| |x - \l_1 \s_1|}{(1+|\xi|^2)^3
  \left(1+\frac{\l_2^2}{\l_1^2}|\xi-\l_1 \s|^2\right)} \leq
  C |\O_{\l_1,a}| \frac{1}{\l_1^5 |\s|^5} \frac{1}{\l_1 |\s|} \leq
  C \frac{1}{\l_1^4 |\s|^6}.
\end{eqnarray}
Let us now treat the second term in \eqref{eq:c112}. Reasoning as
above we find
\begin{equation}\label{eq:v2}
  32 \int_{\O_{\l_1,a}} \frac{(\l_2^{-1}h_1(\l_1^{-1}\xi + a,b) +
  O(\l_2^{-2}))x}{(1+|\xi|^2)^3} = \frac{1}{\l_1 \l_2}
  \frac{A_0}{2} \frac{\partial h_1}{\partial x}(a,b) + e(\l_1,\l_2).
\end{equation}
Hence, using formulas \eqref{eq:ue3}-\eqref{eq:v2}, we are able to
estimate \eqref{eq:c11}, and we find
\begin{equation}\label{eq:v3}
  2 \int_\O (P \d_2)_1 \, ((\d_1)_x \wedge (\d_1)_y)_1 =
  \frac{16 A_0}{\l_1 \l_2} \left( \frac{\s_1^2 - \s_2^2}{|\s|^4}
  + \frac{\partial h_1}{\partial x}(a,b) \right) + e(\l_1,\l_2).
\end{equation}

\

\noindent {\bf Estimate of \eqref{eq:c12}.} The proofs of the
estimates of this and the remaining terms will only be sketched,
since they are similar to that of \eqref{eq:c11}. Using the usual
change of variables, equation \eqref{eq:c12} becomes
\begin{equation}\label{eq:c122}
  - 32 \frac{\l_2}{\l_1} \int_{\O_{\l_1,a}} \frac{x (y - \l_1
  \s_2)}{(1+|\xi|^2)^3 \left(1+\frac{\l_2^2}{\l_1^2}|\xi-\l_1 \s|^2\right)}
  + 32 \int_{\O_{\l_1,a}} \frac{(\l_2^{-1}h_1(\l_1^{-1}\xi + a,b) +
  O(\l_2^{-2}))y}{(1+|\xi|^2)^3}.
\end{equation}
To treat the first integral in \eqref{eq:c122} we begin by
dividing again $\O_{\l_1,a}$ into the above sets $\mathcal{B}_1$,
$\mathcal{B}_2$, $\mathcal{B}_3$. Reasoning as before and
neglecting the higher-order terms we find
\begin{eqnarray*}
  \frac{\l_2}{\l_1} \int_{\mathcal{B}_1} \frac{x (y - \l_1 \s_2)}{(1+|\xi|^2)^3
  \left(1+\frac{\l_2^2}{\l_1^2}|\xi-\l_1 \s|^2\right)} & = &
  - A_0 \frac{\s_1 \s_2}{\l_1 \l_2 |\s|^4} + e(\l_1,\l_2); \\
  \int_{\O_{\l_1,a}} \frac{(\l_2^{-1}h_1(\l_1^{-1}\xi + a,b) +
  O(\l_2^{-2}))y}{(1+|\xi|^2)^3} & = & \frac{1}{\l_1 \l_2}
  \frac{A_0}{2} \frac{\partial h_1}{\partial y}(a,b) + e(\l_1,\l_2).
\end{eqnarray*}
Hence, using the last two equations we deduce
\begin{equation}\label{eq:v4}
  2 \int_\O (P \d_2)_1 \, ((\d_1)_x \wedge (\d_1)_y)_2 =
  \frac{A_0}{\l_1 \l_2} \left( 32 \frac{\s_1 \s_2}{|\s|^4}
  + 16 \frac{\partial h_1}{\partial y}(a,b) \right) + e(\l_1,\l_2).
\end{equation}

\

\noindent {\bf Estimate of \eqref{eq:c13}.} We turn now to the
term involving the coefficient $r_{13}$, \eqref{eq:c13}, which can
be written as
\begin{equation}\label{eq:c132}
  16 \frac{\l_2}{\l_1} \int_{\O_{\l_1,a}} \frac{(1 - |\xi|^2)
  (x_1 - \l_1 \s_1)}{(1+|\xi|^2)^3 \left( 1 + \frac{\l_2^2}{\l_1^2}
  |\xi - \l_1 \s|^2 \right)} - 16 \int_{\O_{\l_1,a}}
  \frac{(\l_2^{-1}h_1(\l_1^{-1}\xi + a,b) + O(\l_2^{-2}))
  (1 - |\xi|^2)}{(1+|\xi|^2)^3 \left( 1 + \frac{\l_2^2}{\l_1^2}
  |\xi - \l_1 \s|^2 \right)}.
\end{equation}
Using \eqref{eq:id}, \eqref{eq:ue2} and reasoning as in
\eqref{eq:ue3} one finds
$$
  \int_{\mathcal{B}_1} \frac{(1 - |\xi|^2) (x_1 -
  \l_1 \s_1)}{(1+|\xi|^2)^3 \left( 1 + \frac{\l_2^2}{\l_1^2}
  |x - \l_1 \s|^2 \right)} = e(\l_1,\l_2).
$$
Similar estimates hold if one integrates on the sets
$\mathcal{B}_2$ and $\mathcal{B}_3$. Moreover, using \eqref{eq:id}
and elementary computations one finds
$$
  \int_{\O_{\l_1,a}} \frac{(\l_2^{-1}h_1(\l_1^{-1}\xi + a,b) +
  O(\l_2^{-2})) (1 - |\xi|^2)}{(1+|\xi|^2)^3 \left( 1 +
  \frac{\l_2^2}{\l_1^2} |\xi - \l_1 \s|^2 \right)} = e(\l_1,\l_2).
$$
From the last two equations we deduce
\begin{equation}\label{eq:v5}
  2 \int_\O (P \d_2)_1 \, ((\d_1)_x \wedge (\d_1)_y)_3 = e(\l_1,\l_2).
\end{equation}

\

\noindent {\bf Estimate of \eqref{eq:c31}.} The expression in
\eqref{eq:c31} becomes
$$
32 \int_{\O_{\l_1,a}} \left( \frac{1}{1 + \frac{\l_2^2}{\l_1^2}
|\xi - \l_1 \s|^2} - \frac{1}{\l_2^2} h_3(\l_1^{-1}\xi + a,b) +
O(\l_2^{-3}) \right) \frac{x}{(1+|\xi|^2)^3}.
$$
Reasoning as above, we obtain
$$
\int_{\mathcal{B}_1} \left( \frac{2}{1 + \frac{\l_2^2}{\l_1^2}
|\xi - \l_1 \s|^2} - \frac{2}{\l_2^2} h_3(\l_1^{-1}\xi + a,b) +
O(\l_2^{-3}) \right) \frac{x}{(1+|\xi|^2)^3} = e(\l_1,\l_2),
$$
and that the integrals on the sets $\mathcal{B}_2$ and
$\mathcal{B}_3$ are also of order $e(\l_1,\l_2)$. Hence we find
\begin{equation}\label{eq:v6}
  2 \int_\O (P \d_2)_3 \, ((\d_1)_x \wedge (\d_1)_y)_1 = e(\l_1,\l_2).
\end{equation}

\

\noindent {\bf Estimate of \eqref{eq:c33}.} We turn now to the
term involving $r_{33}$. Using the above change of variables,
\eqref{eq:c33} becomes
$$
- 16 \int_{\O_{\l_1,a}} \left( \frac{1}{1 + \frac{\l_2^2}{\l_1^2}
|\xi - \l_1 \s|^2} - \frac{1}{\l_2^2} h_3(\l_1^{-1}\xi + a,b) +
O(\l_2^{-3}) \right) \frac{1 - |\xi|^2}{(1+|\xi|^2)^3}.
$$
Using equation \eqref{eq:id} and reasoning as above one finds
\begin{equation}\label{eq:v7}
  2 \int_\O (P \d_2)_3 \, ((\d_1)_x \wedge (\d_1)_y)_3 = e(\l_1,\l_2).
\end{equation}

\

\noindent {\bf Other estimates.} From the estimates of the terms
\eqref{eq:c11}-\eqref{eq:c33} one can deduce also those involving
the coefficients $r_{22}, r_{12}, r_{23}$ and $r_{32}$. In fact,
it is sufficient to permute the coordinates $x$ and $y$ in a
suitable way. Thus one finds
\begin{eqnarray}\label{eq:oe}
  2 \int_\O (P \d_2)_2 \, ((\d_1)_x \wedge (\d_1)_y)_2 =
  \frac{16 A_0}{\l_1 \l_2} \left( \frac{\s_2^2 - \s_1^2}{|\s|^4}
  + \frac{\partial h_2}{\partial y}(a,b) \right) + e(\l_1,\l_2). \\
  \label{eq:oe2}
  2 \int_\O (P \d_2)_2 \, ((\d_1)_x \wedge (\d_1)_y)_1 =
  \frac{16 A_0}{\l_1 \l_2} \left( 2 \frac{\s_1 \s_2}{|\s|^4}
  +  \frac{\partial h_2}{\partial x}(a,b) \right) + e(\l_1,\l_2). \\
  \label{eq:oe3}
  2 \int_\O (P \d_2)_3 \, ((\d_1)_x \wedge (\d_1)_y)_2 = e(\l_1,\l_2).
  \qquad
  2 \int_\O (P \d_2)_2 \, ((\d_1)_x \wedge (\d_1)_y)_3 = e(\l_1,\l_2).
\end{eqnarray}
Hence the conclusion follows from \eqref{eq:v3}, \eqref{eq:v4},
\eqref{eq:v5}, \eqref{eq:v6}, \eqref{eq:v7}, and
\eqref{eq:oe}-\eqref{eq:oe3}.
\end{pf}

\subsection{Expansion for $k$ bubbles}

In this subsection we consider the case of $k$ masses. When $k =
2$, from Proposition \ref{p:},  and Lemmas \eqref{l:intg}
\ref{l:apex}, \ref{l:apex2} we find
\begin{eqnarray}\label{eq:fe2}
  I_\e(u) & = & \frac{16}{9} A_0 + 8 A_0 \left( \frac{1}{\l_1^2}
\tilde{H}(a) + \frac{1}{\l_2^2} \tilde{H}(b) - \frac{\e}{\l_1}
d_{Id}g (a) - \frac{\e}{\l_2} d_{R^{-1}}g (a) \right) \nonumber \\
& + & \frac{16 A_0}{\l_1 \l_2} \left[  r_{11} \left( \frac{\s_1^2
- \s_2^2}{|\s|^4}
  + \frac{\partial h_1}{\partial x}(a,b) \right) +
  r_{22} \left( \frac{\s_2^2 - \s_1^2}{|\s|^4}
  + \frac{\partial h_2}{\partial y}(a,b) \right) \right. \\
  & + & \left. r_{12} \left( 2 \frac{\s_1 \s_2}{|\s|^4}
  + \frac{\partial h_1}{\partial y}(a,b) \right) +
  r_{21} \left( 2 \frac{\s_1 \s_2}{|\s|^4}
  + \frac{\partial h_2}{\partial x}(a,b) \right) \right]
  \nonumber + e(\e,\l_1) + e(\e,\l_2) + e(\l_1,\l_2).
\end{eqnarray}
where $u = P \d_1 + P \d_2$ and $r_{ij}$ are the entries of the
matrix $R$.

\

\noindent We consider now the more general case of $k$ masses.
Given two bubbles $\d_i = R_i \pi$ and $\d_k = R_k \pi$ (we recall
the definition of the stereographic projection $\pi$ in Section
\ref{s:pre}), where $R_i, R_k \in SO(3)$, we denote by $R_{ik}$
the matrix $R_i^{-1} \circ R_k$. By invariance under rotation, it
is clear that the interaction between $\d_i$ and $\d_k$ is the
same as the interaction between $\pi$ and $R_i^{-1} \d_k =
R_i^{-1} R_k \pi$.

In the expansion of the Euler functional for $k$ masses, since
$I_\e$ is cubic in $u$, we are going to find mixed terms of the
form $\int_\O P \d_i \cdot (P \d_j \wedge P \d_k)$, where $i, j$
and $k$ are all different. Since we are assuming that the distance
of the points $p_i, p_j$ and $p_k$ is uniformly bounded from
below, there holds
\begin{equation}\label{eq:ijk}
  \int_\O P \d_i \cdot (P \d_j \wedge P \d_k) \leq C \frac{1}{\l_i
  \l_j \l_k}, \qquad i \neq j \neq k, i \neq k.
\end{equation}
It follows that the interaction among three distinct bubbles in
$Z$ is negligible with respect to the interactions with $g, \O$
and the interaction between two bubbles.

We recall the definition of the quantity $e(\e,\l_1,\dots,\l_k)$
in Section \ref{s:pre}. Using equation \eqref{eq:ijk}, and
omitting some straightforward but tedious computations we obtain
the following Proposition.

\begin{pro}\label{p:km}
Let $\ov{C} > 0$, let $k \in \N$ and let $Z$ be defined by
\eqref{eq:manz}. For $i \neq j$ let us set
\begin{eqnarray}\label{eq:ijij}
  & & F_\O(p_i,p_j,R_i,R_j) = 16 A_0 \\
  & \times & \nonumber \left[ (R_{ij})_{11} \left(
  \frac{(p_j-p_i)_1^2 - (p_j-p_i)_2^2}{|p_j-p_i|^4}
  + \frac{\partial h_1}{\partial x}(p_i,p_j) \right) +
  (R_{ij})_{22} \left( \frac{(p_j-p_i)_2^2 -
  (p_j-p_i)_1^2}{|p_j-p_i|^4} + \frac{\partial h_2}{\partial y}(p_i,p_j)
  \right) \right. \\ & + & \left. (R_{ij})_{12} \left( 2 \frac{(p_j-p_i)_1
  (p_j-p_i)_2}{|p_j-p_i|^4} + \frac{\partial h_1}{\partial y}(p_i,p_j) \right)
  + (R_{ij})_{21} \left( 2 \frac{(p_j-p_i)_1 (p_j-p_i)_2}{|p_j-p_i|^4}
  + \frac{\partial h_2}{\partial x}(p_i,p_j) \right) \right]
  \nonumber
\end{eqnarray}
and
$$
\Sig_{\O,g} (\e,p_1,\dots,p_k, \l_1,\dots,\l_k,R_1,\dots,R_k) =
\sum_{i=1}^k F_{\O,g} (p_i, \l_i, R_i) + \sum_{i<j}
\frac{F_\O(\e,p_i,p_j,R_i,R_j)}{\l_i \l_j}.
$$
Then there holds
$$
I_\e\left( u \right) = \frac{8 k}{9} A_0 + \Sig_{\O,g} +
e(\e,\l_1,\dots,\l_k);
$$
and
$$
\frac{\partial I_\e(u)}{\partial p_i} = \frac{\partial
\Sig_{\O,g}}{\partial p_i} + e(\e,\l_1,\dots,\l_k); \qquad
\frac{\partial I_\e(u)}{\partial \l_i} = \frac{\partial
\Sig_{\O,g}}{\partial \l_i} + \frac{1}{\l_i}
e(\e,\l_1,\dots,\l_k); $$
$$
\frac{\partial I_\e(u)}{\partial R_i} = \frac{\partial
\Sig_{\O,g}}{\partial R_i} + e(\e,\l_1,\dots,\l_k),
$$
where $u = \sum_{i=1}^k P \, R_i \d_{p_i,\l_i}$.
\end{pro}

\subsection{Some remarks}\label{s:rem}

In this subsection we consider the expansion for $2$ masses with
zero boundary data. Our goal is to extremize the functional in
\eqref{eq:fe2} with respect to $a, b$, $\frac{\l_1}{\l_2}$ and
$R$. Letting $G$ denote the Green's function of $\O$ and setting
$$
G_1 (a,\xi) = \frac{\partial G}{\partial a_1} (a,\xi); \qquad
\qquad G_2 (a,\xi) = \frac{\partial G}{\partial a_2} (a,\xi)
$$
there holds
\begin{eqnarray*}
  \left( \frac{\s_1^2 - \s_2^2}{|\s|^4}
  + \frac{\partial h_1}{\partial x}(a,b) \right) =
  \frac{\partial G_1}{\partial x} (a,b); \qquad
  \left( 2 \frac{\s_1 \s_2}{|\s|^4} + \frac{\partial h_1}{\partial y}(a,b)
  \right) = \frac{\partial G_1}{\partial y} (a,b); \\
  \left( 2 \frac{\s_1 \s_2}{|\s|^4}
  + \frac{\partial h_2}{\partial x}(a,b) \right) =
  \frac{\partial G_2}{\partial x} (a,b); \qquad
  \left( \frac{\s_2^2 - \s_1^2}{|\s|^4}
  + \frac{\partial h_2}{\partial y}(a,b) \right) =
   \frac{\partial G_2}{\partial y} (a,b).
\end{eqnarray*}
Using these expressions, the expansion of $I_\e(u)$, with $u = P
\d_1 + R P \d_2$ becomes
\begin{eqnarray}\label{eq:fe3}
  I_\e(u) & = & \frac{16}{9} A_0 + 8 A_0 \left( \frac{1}{\l_1^2}
  \left( \frac{\partial h_1}{\partial x} + \frac{ \partial
  h_2}{\partial y}\right)(a,a) + \frac{1}{\l_2^2} \left(
  \frac{\partial h_1}{\partial x} + \frac{ \partial h_2}{\partial
  y}\right)(b,b) \right) \nonumber \\ & + & \frac{16 A_0}{\l_1 \l_2}
  \left[  r_{11} \frac{\partial G_1}{\partial x} (a,b) + r_{12}
  \frac{\partial G_1}{\partial y} (a,b) + r_{21} \frac{\partial
  G_2}{\partial x} (a,b) + r_{22} \frac{\partial G_2}{\partial y}
  (a,b) \right] + e(\l_1, \l_2).
  \nonumber
\end{eqnarray}
The entries $r_{i3}$ and $r_{3i}$ of the matrix $R$ appear as
lower order in the above formula, see Remark \ref{r:pos} $(b)$. We
can write
$$
r_{11} \frac{\partial G_1}{\partial x} (a,b) + r_{12}
  \frac{\partial G_1}{\partial y} (a,b) + r_{21} \frac{\partial
  G_2}{\partial x} (a,b) + r_{22} \frac{\partial G_2}{\partial y}
  (a,b) = {\bf e}_1 \cdot R \n G_1 (a,b) + {\bf e}_2 \cdot R \n G_2
  (a,b).
$$
As in \cite{isa1}, Lemma 5.4, the extremization with respect to
$R$ gives
$$
  {\bf e}_1 \cdot R \n G_1 (a,b) + {\bf e}_2 \cdot R \n G_2
  (a,b) = \pm \left( |\n G_1 (a,b)|^2 + |\n G_1 (a,b)|^2 \pm
  2 | \n G_1 (a,b) \wedge \n G_2 (a,b) | \right)^{\frac{1}{2}}.
$$
Hence, setting $\tilde{H} = (h_1)_x + (h_2)_y$, we are left with
\begin{eqnarray*}
  I_\e(u) & = & \frac{16}{9} A_0 + 8 A_0 \left( \frac{1}{\l_1^2}
  \tilde{H}(a) + \frac{1}{\l_2^2} \tilde{H}(b) \right) \nonumber \\
  & \pm & \frac{16 A_0}{\l_1 \l_2}
  \left[ \left( \frac{\partial G_1}{\partial x} \pm
  \frac{\partial G_2}{\partial y} \right)^2 (a,b) +
  \left( \frac{\partial G_2}{\partial x} \mp
  \frac{\partial G_1}{\partial y} \right)^2 (a,b)
  \right]^{\frac{1}{2}} + e(\l_1,\l_2),
\end{eqnarray*}
where the $+$ and $-$ signs inside the square brackets are
opposite (hence there are four different possibilities). $I_\e(u)$
has the form $c + a_{11} \xi_1^2 + a_{22} \xi_2^2 \pm 2 a_{12}
\xi_1 \xi_2$, with $a_{ij} > 0$. Thus if we consider the case $c +
a_{11} \xi_1^2 + a_{22} \xi_2^2 + 2 a_{12} \xi_1 \xi_2$, we notice
that minimizing $\sum a_{ij} \xi_i \xi_j / |\xi|^2$ we necessarily
need to select $\xi = (\xi_1,\xi_2)$ with $\xi_1 \xi_2 \leq 0$ and
so this case does not arise. Thus the only case that remains after
extremizing is
\begin{eqnarray}\label{eq:fe4}
  I_\e(u) & = & \frac{16}{9} A_0 + 8 A_0 \left( \frac{1}{\l_1^2}
  \tilde{H}(a) + \frac{1}{\l_2^2} \tilde{H}(b) \right) \nonumber \\
  & - & \frac{16 A_0}{\l_1 \l_2} \left[ \left(
  \frac{\partial G_1}{\partial x} \pm \frac{\partial G_2}{\partial y}
  \right)^2 (a,b) + \left( \frac{\partial G_2}{\partial x} \mp
  \frac{\partial G_1}{\partial y} \right)^2 (a,b)
  \right]^{\frac{1}{2}} + e(\l_1,\l_2).
\end{eqnarray}

\section{Proof of Theorem \ref{t:in2}}\label{s:bou}

In this Section we prove Theorem \ref{t:in2}. We begin with the
following Lemma, proved in \cite{s} and which follows from
straightforward computations.

\begin{lem}\label{l:gs}
Let $\o \in (0,1)$, and let $a_\o = (\o, 0) \in D$. Define also
$\tilde{g}_\o : \partial D \to \R^3$ as
$$
\tilde{g}_\o (x,y) = \left( \frac{x-\o}{(x-\o)^2 + y^2},
\frac{y}{(x-\o)^2 + y^2}, 0 \right), \qquad (x,y) \in \partial D.
$$
Then, letting $g_\o$ be the harmonic extension on $D$ of
$\tilde{g}_\o$, there holds
\begin{equation}\label{eq:gas}
  g_\o (x,y) = \left( \frac{x-\o(x^2+y^2)}{
  (1-\o x)^2+(\o y)^2}, \frac{y}{(1-\o x)^2+(\o y)^2}, 0 \right),
  \qquad (x,y) \in D,
\end{equation}
and
\begin{equation}\label{eq:dgas}
  |\n g_\o|^2 + 2 | (g_\o)_x \wedge (g_\o)_y |
  = \frac{4}{((1-\o x)^2+\o^2y^2)^2}, \qquad (x,y) \in D.
\end{equation}
\end{lem}

\

\noindent By equation \eqref{eq:b1fin}, the concentration points
of blowing-up solutions as $\e \to 0$ are critical points of the
function $\frac{|\n g|^2 \pm 2 | g_x \wedge g_y |}{\tilde{H}}$. In
the next Lemma we describe the critical points of this function in
the case of $g_\o$.

\begin{lem}\label{l:dhgs}
Let $\o \in (0,1)$, and let $g_\o$ be as in Lemma \ref{l:gs}. Then
one has
$$
  W_\o := \left( \frac{|\n g|^2 \pm 2 |
  g_x \wedge g_y |}{\tilde{H}} \right)^{\frac{1}{2}}
  = \sqrt{2} \frac{(1 - x^2 - y^2)}{(1-\o x)^2+(\o y)^2}.
$$
The point $(\o,0)$ is a non-degenerate global maximum for $W_\o$
and
\begin{equation}\label{eq:gax}
  W_\o(\o,0) =  \frac{\sqrt{2}}{1 - \o^2}.
\end{equation}
The Hessian of $W_\o$ at $(\o,0)$ is given by
\begin{equation}\label{eq:gaxh}
  D^2 W_\o(\o,0) = - 2 \sqrt{2} \begin{pmatrix}
    \frac{1}{(1 - \o^2)^3} & 0 \\
    0 & \frac{1}{(1 - \o^2)^3}
  \end{pmatrix}.
\end{equation}
\end{lem}

\begin{rem}\label{r:ndiso}
From equation \eqref{eq:gaxh}, the fact that $\n g (\o,0) \neq 0$,
and from Theorem D in \cite{isa2} it follows that problem
\eqref{eq:inge} admits a solution concentrating at $(\o,0)$ as $\e
\to 0$. The image of these solutions converges to a sphere of
radius $1$ centered at $(0,0,-1)$, since $\var_{a,\l} \to
(0,0,-1)$ as $\e \to 0$, see \eqref{eq:estfi}.

Note that, from \eqref{eq:gax} and \eqref{eq:gaxh}, $W_\o$ attains
a sharp maximum with highly non-degenerate hessian when $\o$ is
close to $1$. We will use this fact to {\em glue} $k$ single
bubbles showing that, for a suitable boundary datum, the
interaction of this datum with the bubbles is stronger than the
interaction among different bubbles. In the next Lemma we give
quantitative estimates of the gradient of $F_{D,g_\o}$ (see
Proposition \ref{p:}) in a suitable neighborhood of one of its
critical points.
\end{rem}

\noindent It is classical to represent a rotation $R_0 \in SO(3)$
using the Euler angles in the following way
$$
  R_0 = \begin{pmatrix}
    \cos \psi \, \cos \phi - \cos \th \, \sin \phi \, \sin \psi
    & \cos \psi \, \sin \phi + \cos \th \, \cos \phi \, \sin \psi
    & \sin \psi \, \sin \th \\
    - \sin \psi \, \cos \phi - \cos \th \, \sin \phi \, \cos \psi
    & - \sin \psi \, \sin \phi + \cos \th \, \cos \phi \, \cos \psi
    & \cos \psi \, \sin \th \\
    \sin \th \, \sin \phi & - \sin \th \, \cos \phi &
    \cos \th
  \end{pmatrix},
$$
where $\th \in (0,\pi)$, $\psi, \phi \in (0,2 \pi)$. For us it is
convenient to use coordinates different from the Euler angles, in
order to have a smooth parametrization near the identity matrix. A
rotation $R$ will be parameterized as
\begin{equation}\label{eq:A2}
    R^{-1} = \begin{pmatrix}
    \cos \psi \, \cos \phi - \cos \th \, \sin \phi \, \sin \psi
    & \cos \psi \, \sin \phi + \cos \th \, \cos \phi \, \sin \psi
    & \sin \psi \, \sin \th \\
    - \sin \th \, \sin \phi & \sin \th \, \cos \phi &
    - \cos \th \\
    - \sin \psi \, \cos \phi - \cos \th \, \sin \phi \, \cos \psi
    & - \sin \psi \, \sin \phi + \cos \th \, \cos \phi \, \cos \psi
    & \cos \psi \, \sin \th
  \end{pmatrix},
\end{equation}
namely as
$$
R^{-1} =
\begin{pmatrix}
  1 & 0 & 0 \\
  0 & 0 & -1 \\
  0 & 1 & 0
\end{pmatrix} R_0.
$$
$R$ is the identity matrix for $\th = \frac{\pi}{2}$, and $\phi =
\psi = 0$, and the angles $\th, \psi, \phi$ are smooth coordinates
near the identity. In fact there holds
$$
\frac{\partial R^{-1}}{\partial \th} = \begin{pmatrix}
  0 & 0 & 0 \\
  0 & 0 & -1 \\
  0 & 1 & 0
\end{pmatrix}; \qquad \frac{\partial R^{-1}}{\partial \psi} = \begin{pmatrix}
  0 & 0 & 1 \\
  0 & 0 & 0 \\
  -1 & 0 & 0
\end{pmatrix}; \qquad \frac{\partial R^{-1}}{\partial \phi} = \begin{pmatrix}
  0 & 1 & 0 \\
  -1 & 0 & 0 \\
  0 & 0 & 0
\end{pmatrix},
$$
when $\th = \frac{\pi}{2}$, $\psi = \phi = 0$. We will show that
the identity matrix is critical with respect to the rotations for
the quantity $d_{R^{-1}} g_\o (\o,0)$. There holds
\begin{eqnarray*}
  d_{R^{-1}} g_\o & = & (\cos \psi \, \cos \phi
  - \cos \th \, \sin \phi \, \sin \psi) \frac{\partial
  (g_\o)_1}{\partial x} + (\cos \psi \, \sin \phi + \cos \th \, \cos
  \phi \, \sin \psi) \frac{\partial (g_\o)_2}{\partial x} \\
  & - & (\sin \th \, \sin \phi) \frac{\partial
  (g_\o)_1}{\partial y} + (\sin \th \, \cos \phi)
  \frac{\partial (g_\o)_2}{\partial y}.
\end{eqnarray*}
From simple computations one finds
$$
  \frac{\partial (g_\o)_1}{\partial x} = \frac{\partial
  (g_\o)_2}{\partial y} =
  \frac{(1 - \o x)^2 - \o^2 y^2}{((1 - \o x)^2 + \o^2 y^2)^2};
  \qquad \frac{\partial (g_\o)_1}{\partial y} = - \frac{\partial
  (g_\o)_2}{\partial x} = - 2
  \frac{(1 - \o x) \o y}{((1 - \o x)^2 + \o^2 y^2)^2},
$$
and hence
\begin{eqnarray}\label{eq:dg}
  d_{R^{-1}} g_\o & = & ( \cos \psi \, \cos \phi
  - \cos \th \, \sin \phi \, \sin \psi + \sin \th \, \cos \phi )
  \frac{(1 - \o x)^2 - \o^2 y^2}{((1 - \o x)^2 + \o^2 y^2)^2}
  \nonumber \\ & + & 2 ( \cos \psi \, \sin \phi + \cos \th \, \cos
  \phi \, \sin \psi + \sin \th \, \sin \phi )
  \frac{(1 - \o x) \o y}{((1 - \o x)^2 + \o^2 y^2)^2}.
\end{eqnarray}
In the next Lemma we study the critical points of $F_{D,g_\o}$ for
$\xi \sim (\o,0)$, $R \sim Id$, $\l \sim 2 \e^{-1}$ and $\e$
small. We use below the coordinates $\th, \psi, \phi$ in
\eqref{eq:A2} to parametrize the matrix $R$.

\begin{lem}\label{l:deg2}
Let $\o \in (0,1)$ and let $g_\o$ be as above. Then, for fixed
$\e$, the point $x = \o$, $y = 0$, $\l = \frac{2}{\e}$, $\th =
\frac{\pi}{2}$, $\psi = 0$, $\phi = 0$ is critical for $F_{D,g_\o}
(\e,\xi,\l,\th,\psi,\phi)$. For $\mu > 0$ define the set
$$
\mathcal{T}_\mu = \left\{ |x - \o| \leq \mu (1-\o^2), |y| \leq \mu
(1-\o^2), \left| \l - \frac{2}{\e} \right| \leq \frac{\mu}{\e},
\left| \th - \frac{\pi}{2} \right| \leq \mu, |\psi| \leq \mu,
|\phi| \leq \mu \right\}.
$$
Then for $\mu$ sufficiently small and $\o$ sufficiently close to
$1$, there exists a universal constant $C_0$ independent of $\e$,
$\mu$ and $\o$ such that
\begin{equation}\label{eq:grga2}
\n F_{D,g_\o}(\e,\chi) \cdot \chi \geq C_0^{-1} \frac{\e^2
\mu^2}{(1-\o^2)^2} \quad \hbox{ on } \partial \mathcal{T}_\mu,
\quad \hbox{ and hence } \quad \deg(\n F_{D,g_\o},
\mathcal{T}_\mu, 0) = 1,
\end{equation}
where $\chi$ denotes the set of variables $x,y,\l,\th,\psi,\phi$,
and the gradient is taken with respect to $\chi$.
\end{lem}

\begin{pf}
We recall that the functional $F_{D,g_\o}$ is defined by
$$
F_{D,g_\o} (\e,\xi,\l,R) = \frac{1}{\l^2} \left[ \tilde{H}(\xi) -
\e \l d_{R^{-1}} g_\o(\xi) \right], \qquad \hbox{ where }
\tilde{H}(\xi) = \frac{2}{(1-|\xi|^2)^2},
$$
and where $d_{R^{-1}} g_\o(\xi)$ is given by \eqref{eq:dg}. In
particular there holds
$$
\n F_{D,g_\o} \left(\e,\o,0,\frac{2}{\e},\frac{\pi}{2},0,0 \right)
= 0; \qquad  Hess F_{D,g_\o} \left(\e,
\o,0,\frac{2}{\e},\frac{\pi}{2},0,0 \right) = 2
\frac{\e^2}{(1-\o^2)^2} A_{\o,\e},
$$
where
$$
A_{\o,\e} =
\begin{pmatrix}
  \frac{1}{(1-\o^2)} + \frac{3 \o^2}{(1-\o^2)^2} & 0 & -
  \frac{1}{2} \frac{\e \o}{(1-\o^2)} & 0 & 0 & 0 \\
  0 & \frac{1}{(1-\o^2)} + \frac{3 \o^2}{(1-\o^2)^2} & 0
  & 0 & 0 & - \frac{\o}{(1-\o^2)} \\
  - \frac{1}{2} \frac{\e \o}{(1-\o^2)} & 0 & \frac{1}{8}
  \e^2 & 0 & 0 & 0 \\
  0 & 0 & 0 & \frac{1}{4} & 0 & 0 \\
  0 & 0 & 0 & 0 & \frac{1}{4} & 0 \\
  0 & - \frac{\o}{(1-\o^2)} & 0 & 0 & 0 & \frac{1}{2}
\end{pmatrix}.
$$
We point out that the matrix $A_{\o,\e}$ is positive-definite and
non-degenerate. Using simple but tedious computations, one finds
\begin{equation}\label{eq:hessf}
   \left| \n F_{D,g_\o}(\e,\chi) \cdot \chi - A_{\o,\e} \chi \right|
   \leq C \frac{\mu^3 \e^2}{(1-\o^2)^2}, \qquad \hbox{ for }
   x,y,\l,\th,\psi,\phi \in \mathcal{T}_\mu.
\end{equation}
Then the conclusion follows from the fact that
$\o,0,\frac{2}{\e},\frac{\pi}{2},0,0$ is a critical point of
$F_{D,g_\o}$, from \eqref{eq:hessf}, and from the fact that
$A_{\o,\e}$ is positive definite.
\end{pf}

\

\noindent Now we are in position to prove Theorem \ref{t:in2}. For
the main idea see Remark \ref{r:ndiso}.

\

\begin{pfn} {\sc of Theorem \ref{t:in2}}
Let $A = \left\{ S_1 \cup \dots \cup S_k \right\}$ be as in
Theorem \ref{t:in2}, and let $\{ {\bf v}_1, \dots, {\bf v}_k \}
\subseteq \R^3$ denote the centers of $S_1, \dots, S_k$
respectively. Note that, since all the spheres have radius $1$ and
since they all pass through the origin, one has $|{\bf v}_j| = 1$,
for all $j = 1, \dots, k$. Let $\mathcal{R}_1, \dots,
\mathcal{R}_k \in SO(3)$ satisfy
\begin{equation}\label{eq:rk}
  \mathcal{R}_j (0,0,-1) = {\bf v}_j, \qquad \qquad \hbox{ for all }
  j = 1, \dots, k,
\end{equation}
see Remark \ref{r:ndiso}. Let $\o \in (0,1)$ and define
$\tilde{G}_{k,\o} :
\partial D \to \R^3$ by
$$
  \tilde{G}_{k,\o}(x,y) = \sum_{j=1}^k \mathcal{R}_j \tilde{g}_{j,\o},
  \qquad \qquad (x,y) \in \partial D,
$$
where
$$
  \tilde{g}_{j,\o} (x,y) = \tilde{g}_\o \left( \left( \cos
  \frac{2 \pi j}{k} \right) x + \left( \sin \frac{2 \pi j}{k}
  \right) y, - \left( \sin \frac{2 \pi j}{k} \right) x + \left(
  \cos \frac{2 \pi j}{k} \right) y \right).
$$
It is clear that the harmonic extension $G_{k,\o}$ of
$\tilde{G}_{k,\o}$ to the interior of $D$ is given by
\begin{eqnarray}\label{eq:gka}
  G_{k,\o}(x,y) & = & \sum_{j=1}^k \mathcal{R}_i g_{j,\o}(x,y)
  \nonumber \\ & = & \sum_{j=1}^k \mathcal{R}_i g_\o \left( \left(
  \cos \frac{2 \pi j}{k} \right) x + \left( \sin \frac{2 \pi j}{k}
  \right) y, - \left( \sin \frac{2 \pi j}{k} \right) x + \left( \cos
  \frac{2 \pi j}{k} \right) y \right), \quad (x,y) \in D.
\end{eqnarray}
where $g_\o$ is given by \eqref{eq:gas}. Our goal is now to study
the critical points of the functional $\Sig_{D,G_{k,\o}}$ defined
in Proposition \ref{p:km}.

We introduce coordinates $\th_j$, $\psi_j$ and $\phi_j$
parameterizing a rotation $R_j$ (note that this is a generic
rotation, which differs from the fixed rotation $\mathcal{R}_j$)
in the following way
\begin{equation}\label{eq:A3}
  \begin{pmatrix}
    \cos \psi_j \, \cos \phi_j - \cos \th_j \, \sin \phi_j \, \sin \psi_j
    & \cos \psi_j \, \sin \phi_j + \cos \th_j \, \cos \phi_j \, \sin \psi_j
    & \sin \psi_j \, \sin \th_j \\
    - \sin \th_j \, \sin \phi_j & \sin \th_j \, \cos \phi_j &
    - \cos \th_j \\
    - \sin \psi_j \, \cos \phi_j - \cos \th_j \, \sin \phi_j \, \cos \psi_j
    & - \sin \psi_j \, \sin \phi_j + \cos \th_j \, \cos \phi_j \, \cos \psi_j
    & \cos \psi_j \, \sin \th_j
  \end{pmatrix} = R_j^{-1} \mathcal{R}_j.
\end{equation}
The choice of this parametrization will become clear below. Note
that for $\th_j \sim \frac{\pi}{2}$ and $\psi_j, \phi_j$ close to
0, these angles are a smooth parametrization of $SO(3)$ near
$\mathcal{R}_j$. Define also the set
$$
\mathcal{T}^j_\mu = \left\{ \left| (x_j,y_j) - \o \left( \cos
\frac{2 \pi j}{k}, \sin \frac{2 \pi j}{k} \right)\right| \leq \mu
(1-\o^2), \left| \l_j - \frac{2}{\e} \right| \leq \frac{\mu}{\e},
\left| \th_j - \frac{\pi}{2} \right| \leq \mu, |\psi_j| \leq \mu,
|\phi_j| \leq \mu \right\}.
$$
We are going to prove that for $\e$ sufficiently small and for
$\o$ sufficiently close to $1$, the functional $\Sig_{D,G_{k,\o}}$
has a critical point with $x_j,y_j,\l_j,\th_j,\psi_j,\phi_j \in
\mathcal{T}^j_\mu$ for all $j = 1, \dots, k$.


By Proposition \ref{p:km} and by the definition of
$F_{D,G_{k,\o}}$ we have
\begin{eqnarray}\label{eq:l2}
  \Sig_{D,G_{k,\o}} & = & \sum_{j=1}^k F_{D,G_{k,\o}} (\e,p_j, \l_j, R_j) +
  \sum_{l<j} \frac{F_D(p_l,p_j,R_l,R_j)}{\l_l \l_j} =
  F_{D,\mathcal{R}_k g_{k,\o}} (\e,p_k,\l_k,R_k) \nonumber \\ & - &
  8 \frac{\e}{\l_k} A_0 \sum_{j \neq k} d_{R_k^{-1}} g_{\o,j} (p_k) +
  \sum_{j\neq k} F_{D,G_{k,\o}}
  (\e, p_j, \l_j, R_j) + \sum_{l<j} \frac{F_D(p_l,p_j,R_l,R_j)}{\l_l \l_j}.
\end{eqnarray}
By invariance we can write
\begin{eqnarray}\label{eq:l}
  F_{D,\mathcal{R}_k g_{k,\o}} (\e,p_k,\l_k,R_k) & = & \frac{1}{\l_k^2}
  \left[ \tilde{H}(\xi_k) - \e \l_k d_{R_k^{-1}} \mathcal{R}_k
  g_\o(\xi_k) \right] = \frac{1}{\l_k^2}
  \left[ \tilde{H}(\xi_k) - \e \l d_{R_k^{-1}} \mathcal{R}_k
  g_\o(\xi) \right] \nonumber \\ & = & F_{D,g_\o} (\e,p_k,\l_k,
  \mathcal{R}_k^{-1} R_k) = F_{D,g_\o}
  (\e,x_k,y_k,\l_k,\th_k,\psi_k,\phi_k).
\end{eqnarray}
We remark that the function in \eqref{eq:l} is exactly the one
studied in Lemma \ref{l:deg2}. This justifies the choice of the
coordinates $\th_j$, $\psi_j$, $\phi_j$ in \eqref{eq:A3}.

There holds
\begin{eqnarray*}
  \frac{\partial}{\partial x_k} \Sig_{D,G_{k,\o}} & = &
  \frac{\partial}{\partial x_k}
  F_{D,g_\o} (\e,x_k,y_k,\l_k,\th_k,\psi_k,\phi_k) -
  8 \frac{\e}{\l_k} A_0 \sum_{j \neq k}
  \frac{\partial}{\partial x_k}  d_{R_k^{-1}} g_{\o,j} (p_k) \\ &
  + & \sum_{l \neq k} \frac{\partial}{\partial x_k}
  \frac{F_D(p_l,p_k,R_l,R)}{\l_l \l_k}.
\end{eqnarray*}
Since all the $\l_i$'s are of order $\e^{-1}$, and since the
mutual distance between the points $p_j$'s is bounded from below,
it is easy to check that for $\mu$ sufficiently small
$$
(p_i,\l_i,R_i) \in \mathcal{T}^i_\mu \forall \; i \qquad
\Rightarrow \qquad \frac{\e}{\l_k} \left| \frac{\partial}{\partial
x_k} d_{R_k^{-1}} g_{\o,j} (p_k) \right| + \left|
\frac{\partial}{\partial x_k} \frac{F_D(p_j,p_k,R_j,R_k)}{\l_j
\l_k} \right| \leq C \e^2 \quad j \neq k,
$$
where $C$ is a positive constant independent of $\e$, $\o$ and
$\mu$. Hence the last formula and \eqref{eq:l2} imply
$$
\left| \frac{\partial}{\partial x_k} \Sig_{D,G_{k,\o}} -
\frac{\partial}{\partial x_k} F_{D,g_\o}
  (\e,x_k,y_k,\l_k,\th_k,\psi_k,\phi_k) \right| \leq
C \e^2 \qquad \hbox{ if } (p_j,\l_j,R_j) \in \mathcal{T}^j_\mu
\hbox{ for all } j = 1, \dots, k.
$$
Using similar estimates we find
\begin{eqnarray}\label{eq:l4}
  \begin{cases}
  \left| \n_\zeta \Sigma_{D,G_{k,\o}} - \n_\zeta F_{D,g_\o}
  (\e,x_k,y_k,\l_k,\th_k,\psi_k,\phi_k) \right| \leq C \e^2;
  & \\ \left| \frac{\partial}{\partial \l_k} \Sigma_{D,G_{k,\o}}
  - \frac{\partial}{\partial \l_k} F_{D,g_\o}
  (\e,x_k,y_k,\l_k,\th_k,\psi_k,\phi_k) \right| \leq C \e^3, &
  \end{cases}
\end{eqnarray}
provided $(p_j,\l_j,R_j) \in \mathcal{T}^j_\mu$ for all $j = 1,
\dots, k$,  Here $\zeta$ denotes the set of variables $x_k$,
$y_k$, $\th_k, \psi_k, \phi_k$, and where $C$ is a positive
constant independent of $\e$, $\o$ and $\mu$.

Let us fix $\mu$ and $\e$ sufficiently small such that
\eqref{eq:grga2} and \eqref{eq:l4} hold. Then we have
$$
\n_\chi \Sigma_{D,G_{k,\o}} \cdot \chi \geq C_0^{-1} \frac{\e^2
\mu^2}{(1-\o^2)^2} - C \e^2, \quad \chi \in \partial
\mathcal{T}_\mu^k; \qquad \hbox{ if } (p_j,\l_j,R_j) \in
\mathcal{T}^j_\mu \hbox{ for all } j,
$$
where $C$ and $C_0$ are  independent of $\e$, $\o$ and $\mu$. Now,
choosing $\o$ sufficiently close to $1$, depending on $C$, $C_0$
and $\mu$, and reasoning in the same way for the indexes different
from $k$ we obtain
$$
\n_{\chi_j} \Sigma_{D,G_{k,\o}} \cdot \chi_j \geq
\frac{C_0^{-1}}{2} \e^2, \quad \deg(\n_{\chi_j} \Sig_{D,G_{k,\o}},
\mathcal{T}_\mu^j, 0) = 1, \qquad \hbox{ if } (p_j,\l_j,R_j) \in
\mathcal{T}^j_\mu \hbox{ for all } j,
$$
where $\chi_j$ denotes the set of variables $x_j,y_j,\l_j$,
$\th_j,\psi_j,\phi_j$. For the above choices of $\mu$ and $\o$,
let $I_{G_{k,\o(\mu)},\e}$ denote the Euler functional $I_\e$
corresponding to the boundary datum $\tilde{G}_{k,\o(\mu)}$. By
Proposition \ref{p:exp}, for $\e$ sufficiently small we obtain
$$
  \n_{\chi_j} \tilde{I}_{D,G_{k,\o},\e} \cdot \chi_j \geq
\frac{C_0^{-1}}{2} \e^2, \quad \deg(\n_{\chi_j}
\tilde{I}_{D,G_{k,\o},\e}, \mathcal{T}_\mu^j, 0) = 1, \qquad
\hbox{ if } (p_j,\l_j,R_j) \in \mathcal{T}^j_\mu \hbox{ for all }
j.
$$
Then, by Proposition \ref{p:ambbad} and Lemma \ref{l:li} below,
letting $\e \to 0$, we find a family of solutions $u_{\e,\mu}$ of
$I'_{G_{k,\o(\mu)},\e} = 0$ satisfying, up to a subsequence
$$
u_{\e,\mu}(D) \to A_\mu = \left\{ S_{1,\mu}, \dots, S_{k,\mu}
\right\} \hbox{ in the Hausdorff sense } \qquad \hbox{ as } \e \to
0,
$$
where $S_{1,\mu}, \dots, S_{k,\mu}$ are spheres of radius $1$
passing through the origin and lying in a neighborhood of order
$\mu$ of $S_{1}, \dots, S_{k}$ respectively. Now we can choose
$\mu(\e) \to 0$ sufficiently small as $\e \to 0$, and find a
corresponding $\o(\e) \to 1$ such that the solution
$u_{\e,\mu(\e)}$ of $I'_{G_{k,\o(\mu(\e))},\e} = 0$ obtained with
the above method satisfies
$$
u_{\e,\mu}(D) \to A = \left\{ S_1, \dots, S_k \right\} \hbox{ in
the Hausdorff sense } \qquad \hbox{ as } \e \to 0.
$$
This concludes the proof of the Theorem.
\end{pfn}

\begin{lem}\label{l:li}
Let $\tilde{g} : \partial \O \to \R^3$ be a smooth function, let
$k \in \N$, $\ov{C} > 0$, and let $Z$ be defined as in
\eqref{eq:manz}. Let $u$ be a solution of \eqref{eq:pe} of the
form
$$
u = \sum_{i = 1}^k P R_i \d_{p_i,\l_i} + w; \qquad \hbox{ with }
\sum_{i = 1}^k P R_i \d_{p_i,\l_i} \in Z, \hbox{ and }
\|w\|_{H^1_0(\O)} \to 0 \hbox{ as } \e \to 0.
$$
Then $\|w\|_{L^\infty(\O)} \to 0$ as $\e \to 0$.
\end{lem}

\begin{pf} In the following we simply write $\d_i$ for $R_i \d_{p_i,\l_i}$,
and we let $\var_i$ be the function in \eqref{eq:fi} corresponding
to $\d_i$. The function $w$ satisfies
$$
  \begin{cases}
   \D w = 2 \left( \sum_i (\d_i - \var_i) + w + \e g \right)_x
   \wedge \left( \sum_j (\d_j - \var_j) + w + \e g \right)_y -
   \sum_i (\d_i)_x \wedge (\d_i)_y, & \text{ in } \O, \\
   w = 0 \text{ on } \partial \O. &
  \end{cases}
$$
where $g$, as before, denotes the harmonic extension of
$\tilde{g}$ to $\O$. Expanding the wedge produce on the right-hand
side we obtain (as before $P \d_i = \d_i - \var_i$)
\begin{eqnarray}\label{eq:numb}
  \D w & = & 2 \sum_{i \neq j} (P \d_i)_x \wedge (P \d_j)_y + 2
  \sum_i \left[ (P \d_i)_x \wedge w_y + w_x \wedge (P \d_i)_y
  \right] + 2 \e \sum_i \left[ (P \d_i)_x \wedge g_y + g_x \wedge
  (P \d_i)_y \right] \nonumber \\ & - & 2 \sum_i \left[ (\d_i)_x
  \wedge (\var_i)_y + (\var_i)_x \wedge (\d_i)_y \right] + \sum_i
  (\var_i)_x \wedge (\var_i)_y + 2 w_x \wedge w_y + 2 \e (w_x \wedge
  g_y + g_x \wedge w_y) \\ & + & \e^2 g_x \wedge g_y. \nonumber
\end{eqnarray}
Using \eqref{eq:est}, \eqref{eq:est2} and some elementary
computations, for any $p > 1$ the first term in the right hand
side can be estimated in the following way
$$
\|(P \d_i)_x \wedge (P \d_j)_y\|_{L^p(\O)} \leq C(\ov{C}, p)
\left( \frac{1}{\l_i\l_j} + \frac{\l_j^{\frac{p-2}{p}}}{\l_i} +
\frac{\l_i^{\frac{p-2}{p}}}{\l_j} \right) \leq C(\ov{C}, p) \,
\e^{\frac{2}{p}}, \quad i \neq j.
$$
From standard elliptic estimates it follows that
$$
\left\| (\D)^{-1} \left( (P \d_i)_x \wedge (P \d_j)_y \right)
\right\|_{L^\infty(\O)} \leq C(\ov{C}, p) \, \e^{\frac{2}{p}},
\quad i \neq j,
$$
where $(\D)^{-1}$ denotes the Green's operator for $\D$ in $\O$
with Dirichlet boundary conditions. Let us focus now on the second
term in \eqref{eq:numb}. Writing for brevity $\psi = P \d_i$, one
has
$$
\left( \psi_x \wedge w_y + w_x \wedge \psi_y \right) = [
J(w_2,\psi_3) + J(\psi_2,w_3) ] {\bf i} + [ J(w_1,\psi_3) +
J(\psi_1,w_3) ] {\bf j} + [ J(w_1,\psi_2) + J(\psi_1,w_2) ] {\bf
k},
$$
where $J(F,G) = F_x G_y - F_y G_x$ is the Jacobian function. By
the result in \cite{cl} there holds
\begin{equation}\label{eq:tilde}
  \left\| (\D)^{-1} \left( \psi_x \wedge w_y + w_x \wedge \psi_y
\right) \right\|_{L^\infty(\O)} \leq C \| P \d_i \|_{H^1_0(\O)} \,
\| w \|_{H^1_0(\O)} \to 0 \hbox{ as } \e \to 0.
\end{equation}
The remaining terms in \eqref{eq:numb} can be estimated as in
\eqref{eq:tilde}.
\end{pf}

\begin{rem}\label{r:fixg}
With an easy modification of the above arguments we can easily
obtain the limit configuration $\left\{ S_1, \dots, S_k \right\}$
with a boundary datum of the form $\e \tilde{G}$, for some fixed
function $\tilde{G}$ on $\partial D$ independent of $\e$.
\end{rem}

\begin{rem}\label{r:cm}
We remark that to obtain $L^\infty$ estimates on the solutions of
\eqref{eq:pe}, we use in a crucial way that these solutions
satisfy the $H$-surface equation with $H \equiv constant$. Such
estimates are not available for general Palais-Smale sequences, as
exhibited in \cite{cm2}.
\end{rem}

\begin{figure}
\centerline{{\psfig{figure=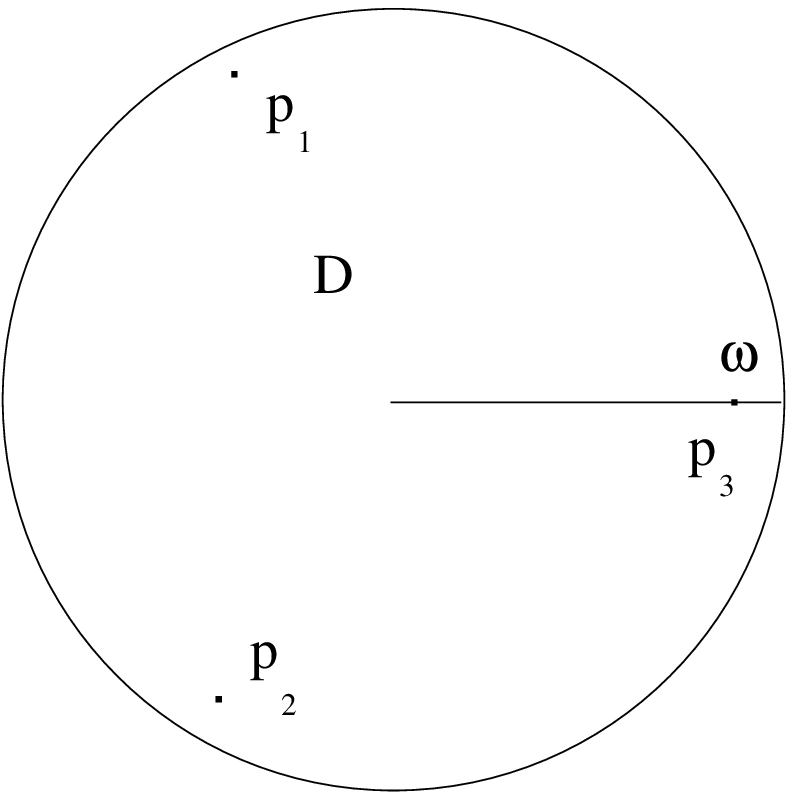,width=6.5cm}} \qquad
{\psfig{figure=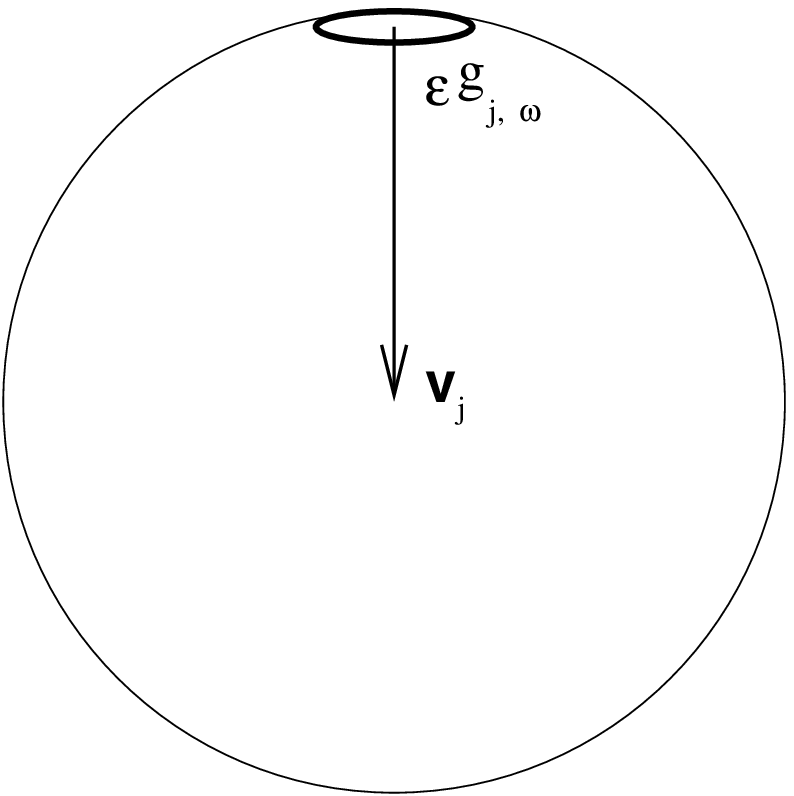,width=6.5cm}}}

\caption{the points $p_j$ on the disk $D$ and the bubble generated
by $g_{\o,j}$} \label{nomelabel3}
\end{figure}

\begin{figure}
\centerline{{\psfig{figure=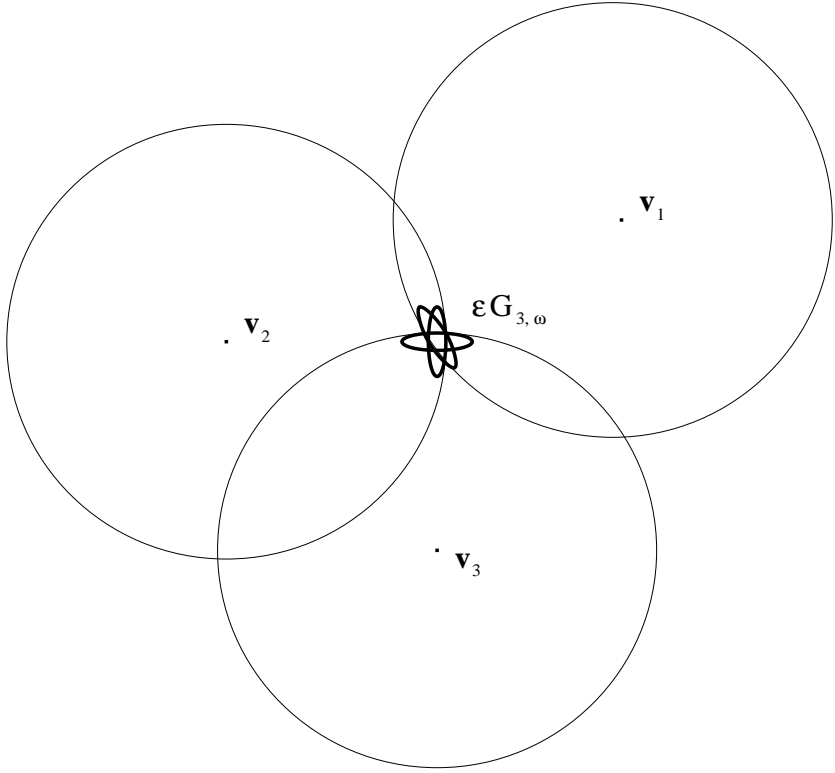,width=11cm}} }

\caption{the boundary datum $\e G_{k,\o}$ and the corresponding
configuration of spheres} \label{nomelabel4}
\end{figure}

\noindent In Figure \ref{nomelabel3} we indicate the location of
the points $p_j$ in $D$ when $\o$ is close to $1$, see the
definition of $\mathcal{T}^j_\mu$. We also plot the boundary datum
$\e g_{\o,j}$, which lies in a plane, and the corresponding bubble
(as in Remark \ref{r:ndiso}) whose {\em center} ${\bf v}_j$ which
is, roughly, perpendicular to the plane of $g_{\o,j}$. We note
that the image of $\tilde{g}_\o$ is a great circle (the Kelvin
inversion of $\partial D$ w.r.t. the point $(\o,0)$). In Figure
\ref{nomelabel4} we plot the configuration of bubbles generated by
the function $\e G_{k,\o}$. Each bubble is nearly perpendicular to
some $g_j,\o$ (whose sum is $G_{k,\o}$).

\section{Appendix}

\noindent This Appendix is devoted to the characterization of the
solutions of the equation $\ov{I}''(\d_{a,\l})[w] = 0$. We can
suppose by invariance that $a = 0$ and $\l = 1$, and we set
$\d_{a,\l} = \pi$, see Section \ref{s:pre}. If a function $w \in
\mathcal{D}$ satisfies $\ov{I}''(\d)[w] = 0$, then it solves the
linearization of \eqref{eq:pinf}, namely
\begin{equation}\label{eq:lin}
  \D w = 2 \left( w_x \wedge \d_y + \d_x \wedge w_y \right),
  \qquad \hbox{ in } \R^2, \qquad w \in \mathcal{D}.
\end{equation}
After inverse stereographic projection, equation \eqref{eq:lin}
can be equivalently viewed on $S^2$ as follows
\begin{equation}\label{eq:lins}
  \D_{g_0} w = 2 (\sin \var)^{-1}
  \left( w_\th \wedge \d_\var + \d_\th \wedge w_\var \right),
  \qquad \hbox{ in } S^2, \qquad w \in H^1(S^2;\R^3).
\end{equation}
where $(\th, \var)$, $0 \leq \th \leq 2 \pi$, $0 \leq \var \leq
\pi$, are spherical coordinates on $S^2$ and $\D_{g_0}$ is the
Laplacian with respect to standard metric on $S^2$.

To analyze \eqref{eq:lins} we shall use some properties of
spherical harmonics that we now recall. Let $P_n(x)$, $x \in
(-1,1)$, denote the $n$-th Legendre function. We define the
associated Legendre function $P^k_n(x)$ by
\begin{equation}\label{eq:a3}
  P^k_n(x) = (-1)^k (1-x^2)^{\frac{k}{2}} \frac{d^k}{d x^k}
  P_n(x), \qquad k \geq 0.
\end{equation}
The spherical harmonics are defined by
\begin{equation}\label{eq:a4}
  Y_{n,k}(\th,\var) = c_{n,|k|} P^{|k|}_n (\cos \var) e^{i k \th},
  \qquad - n \leq k \leq n,
\end{equation}
where the normalization constant $c_{n,|k|}$, see \cite{erd}
equation (21), p.171, is given by
\begin{equation}\label{eq:a5}
  c_{n,|k|} = \left( \frac{2n+1}{4\pi} \right)^{\frac{1}{2}}
  \sqrt{\frac{(n-|k|)!}{(n+|k|)!}}.
\end{equation}
From \eqref{eq:a3} we easily have, by differentiation
\begin{equation}\label{eq:sh3}
P^{k+1}_n(\cos \var) = - \sin \var (P^{k}_n)' (\cos \var) - m \cot
\var P^k_n(\cos \var),
\end{equation}
and from equation (41), p.107 in \cite{hob},
\begin{equation}\label{eq:sh4}
P^{k+2}_n(\cos \var) + 2(m+1) \cot \var P^{k+1}_n (\cos \var) +
(n-k) (n+k+1) P^k_n(\cos \var) = 0.
\end{equation}
From \eqref{eq:a5} we also have, for $0 \leq k \leq n$, $n \geq
2$,
\begin{equation}\label{eq:a8}
  d_{n,k} := \frac{c_{n,k}}{c_{n,k+1}} \leq \sqrt{3/2} n, \qquad
  e_{n,k} := (n-k) (n+k+1) \frac{c_{n,k}}{c_{n,k-1}} \leq \sqrt{3/2}
  n.
\end{equation}
In the sequel we will often write $P^k_n(\var)$ for $P^k_n(\cos
\var)$.

We define the finite-dimensional subspace $\mathcal{H}_n$ to be
the linear span of $Y_{n,k}$, $-n \leq k \leq n$, and
$\mathcal{H}_n$ the subspace of $L^2(S^2;\R^3)$ consisting of
vectors $w = (f_1,f_2,f_3)$ with $f_i \in H_n$, $i = 1, 2, 3$.
Recall that any function $w \in L^2(S^2;\R^3)$ can be decomposed
orthonormally as $w = \sum_{n=0}^\infty w_n$, with $w_n \in
\mathcal{H}_n$. We have

\begin{lem}\label{l:la}
Let $w : S^2 \to \R^3$ be a solution of \eqref{eq:lins}, and let
$w = \sum_{n=0}^\infty w_n$, with $w_n \in \mathcal{H}_n$. Then
$w_n = 0$ for $n \geq 4$.
\end{lem}

\begin{pf}
We first claim that for any $F \in \mathcal{H}_n$ there holds
\begin{equation}\label{eq:a9}
  \G (F) := \D_{g_0} F - \frac{2}{\sin \var} \left( F_\th \wedge
  \d_\var + \d_\th \wedge F_\var \right) \in \mathcal{H}_n.
\end{equation}
If our claim is verified, to prove the Lemma it will be enough to
pick a solution $w$ to \eqref{eq:lins} in $\mathcal{H}_n$ and to
show that $w_n = 0$ for $n \geq 4$.

We now prove our claim. W.l.o.g. pick $F \in \mathcal{H}_n$ of the
form
\begin{equation*}
  F = \left( \a_{k_1} c_{n,k_1} P^{k_1}_n (\var) e^{i k_1 \th},
  \b_{k_2} c_{n,k_2} P^{k_2}_n (\var) e^{i k_2 \th},
  \g_{k_3} c_{n,k_3} P^{k_3}_n (\var) e^{i k_3 \th} \right)
  = \left( \a_{k_1} Y_{n,k_1}, \b_{k_2} Y_{n,k_2},
  \g_{k_3} Y_{n,k_3} \right)
\end{equation*}
Next we have $\d_\var = (\cos \var \cos \th, \cos \var \sin \th, -
\sin \var)$ and $\d_\th = (- \sin \var \sin \th, \sin \var \cos
\th, 0)$. We will show that $\G(F) = - n(n+1)F - 2 {\bf v}$, where
\begin{equation}\label{eq:a10}
  {\bf v} = \begin{pmatrix}
    - i k_2 \b_{k_2} Y_{n,k_2} + 1/2 d_{n,k_3} \g_{k_3} Y_{n,k_3+1} -
    1/2 e_{n,k_3} \g_{k_3} Y_{n,k_3-1} \\
    i k_1 \a_{k_1} Y_{n,k_1} - i/2 d_{n,k_3} \g_{k_3} Y_{n,k_3+1} -
    i/2 e_{n,k_3} \g_{k_3} Y_{n,k_3-1} \\
    i/2 d_{n,k_2} \b_{k_2} Y_{n,k_2+1} - i/2 e_{n,k_2} \b_{k_2} Y_{n,k_2-1}
    - 1/2 d_{n,k_1} \a_{k_1} Y_{n,k_1+1} + 1/2 e_{n,k_1} \a_{k_1} Y_{n,k_1-1}
    \ \end{pmatrix}.
\end{equation}
Since ${\bf v} \in \mathcal{H}_n$, our claim follows. It is
evident that $\D_{g_0} F = - n (n+1) F$, thus it is enough to show
that the second expression in \eqref{eq:a9} is $2 {\bf v}$. This
follows by noting that
\begin{eqnarray*}
  F_\th \wedge \d_\var & = & \begin{pmatrix}
    - i k_2 \b_{k_2} c_{n,k_2} \sin \var P_n^{k_2} (\var) e^{i k_2 \th} -
    i k_3 g_{k_3} c_{n,k_3} \cos \var \sin \th P_n^{k_3} (\var) e^{i k_3 \th}
    \\ i k_3 g_{k_3} c_{n,k_3} \cos \var \sin \th P_n^{k_3} (\var) e^{i k_3
    \th} + i k_1 \a_{k_1} c_{n,k_1} \sin \var P_n^{k_1} (\var)
    e^{i k_1 \th} \\ i k_1 \a_{k_1} c_{n,k_1} \cos \var \sin \th
    P_n^{k_1} (\var) e^{i k_1 \th} - i k_2 \b_{k_2} c_{n,k_2} \sin \var
    P_n^{k_2} (\var) e^{i k_2 \th} \
  \end{pmatrix}; \\
  \d_\th \wedge F_\var & = & \begin{pmatrix}
    -g_{k_3} c_{n,k_3} \sin^2 \var \cos \th (P_n^{k_3})'(\var) e^{i k_3
    \th} \\ - \g_{k_3} c_{n,k_3} \sin^2 \var \sin \th (P_n^{k_3})'(\var)
    e^{i k_3 \th} \\ \b_{k_2} c_{n,k_2} \sin^2 \var \sin \th
    (P_n^{k_2})'(\var) e^{i k_2 \th} + \a_{k_1} c_{n,k_1} \sin^2
    \var \cos \th (P_n^{k_1})'(\var) e^{i k_1 \th} \
  \end{pmatrix}.
\end{eqnarray*}
Using \eqref{eq:sh3} and \eqref{eq:sh4} it is easily verified that
$(\sin \var)^{-1} \left[ F_\th \wedge \d_\var + \d_\th \wedge
F_\var \right] = {\bf v}$. Then from \eqref{eq:a10} and
integrating the equation $\G(w) \cdot \ov{w} = 0$ on $S^2$ we find
\begin{equation}\label{eq:a11}
  - n (n+1) \begin{pmatrix}
    A^2 \\
    B^2 \\
    C^2 \
  \end{pmatrix} = - 2 \sum_k \begin{pmatrix}
    - i k \b_k \ov{\a}_k + 1/2 d_{n,k} \g_k \ov{\a}_{k+1} -
    1/2 e_{n,k} \g_k \ov{\a}_{k-1} \\
    i k \a_k \ov{\b}_k - i/2 d_{n,k} \g_k \ov{\b}_{k+1} - i/2 e_{n,k} \g_k
    \ov{\b}_{k-1} \\ i/2 d_{n,k} \b_k \ov{\g}_{k+1} - i/2 e_{n,k}
    \b_k \ov{\g}_{k-1} - 1/2 d_{n,k} \a_k \ov{\g}_{k+1} + 1/2
    e_{n,k} \a_k \ov{\g}_{k-1} \
  \end{pmatrix},
\end{equation}
where $w = \sum_k (\a_k Y_{n,k}, \b_k Y_{n,k}, \g_k Y_{n,k})$ and
where $A^2 = \sum |\a_k|^2$, $B^2 = \sum |\b_k|^2$, $C^2 = \sum
|\g_k|^2$. Using \eqref{eq:a8}, \eqref{eq:a11} and the
Cauchy-Schwartz inequality we find
\begin{eqnarray*}
  n(n+1) |(A^2,B^2,C^2)| & \leq & 2n \left| \left( AB + \sqrt{3/2} AC,
  AB + \sqrt{3/2} BC \sqrt{3/2} (BC + AC) \right) \right| \\ &
  \leq & 2 n \left( \frac{9+\sqrt{6}}{2} \right)^{\frac{1}{2}}
  \left( A^4 + B^4 + C^4 \right)^{1/2}.
\end{eqnarray*}
Thus $n + 1 \leq \left( 18 + \sqrt{24} \right)^{\frac{1}{2}}$,
which implies $n \leq 3$.
\end{pf}

\begin{lem}\label{l:pol}
The solutions of equation \eqref{eq:lins} are of the form
$$
w = c + \begin{pmatrix}
  \a x_2 + \b x_3 \\ - \a x_1 + \g x_3 \\ - \b x_1 - \g x_2
\end{pmatrix} + (\a' x_1 + \b' x_2 + \g' x_3) \begin{pmatrix}
  x_1 \\ x_2 \\ x_3
\end{pmatrix}
$$
where $c \in \R^3$ and $\a, \b, \g$, $\a', \b', \g' \in \R$ are
arbitrary.
\end{lem}

\begin{pf}
We denote by $J(w)$ the r.h.s. of \eqref{eq:lins}. From the proof
of Lemma \ref{l:la} it follows that $J$ preserves the degree of
spherical harmonic functions. Equivalently, $J$ preserves the
degree of polynomial functions in $\R^3$ restricted to $S^2$. By
this reason and by Lemma \ref{l:la}, we can confine ourselves to
study $J$ just on polynomials of order $1, 2$ and $3$. Since the
computations involved in the proof are straightforward, we just
give a simple sketch below, omitting some details.

Using simple computations, we obtain
\begin{equation}\label{eq:f1}
  J (x_1,0,0) = (0,2x_2,2x_3); \qquad
  J (x_2,0,0) = (0,-2x_1,0); \qquad J (x_3,0,0) = (0,0,-2x_1).
\end{equation}
With a permutation of coordinates one also finds
\begin{equation}\label{eq:f2}
  J (0,x_2,0) = (2x_1,0,2x_3); \qquad
  J (0,x_3,0) = (0,0,-2x_2); \qquad F (0,x_1,0) = (-2x_2,0,0).
\end{equation}
\begin{equation}\label{eq:f3}
  F (0,0,x_3) = (2x_1,2x_2,0); \qquad
  F (0,0,x_1) = (-2x_3,0,0); \qquad F (0,0,x_2) = (0,-2x_3,0).
\end{equation}
Hence, letting $w = (a_1 x_1 + a_2 x_2 + a_3 x_3,
    b_1 x_1 + b_2 x_2 + b_3 x_3, c_1 x_1 + c_2 x_2 + c_3 x_3)$ we
    find
$$
  J (w) = 2 \begin{pmatrix}
    - b_1 x_2 + b_2 x_1 - c_1 x_3 + c_3 x_1 \\
    a_1 x_2 - a_2 x_1 - c_2 x_3 + c_3 x_2 \\
    a_1 x_3 - a_3 x_1 + b_2 x_3 - b_3 x_2
  \end{pmatrix}; \qquad \D w = - 2 \begin{pmatrix}
  a_1 x_1 + a_2 x_2 + a_3 x_3 \\
    b_1 x_1 + b_2 x_2 + b_3 x_3 \\
    c_1 x_1 + c_2 x_2 + c_3 x_3 \
  \end{pmatrix},
$$
The system of equations $\D w = J(w)$
admits the following solutions
\begin{equation}\label{eq:f6}
  \begin{pmatrix}
    a_1 & a_2 & a_3 \\
    b_1 & b_2 & b_3 \\
    c_1 & c_2 & c_3 \
  \end{pmatrix} =
  \begin{pmatrix}
    0 & \a & \b \\
    - \a & 0 & \g \\
    - \b & - \g & 0 \
  \end{pmatrix},
\end{equation}
with $\a, \b$, $\g$ arbitrary real numbers.

Let us now consider the homogeneous second order polynomials. We
have $\D x_i^2 = 2 (1 - 3 x_i^2)$ and $\D (x_i x_j) = - 6 x_i
x_j$. Using the Leibnitz rule and \eqref{eq:f1}-\eqref{eq:f3}, we
can compute $J(w)$ when $w$ has the form
$$
w = \begin{pmatrix}
  a_1 x_1^2 + a_2 x_1 x_2 + a_3 x_2^2 + a_4 x_1 x_3 + a_5 x_2 x_3 + a_6 x_3^2
\\
  b_1 x_1^2 + b_2 x_1 x_2 + b_3 x_2^2 + b_4 x_1 x_3 + b_5 x_2 x_3 + b_6 x_3^2
\\
  c_1 x_1^2 + c_2 x_1 x_2 + c_3 x_2^2 + c_4 x_1 x_3 + c_5 x_2 x_3 + c_6 x_3^2
\end{pmatrix}
$$
From the relation $\D w = J(w)$, and using elementary computations
we obtain
$$
w = (\a x_1 + \b x_2 + \g x_3) + (\d, \eta, \sigma),
$$
where $\a, \b, \g$, $\d, \eta, \sigma$ are arbitrary real numbers.

Let us now turn to the third order polynomials. We have
$$
\D x_i^3 = 6 x_i (1 - 2 x_1^2); \qquad \D (x_i^2 x_j) = 2 x_j (1 -
6 x_i^2); \qquad \D (x_i x_j x_k) = - 12 x_i x_j x_k.
$$
Again, the values of $F$ on the third order polynomials can be
computed with the Leibnitz rule and \eqref{eq:f1}-\eqref{eq:f3}.
Letting
\begin{eqnarray*}
w = \begin{pmatrix} a_1 x_1^3 + a_2 x_1^2 x_2 + a_3 x_1^2 x_3 +
a_4 x_1 x_2 x_3 + a_5 x_1 x_2^2 + a_6 x_1 x_3^2 + a_7 x_2^3 + a_8
x_3^3
+ a_9 x_2^2 x_3 + a_{10} x_2 x_3^2 \\
b_1 x_1^3 + b_2 x_1^2 x_2 + b_3 x_1^2 x_3 + b_4 x_1 x_2 x_3 + b_5
x_1 x_2^2 + b_6 x_1 x_3^2 + b_7 x_2^3 + b_8 x_3^3 + b_9 x_2^2 x_3
+ b_{10} x_2 x_3^2
\\
c_1 x_1^3 + c_2 x_1^2 x_2 + c_3 x_1^2 x_3 + c_4 x_1 x_2 x_3 + c_5
x_1 x_2^2 + c_6 x_1 x_3^2 + c_7 x_2^3 + c_8 x_3^3 + c_9 x_2^2 x_3
+ c_{10} x_2 x_3^2
\end{pmatrix}
\end{eqnarray*}
Hence, equating the coefficients in the above two expressions we
find the following system
and considering the equation $\D w = J(w)$, we find a system
decoupled in four parts. The first part consists of seven
equations involving the seven terms $a_2, a_7, a_{10}$, $b_1, b_5,
b_6$ and $c_4$. Using simple computations one finds $a_2 = a_7 =
a_{10} = \a$, $b_1 = b_5 = b_6 = -\a$, $c_4 = 0$, for some $\a \in
\R$.

The second part consists of seven equations involving the seven
terms $a_3, a_8, a_9$, $b_4$ and $c_1, c_5, c_6$. Using simple
computations one finds $a_3 = a_8 = a_9 = \b$, $c_1 = c_5 = c_6 =
-\b$, $b_4 = 0$, for some $\b \in \R$.

The third part consists of seven equations involving the seven
terms $a_4$, $b_3, b_8, b_9$, and $c_2, c_7, c_{10}$. Using simple
computations one finds $b_3 = b_8 = b_9 = \g$, $c_2 = c_7 = c_{10}
= -\g$, $a_4 = 0$, for some $\g \in \R$.

The fourth part consists in nine equations involving the terms
$a_1, a_5, a_6$, $b_2, b_7, b_{10}$ and $c_3, c_8, c_9$. Using
simple computations one finds $a_1 = a_5 = a_6 = b_2 = b_7 =
b_{10} = c_3 = c_8 = c_9 = 0$.

The solution obtained in this way represent just the linear
functions in \eqref{eq:f6}, taking into account of the identity
$x_1^2 + x_2^2 + x_3^2 = 1$ on $S^2$. This concludes the proof.
\end{pf}

\begin{pfn} {\sc of Proposition \ref{p:up}.} Coming back to the
space $\mathcal{D}$, and using some elementary computation, the
proof of the last statement follows immediately from Lemma
\ref{l:pol}. The first inequality is immediate to check. The
second inequality follows from the proof of Lemma \ref{l:pol} when
$v \in \oplus_{n \geq 4} \mathcal{H}_n$. When $v$ has some
non-zero components in $\oplus_{n \leq 3} \mathcal{H}_n$, then it
is sufficient to use straightforward computations, since we have
to deal with finite combinations of spherical harmonics.
Alternatively note that, since $\d$ is a mountain-pass critical
point of $\ov{I}$, the linearized operator possesses only one
negative eigenvalue (with corresponding eigenvector $\d$), hence
if $v \perp \d$ and $v \perp Ker \ov{I}''(\d)$, $v$ must be a
combination of positive eigenvectors.
\end{pfn}


\begin{thebibliography}{99}



\bibitem{ab} Ambrosetti A., Badiale
M., Homoclinics: Poincar\'e-Melnikov type results via a
variational approach, Ann. Inst. Henri. Poincar\'e Analyse Non
Lin\'eaire 15 (1998), 233-252. Preliminary note in C. R. Acad.
Sci. Paris 323, S\'erie I (1996), 753-758.




\bibitem{bab} Bahri A. Critical points at infinity in some variational
problems. Pitman Research Notes in Mathematics Series, 182.
Longman Scientific \& Technical, Harlow; copublished in the United
States with John Wiley \& Sons, Inc., New York, (1989). vi+I15+307
pp.


\bibitem{bach} Bahri A., Chanillo S., The difference of
topology at infinity in changing-sign Yamabe problems on $S\sp 3$
(the case of two masses). Comm. Pure Appl. Math. 54 (2001), no. 4,
450--478.


\bibitem{bc} Bahri A., Coron J.M. On a nonlinear elliptic equation
involving the critical Sobolev exponent: the effect of the
topology of the domain. Comm. Pure Appl. Math. 41 (1988), no. 3,
253--294.



\bibitem{blr} Bahri A., Li Y.Y., Rey O.,
On a variational problem with lack of compactness: the topological
effect of the critical points at infinity. Calc. Var. Partial
Differential Equations 3 (1995), no. 1, 67--93.


\bibitem{baf}  Bandle C.; Flucher M., Harmonic radius and
concentration of energy; hyperbolic radius and Liouville's
equations $\Delta U=e\sp U$ and $\Delta U=U\sp {(n+2)/(n-2)}$.
SIAM Rev. 38 (1996), no. 2, 191--238.


\bibitem{bk}  Baraket S., Pacard F., Construction of singular
limits for a semilinear elliptic equation in dimension $2$. Calc.
Var. Partial Differential Equations 6 (1998), no. 1, 1--38.



%
%
%


\bibitem{brc} Brezis H., Coron J.M., Multiple solutions
of $H$-systems and Rellich's conjecture. Comm. Pure Appl. Math.
37 (1984), no. 2, 149--187.


\bibitem{brc2} Brezis H., Coron J.M.,
Convergence of solutions of $H$-systems or how
to blow bubbles. Arch. Rational Mech. Anal. 89 (1985), no. 1, 21--56.



\bibitem{cm} Caldiroli P., Musina R., Existence of minimal H-bubbles,
to appear in Comm. Contemp. Math.

\bibitem{cm2} Caldiroli P., Musina R., A remark on Palais-Smale sequences
for $H$-systems, unpublished paper.

\bibitem{cl} Chanillo S., Li Y.Y., Continuity of solutions of
uniformly elliptic equations in $R\sp 2$. Manuscripta Math. 77
(1992), no. 4, 415--433.

%


\bibitem{erd} Erdelyi A., Magnus W., Oberhettinger F., Tricomi F.G.
Higher transcendental functions. Vols. I, II. Based, in part, on
notes left by Harry Bateman. McGraw-Hill Book Company, Inc., New
York-Toronto-London, (1953). xxvi+302, xvii+396 pp.

\bibitem{fl} Flucher M., Extremal functions for the Trudinger-Moser
inequality in $2$ dimensions. Comment. Math. Helv. 67 (1992), no.
3, 471--497.

\bibitem{han} Han Z.C., Asymptotic approach to singular
solutions for nonlinear elliptic equations involving critical
Sobolev exponent. Ann. Inst. H. Poincaré Anal. Non Linéaire 8
(1991), no. 2, 159--174.


\bibitem{hw} Hartman P., Wintner A., On the local behavior of
solutions of non-parabolic Partial differential equations,  Am J.
Math, 75,(1953), 449-476.


\bibitem{hei} Heinz E., On the nonexistence of a surface of
constant mean curvature with finite area and prescribed
rectifiable boundary. Arch. Rational Mech. Anal. 35, (1969),
249--252.

\bibitem{hil} Hildebrandt S., On the Plateau problem for surfaces of
constant mean curvature. Comm. Pure Appl. Math. 23, (1970),
97--114.


\bibitem{hob} Hobson E.W., The theory of spherical and ellipsoidal harmonics.
Chelsea Publishing Company, New York, (1955). xi+500 pp.


\bibitem{isc} Isobe T., Classification of blow-up points and
multiplicity of solutions for $H$-systems. Comm. Partial
Differential Equations 25 (2000), no. 7-8, 1259--1325.

\bibitem{isa1} Isobe T., On the asymptotic analysis of
$H$-systems. I. Asymptotic behavior of large solutions. Adv.
Differential Equations 6 (2001), no. 5, 513--546.


\bibitem{isa2} Isobe T., On the asymptotic analysis of
$H$-systems. II. The construction of large solutions. Adv.
Differential Equations 6 (2001), no. 6, 641--700.




\bibitem{pas} Passaseo D., The effect of the domain shape on
the existence of positive solutions of the equation
$\Delta u+u\sp {2\sp *-1}=0$. Topol. Methods Nonlinear Anal. 3
(1994), no. 1, 27--54.

\bibitem{r0} Rey O., The role of the Green's function in a nonlinear
elliptic equation involving the critical Sobolev exponent. J.
Funct. Anal. 89 (1990), no. 1, 1--52.

\bibitem{r1} Rey O., Concentration of solutions to elliptic
equations with critical nonlinearity. Ann. I.H.P. Anal. Non
Lineaire 9 (1992), no. 2, 201--218.

\bibitem{s} Sasahara Y., An asymptotic analysis for large solutions
of $H$-systems. Adv. Math. Sci. Appl. 5 (1995), no. 1, 219--237.


\bibitem{st0} Steffen K., On the nonuniqueness of surfaces with constant
mean curvature spanning a given contour. Arch. Ration. Mech. Anal.
94, 101-122 (1986).


\bibitem{st} Steffen K., Parametric surfaces of prescribed mean
curvature. Calculus of variations and geometric evolution problems
(Cetraro, 1996), 211--265, Lecture Notes in Math., 1713,
Springer, Berlin, 1999.


\bibitem{str} Struwe M., Plateau's problem and the calculus of
variations. Mathematical Notes, 35. Princeton University Press,
Princeton, NJ, (1988). x+148 pp.

\bibitem{war} Warschawski S. On the differentiability at the
boundary in conformal mapping, Proc. Amer. Math. Soc., 12, (1961),
614-620.


\bibitem{w1} Wente H., The differential equation $\D x = 2 H x_u
\wedge x_v$ with vanishing boundary values, Proc. Amer. Math. Soc.
50 (1975), 131-137.




\end{thebibliography}
\end{document}